\documentclass[sn-mathphys,Numbered]{sn-jnl}

\usepackage{graphicx}%
\usepackage{multirow}%
\usepackage{amsmath,amssymb,amsfonts}%
\usepackage{amsthm}%
\usepackage{mathrsfs}%
\usepackage[title]{appendix}%
\usepackage{xcolor}%
\usepackage{textcomp}%
\usepackage{manyfoot}%
\usepackage{booktabs}%
\usepackage{algorithm}%
\usepackage{algorithmicx}%
\usepackage{algpseudocode}%
\usepackage{listings}%
\usepackage{amsmath}
\usepackage{amssymb}
\usepackage{amsfonts}
\usepackage{graphicx}
\usepackage{subfig}

\begin{document}

\title[Stabilized Isogeometric Collocation Methods for Hyperbolic Conservation Laws]{Stabilized Isogeometric Collocation Methods for Hyperbolic Conservation Laws}


\author*[1]{\fnm{Ryan M.} \sur{Aronson}}\email{rmaronso@stanford.edu}

\author[2]{\fnm{John A.} \sur{Evans}}\email{john.a.evans@colorado.edu}

\affil*[1]{\orgdiv{Institute for Computational and Mathematical Engineering}, \orgname{Stanford University}, \orgaddress{ \city{Stanford}, \state{CA}, \postcode{94305}}}

\affil[2]{\orgdiv{Ann and H.J. Smead Aerospace Engineering Sciences}, \orgname{University of Colorado Boulder}, \orgaddress{ \city{Boulder}, \state{CO}, \postcode{80303}}}


\abstract{We introduce stabilized spline collocation schemes for the numerical solution of nonlinear, hyperbolic conservation laws. A nonlinear, residual-based viscosity stabilization is combined with a projection stabilization-inspired linear operator to stabilize the scheme in the presence of shocks and prevent the propagation of spurious, small-scale oscillations. Due to the nature of collocation schemes, these methods possess the possibility for greatly reduced computational cost of high-order discretizations. Numerical results for the linear advection, Burgers, Buckley-Leverett, and Euler equations show that the scheme is robust in the presence of shocks while maintaining high-order accuracy on smooth problems.}


\maketitle

\section{Introduction}

While originally introduced as an attempt to close the gap between design and analysis, a number of other interesting properties and advantages of Isogeometric Analysis (IGA) \cite{hughes_IGA, cottrell_IGAbook} have been discovered. The smooth spline basis functions used in IGA have been shown to have improved approximation power when compared to the $C^0$ polynomial basis functions used in standard Finite Element Analysis (FEA) \cite{evans_nwidths, sande2020explicit, bressan2019approximation}. Moreover, smooth basis functions enable the development of collocation methods based on the strong form of the governing partial differential equation (PDE) \cite{reali_intro, auricchio_isogeometric_coll}. Unlike standard Galerkin methods, collocation schemes do not require computations of element integrals, which are typically evaluated using expensive Gauss quadrature schemes. Instead, the strong form of the equations is evaluated at a set of collocation points to determine the numerical solution. The potential of collocation for reducing the cost of simulations has already been demonstrated \cite{schillinger_cost_comp}, and these schemes are ideally suited for explicit dynamics simulations where the cost of computing the residual integrals can dominate the computational time. This is especially true for the simulation of nonlinear PDEs using high-order bases.

In the context of incompressible fluid dynamics, there was interest in using spline basis functions for analysis even before the invention of IGA \cite{botella_collocation, kravchenko1999b}. The improved dispersion properties of spline bases when compared to standard, high-order bases make these types of schemes attractive for simulations of turbulent flows, and collocation schemes can be used to efficiently evaluate the nonlinear convection term. Recently, B-spline collocation has been used in conjunction with spectral methods for simulations of turbulent channel flows \cite{lee2015DNS}, and the IGA community has continued to investigate both Galerkin and collocation discretization schemes based purely on splines. For example, in \cite{aronsonDivConforming} we developed divergence-conforming B-spline collocation schemes and in \cite{aronson2023stabilized} we described the construction of stabilized spline collocation schemes which combat both advection and pressure instabilities. 

Much less common in the literature are spline-based discretization schemes developed for the simulation of compressible flow problems. In this case a nonlinear stabilization is typically required to ensure high-order discretization schemes remain well behaved in the presence of discontinuities, such as shocks. For example, in \cite{jaeschke2020CompressibleIGA, moller2020CompressibleIGA} the Algebraic Flux Correction method was adapted to spline-based Galerkin discretizations of conservation laws. Alternatively, \cite{duvigneau2018isogeometric} used a discontinuous Galerkin IGA method coupled with an artificial viscosity-based shock capturing scheme of the form presented in \cite{persson2006sub} for nonlinear stabilization. As both of these schemes are based on the Galerkin method, however, the integrals of the nonlinear flux terms, which must be computed at every time step, can render the methods quite expensive. 

In this work, we extend the use of spline collocation discretization schemes to hyperbolic conservation laws, as a first step towards modeling compressible flows. Much like the setting of explicit, nonlinear dynamics, the use of collocation schemes for conservation laws has the potential to greatly the improve computational complexity of high-order methods, both due to the lack of required numerical integrations as well as the tensor product structure of the mass matrix that arises from collocation. We stabilize the collocation schemes in the presence of shocks using an artificial viscosity, like in \cite{duvigneau2018isogeometric}. Unlike that work, however, this viscosity is constructed proportionally to the residual of the governing PDE, providing stabilization in regions near shocks while adding little viscosity in smooth regions where the solution is well resolved. This is inspired by the work in \cite{guermond2011entropy, stiernstrom2021SBP, tominec2023RBF, nazarov2017investigation}, where similar ideas have found success in creating high-order, stable schemes based on a variety of discretizations. We also introduce a novel linear stabilization technique inspired by projection stabilization \cite{braack2006lps}, which effectively suppresses small-scale oscillations that can sometimes appear upwind of the shock, even when the artificial viscosity is included. This stabilization is essentially an alternative to the Streamline-Upwind-Petrov-Galerkin inspired stabilization developed in \cite{aronson2023stabilized}, but with the advantage of having no unsteady term in the stabilization which could make explicit time stepping more complex.

The rest of the work is laid out as follows: First we describe the details of the basic discretization, starting with the strong forms of the PDEs we seek to solve. In this work we examine multiple scalar conservation laws but only one system of conservation laws, namely the compressible Euler equations. We then discuss the details of the spatial discretization process, including the construction of B-spline basis functions and appropriate collocation points. The remainder of Section 2 discusses explicit time integration and the speedups that are possible from the choice of collocation in space. Section 3 then discusses the stabilization schemes, starting with the nonlinear artificial viscosity construction, followed by our novel linear stabilization technique. Finally, Section 4 summarizes the results obtained for a variety of canonical test problems and shows the promise of the method for high-order Computational Fluid Dynamics (CFD) applications.

\section{Isogeometric Collocation Formulation of Conservation Laws}

To begin, we state the strong formulations of the partial differential equations we are interested in solving in this work. A full problem statement for an arbitrary scalar conservation law posed on the domain $\Omega \subset \mathbb{R}^d$ can be stated as

\bigskip

$$
\left\{ \hspace{5pt}
\parbox{5in}{
\noindent Given flux function $\mathbf{f} : \mathbb{R} \rightarrow \mathbb{R}^d$, boundary operator and data $\mathcal{G}$ and $g$, initial condition $\phi_0 : \Omega \rightarrow \mathbb{R}$, and final time $t_f \in \mathbb{R}$, find $\phi : \Omega \times [0, t_f] \rightarrow \mathbb{R}$ such that:
\begin{equation}
   \frac{\partial \phi}{\partial t}(\mathbf{x}, t) + \nabla \cdot \mathbf{f}(\phi(\mathbf{x}, t))  = 0  \textup{ in } \Omega,
\end{equation}
\begin{equation}
    \mathcal{G}\phi(\mathbf{x}, t) = g(\mathbf{x}, t) \textup{ on }  \partial \Omega,
\end{equation}
\begin{equation}
    \phi(\mathbf{x}, 0) = \phi_0(\mathbf{x}) \textup{ in } \Omega.
\end{equation}
}
\right.
$$

\bigskip

\noindent The first equation represents the conservation law itself, while the second and third specify the boundary and initial conditions, respectively. The specification of the flux function defines the conservation law being solved, and the conserved quantity $\phi$ can take on a different physical meaning in each case. In addition to scalar conservation laws, we also consider the compressible Euler equations, which can be expressed in conservative form as

\bigskip

$$
\left\{ \hspace{5pt}
\parbox{5in}{
\noindent Given boundary operator and data $\mathcal{G}$ and $\mathbf{g}$, initial conditions $\rho_0 : \Omega \rightarrow \mathbb{R}$, $\rho \mathbf{u}_0 : \Omega \rightarrow \mathbb{R}^d$, and $E_0 : \Omega \rightarrow \mathbb{R}$, and final time $t_f \in \mathbb{R}$, find $\rho : \Omega \times [0, t_f] \rightarrow \mathbb{R}$, $\rho \mathbf{u} : \Omega \times [0, t_f] \rightarrow \mathbb{R}^d$, and $E : \Omega \times [0, t_f] \rightarrow \mathbb{R}$ such that:
\begin{equation}
    \frac{\partial}{\partial t} 
    \begin{bmatrix} \rho \\ \rho \mathbf{u}\\ E \end{bmatrix}(\mathbf{x}, t)
    + \nabla \cdot
    \begin{bmatrix} \rho \mathbf{u} \\ \rho \mathbf{u} \otimes \mathbf{u} + p \mathbf{I} \\ \mathbf{u} (E + p) \end{bmatrix}(\mathbf{x},t)
    = \mathbf{0} \textup{ in } \Omega,
\end{equation}
\begin{equation}
    \mathcal{G}\begin{bmatrix} \rho \\ \rho \mathbf{u}\\ E \end{bmatrix}(\mathbf{x}, t) = \mathbf{g}(\mathbf{x}, t) \textup{ on }  \partial \Omega,
\end{equation}
\begin{equation}
    \begin{bmatrix} \rho \\ \rho \mathbf{u}\\ E \end{bmatrix}(\mathbf{x}, 0) = \begin{bmatrix} \rho_0 \\ \rho \mathbf{u}_0\\ E_0 \end{bmatrix}(\mathbf{x}) \textup{ in } \Omega,
\end{equation}
\begin{equation}
    p(\mathbf{x},t) = (\gamma - 1) (E(\mathbf{x},t) - \frac{1}{2} \rho \mathbf{u}(\mathbf{x},t) \cdot \mathbf{u}(\mathbf{x},t)),
\end{equation}
\begin{equation}
    T(\mathbf{x},t) = \frac{p(\mathbf{x},t)}{\rho(\mathbf{x},t)}.
\end{equation}
}
\right.
$$

\bigskip

\noindent Here $\rho$ is the density of the fluid, $\mathbf{u}$  is the fluid velocity, $E$ is the total energy, and $p$ is the pressure. The first three equations represent the conservation of mass, momentum, and energy, respectively. The system of equations is closed by assuming the fluid is an ideal gas with constant heat capacity, which yields the last 2 constitutive laws.

In order to discretize the problems above in space with a collocation scheme, we require a set of proper basis functions $\{N_i\}_{i=1}^{n_d}$ and collocation points $\{\mathbf{x}_i\}_{i=1}^{n_d}$, where $n_d$ is the total number of degrees of freedom within the discrete space. A system of ordinary differential equations in time is formed by assuming the discrete solution is represented by a linear combination of these basis functions at any instant in time, with coefficients determined by requiring the strong form of the PDE to hold at every collocation point. We discuss these steps in further detail below.

\subsection{B-Spline Basis Functions}

For a collocation scheme to be well-defined, the basis functions $\{N_i\}_{i=1}^{n_d}$ must be sufficiently smooth so that the strong form residual can be evaluated. Note that the standard $C^0$ piecewise polynomial basis functions used in FEA do not satisfy this requirement. The high-order spline basis functions used in IGA, however, can easily be constructed with sufficient regularity. In this work we will only utilize B-spline basis functions, but the methods trivially extend to NURBS basis functions as well. 

In 1D, the specification of the polynomial order $k$ and the knot vector $\Xi = \{\xi_1, ... \xi_{n+k+1}\}$ fully describes the B-spline basis functions. These functions are evaluated through the Cox-de Boor recursion, starting with the $k = 0$ basis functions

\begin{equation}
    N_{i,0}(\xi) = \begin{cases} 
                1 & \xi_i \leq \xi \leq \xi_{i+1} \\
                0 & \text{otherwise},
              \end{cases}
\end{equation}

\noindent and higher-order bases defined via

\begin{equation}
    N_{i,k}(\xi) = \frac{\xi - \xi_i}{\xi_{i+k} - \xi_i}N_{i,k-1}(\xi) + \frac{\xi_{i+k+1} - \xi}{\xi_{i+k+1} - \xi_{i+1}}N_{i+1,k-1}(\xi).
\end{equation}

\noindent At any knot in the knot vector $\Xi$, the basis functions are $C^{k-\ell}$-continuous, where $\ell$ is the multiplicity of the knot. Thus it is straightforward to define basis functions with suitable regularity for collocation. In this work we will always utilize maximal continuity splines, thus our knot vectors never contain repeated interior knots. We also utilize open knot vectors, meaning that the first and last knot values are repeated $k+1$ times, making the spline interpolatory at these locations. 

B-spline basis functions in higher spatial dimensions are constructed by simply taking the tensor product of one dimensional B-spline bases in each direction \cite{cottrell_IGAbook}. In general, the degrees and knot vectors defining the basis in each parametric direction can differ, but for simplicity, most of the results presented in this work only consider tensor product bases with the same 1D basis in each direction. For notational convenience we denote a spline space defined on a domain $\Omega \subset \mathbb{R}^d$ as $S^h = \textup{span}\{ N_i\}_{i = 1}^{n_d}$, where the functions $N_i : \Omega \rightarrow \mathbb{R}$ are B-spline basis functions and $n_d$ is the total number of functions in the discrete space.

\subsection{The Greville Abcissae}

With the choice of basis functions made, we must now select a set of collocation points to obtain a well-posed scheme. By well-posed, we mean that we desire a set of collocation points which produce a non-singular interpolation operator. There are several possible choices for these points in the isogeometric collocation setting, such as the Greville abscissae \cite{auricchio_isogeometric_coll, johnson_Greville_coll}, the Demko abscissae \cite{demko_abs}, or the Cauchy-Galerkin points \cite{montardini_optimal, anitescu2015isogeometric}. The Greville abscissae are the simplest choice, as they are computed in 1D simply through averages of knot values via

\begin{equation}
    \hat{\xi_i} = \frac{\xi_{i+1} + ... + \xi_{i+k}}{k}.
\end{equation}

\noindent In higher dimensions, the Greville points of the tensor product space are simply the tensor product of the 1D Greville abscissae. In this work we only consider collocation at Greville points, though we note that there are specific cases with large polynomial degrees of 20 or higher and certain non-uniform meshes where the scheme will no longer be well-posed \cite{jia1988spline}. The Demko abscissae provably produce a non-singular interpolation operator for all cases, but require an iterative approximation to compute. Finally, we note that the main advantage of the use of the Cauchy-Galerkin points is recovery of optimal $L^2$ convergence rates for odd basis degrees. However, in this work we only consider first-order conservation laws, and \cite{aronsonDivConforming, aronson2023stabilized} have shown that collocating first-order PDEs  results in the same optimal asymptotic convergence rates, even when the Greville points are utilized.

\subsection{Spatial Discretization}

With a complete description of the elements required for spline collocation schemes, we move on to apply this technique for the spatial discretization of conservation laws. We consider the specific case of scalar conservation laws and assume fully Dirichlet boundary conditions for simplicity. Then the semi-discrete problem is obtained by assuming the numerical solution lies in a predefined spline space $S^h$ with Greville abscissae $\{\mathbf{x}_i\}_{i = 1}^{n_d}$, and the strong form of the problem holds at every collocation point. Defining the space $S^h_T = \{ \phi^h : \Omega \times [0, t_f] \rightarrow \mathbb{R} : \phi^h(\cdot, t) \in S^h \textit{ }  \forall t \in (0,t_f) \}$, the semi-discrete form can be written as

\bigskip

$$
\left\{ \hspace{5pt}
\parbox{5in}{
\noindent Given $\mathbf{f} : \mathbb{R} \rightarrow \mathbb{R}^d$, $g: \partial \Omega \times [0,t_f] \rightarrow \mathbb{R}$, initial condition $\phi_0 : \Omega \rightarrow \mathbb{R}$, and final time $t_f \in \mathbb{R}$, find $\phi^h \in S^h_T$ such that:
\begin{equation}
   \frac{\partial \phi^h}{\partial t}(\mathbf{x}_i, t) + \nabla \cdot \mathbf{f} (\phi^h(\mathbf{x}_i, t)) = 0 \quad \forall \mathbf{x}_i \in \Omega,
\end{equation}
\begin{equation}
    \phi^h(\mathbf{x}_i, t) = g(\mathbf{x}_i, t) \quad \forall \mathbf{x}_i \in \partial \Omega.
\end{equation}
\begin{equation}
    \phi^h(\mathbf{x}_i, 0) = \phi_0(\mathbf{x}_i) \quad \forall \mathbf{x}_i \in \Omega
\end{equation}
}
\right.
$$

\bigskip

This semi-discrete form is simply a system of ordinary differential equations, with one equation per collocation point. Integration in time can then be done with any standard scheme, and is discussed further below. For completeness we also include the semi-discrete form of the Euler equations subjected to Dirichlet boundary conditions. Denoting by $\mathbf{S}^h_T$ the vector-valued space wherein each component belongs to $S^h_T$, we arrive at

\bigskip

$$
\left\{ \hspace{5pt}
\parbox{5in}{
\noindent Given $\mathbf{g} : \partial \Omega \times [0,t_f] \rightarrow \mathbb{R}^{d+2}$, initial conditions $\rho_0 : \Omega \rightarrow \mathbb{R}$, $\rho \mathbf{u}_0 : \Omega \rightarrow \mathbb{R}^d$, and $E_0 : \Omega \rightarrow \mathbb{R}$, and final time $t_f \in \mathbb{R}$, find $\rho^h \in S^h_T$, $\rho \mathbf{u}^h \in \mathbf{S}^h_T$, and $E^h \in S^h_T$ such that:
\begin{equation}
    \frac{\partial}{\partial t} 
    \begin{bmatrix} \rho^h \\ \rho \mathbf{u}^h\\ E^h \end{bmatrix}(\mathbf{x}_i, t)
    + \nabla \cdot
    \begin{bmatrix} \rho \mathbf{u}^h \\ \rho \mathbf{u}^h \otimes \mathbf{u}^h + p^h \mathbf{I} \\ \mathbf{u}^h (E^h + p^h) \end{bmatrix}(\mathbf{x}_i,t)
    = \mathbf{0} \quad \forall \mathbf{x}_i\in \Omega,
\end{equation}
\begin{equation}
    \begin{bmatrix} \rho^h \\ \rho \mathbf{u}^h\\ E^h \end{bmatrix}(\mathbf{x}_i, t) = \mathbf{g}(\mathbf{x}_i, t) \quad \forall \mathbf{x}_i \in \partial \Omega,
\end{equation}
\begin{equation}
    \begin{bmatrix} \rho^h \\ \rho \mathbf{u}^h\\ E^h \end{bmatrix}(\mathbf{x}_i, 0) = \begin{bmatrix} \rho_0 \\ \rho \mathbf{u}_0\\ E_0 \end{bmatrix}(\mathbf{x}_i) \quad \forall \mathbf{x}_i \in \Omega,
\end{equation}
\begin{equation}
    p^h(\mathbf{x}_i,t) = (\gamma - 1) (E^h(\mathbf{x}_i,t) - \frac{1}{2} \rho \mathbf{u}^h(\mathbf{x}_i,t) \cdot \mathbf{u}^h(\mathbf{x}_i,t)),
\end{equation}
\begin{equation}
    T^h(\mathbf{x}_i,t) = \frac{p^h(\mathbf{x}_i,t)}{\rho^h(\mathbf{x}_i,t)}.
\end{equation}
}
\right.
$$

\bigskip

The enforcement of initial conditions deserves slightly more attention in the case of discontinuous initial data. Simply collocating this data as stated above will yield an oscillatory solution due to Gibb's phenomena near the discontinuities. However, we have found in many of our studies that the inclusion of the nonlinear stabilization described in the following section effectively suppressed these initial oscillations in relatively few time steps, seemingly without adversely affecting the solution quality at later times. One can also collocate a smoothed initial condition to remove these oscillations, or use the initial data to set the control coefficients of the spline solution directly, bypassing the need for interpolation altogether. This is also how discontinuous boundary conditions were enforced in \cite{aronson2023stabilized}.

\subsection{Time Integration - Residual Computation}

In this work we are interested in explicit time integration schemes, such as RK4, which do not require a nonlinear solution procedure at every time step in the presence of nonlinear flux functions. The first main step of explicit time integration is the construction of the spatial residual vector, and this process can be quite costly. In Galerkin methods, this involves the computation of terms of the form

\begin{equation}
    \int_{\Omega_E} (\nabla \cdot \mathbf{f}) \psi^h d\Omega, 
\end{equation}

\noindent over each element $\Omega_E$, for all $\psi^h$ in the discrete test function space. As these integrals are usually approximated using Gauss quadrature rules, the integral evaluation cost scales like $O(k+1)^{2d}$ per element.

A major cost advantage of collocation methods is the replacement of these costly integrations with strong form evaluations. For a standard spline collocation scheme the number of operations per element instead scales like $O(k+1)^d$ per element, as most elements include only $O(1)$ collocation points. 

Inspired by finite difference schemes for conservation laws, in which the conservative form of the PDE is discretized to capture correct shock speeds, for conservation laws we also choose to collocate the conservative form. Thus we require a method to evaluate the divergence of the flux. For our method we use the flux evaluated on the collocation points to fit a spline surface in the same discrete space as the unknowns, allowing for the evaluation of partial derivatives. This process involves the application of the inverse of the mass matrix, which can be done efficiently in the same asymptotic complexity as the residual evaluations, as described below. 

\subsection{Time Integration - The Mass Matrix}

Like a standard Galerkin finite element scheme, isogeometric collocation schemes will lead to non-diagonal mass matrices. Thus linear solves are required at every time step, not only to fit a spline surface to the flux values but also to obtain the updates to the solution coefficients in time. 

While this can limit the efficiency of timestepping, collocation schemes possess convenient properties to help mitigate this cost. As the mass matrix is not created via spatial integration, the mass matrix in two or three dimensional cases is exactly the tensor product of the 1D mass matrices in each parametric direction, even when the geometric mapping in not affine. This is not the case in Galerkin schemes \cite{chan2018multi}. In particular, this also means that the inverse of the mass matrix is equal to the tensor product of the 1D inverses. Thus by storing a factorization of the 1D mass matrices, we can efficiently perform the linear solves required during explicit time stepping. 

In particular, the application of the full inverse can be written as $O(n^{d-1})$ applications of the 1D LU factorization, where $n$ is the number of degrees of freedom in each direction and $d$ is the number of spatial dimensions. The LU decomposition matrices have the same bandwidth as the original mass matrices \cite{kilicc2013inverse}, meaning that for the 1D matrices each of the factors has bandwidth $k+1$. Using this, we can estimate that the cost of a 1D mass matrix solve is $O(n(k+1))$ once we have the factorization. Therefore, the cost of application of the full mass inverse requires $O(n^d(k+1))$ operations. For completeness we show the specific forms of these operations in the 2D and 3D settings below. 

For a 3D problem, we let $M^{-1} = M_x^{-1} \otimes M_y^{-1} \otimes M_z^{-1}$ be the inverse of the collocation mass matrix. If we assume for simplicity that the number of basis functions is the same in each direction and equal to $n$, then $M_x^{-1} \in \mathbb{R}^{n \times n}$, $M_y^{-1} \in \mathbb{R}^{n \times n}$, and $M_z^{-1} \in \mathbb{R}^{n \times n}$ are the corresponding 1D mass matrix inverses. Algorithm \ref{alg:mass_inv} details the process of applying the full mass matrix inverse to a vector $\mathbf{x}$ using the 1D inverses. Note that for 2D problems this process simplifies to reshaping $\mathbf{x}$ into a matrix $X$, computing $M_x^{-1} X M_y^{-T}$, and reshaping back to a vector.

\begin{algorithm}
\caption{Mass Matrix Solution}
\label{alg:mass_inv}
\begin{algorithmic}[1]
        \State $X \gets reshape(\mathbf{x}, [n,n,n])$
        \State $u_{mjl} \gets M_{z, mn}^{-1} x_{jln}$
        \State $v_{kmj} \gets M_{y, kl}^{-1} u_{mjl}$
        \State $w_{ikm} \gets M_{x, ij}^{-1} v_{kmj}$
        \State $\mathbf{w} = reshape(W, [n^3,1])$
        \State \Return $\mathbf{w}$
 \end{algorithmic}
\end{algorithm}

This is faster than utilizing a naïve LU factorization of a full, 2D or 3D Galerkin mass matrix at every time step, which would cost approximately $O(n^{d+1} (k+1)^{d-1} )$ due to increased bandwidth of the matrix ($n^d$ rows with bandwidth of $n(k+1)^{d-1}$). Note that this cost can also be improved in the Galerkin setting through the use of sum factorization \cite{bressan2019sum}, though this is much more complex for non-affine domains than the collocation scheme. Also worth mentioning is the Spectral Element Method (SEM) \cite{patera1984spectral, sherwin1995triangular}, which is obtained by using a Gauss-Legendre-Lobatto numerical integration scheme and a Lagrange polynomial basis with nodes located at the quadrature points. This too can be written as a collocation scheme with the added advantage of resulting in a diagonal mass matrix. In some ways the spline collocation technique presented in this work can be thought of as the generalization of SEM to smooth, spline bases.  

Finally, one can generate schemes that do not require the solution of a mass matrix system at every time step, even for non-diagonal mass matrices, using lumping techniques \cite{voet2023mathematical}. These techniques have also been applied in isogeometric collocation schemes and result in efficient schemes that can be comparable to finite-difference-time-domain methods in terms of cost \cite{evans2018explicit}, though we do not further pursue this in this work.  

Combining the results presented in this section with the previous, the total cost of the collocation schemes approximately scales like $O(n^d(k+1)^d)$ per time step, with the asymptotic scaling being dominated by the cost of evaluation of the flux and its derivatives on the collocation points dominates the scaling, not the applications of inverse mass matrices.

\subsection{A Note on Maximum Stable Time Step Sizes}

In the previous sections we have focused on the operations and cost required for one explicit time step. Perhaps equally important, however, is the largest stable time step that can be taken using the scheme. Common existing high order methods using $C^0$ basis functions, like SEM and discontinuous Galerkin (DG), have critical time steps which scale like $h / k^2$, where $h$ is the mesh size and $k$ is the polynomial order \cite{warburton2008taming}. Spline-based discretizations, on the other hand, have been shown to have critical time steps that scale as $h/k$ \cite{chan2018multi, zampieri2021isogeometric}. Thus spline discretizations also have potential for improved efficiency through a factor of $k$ larger time steps when compared with standard, high-order schemes. 

\section{Stabilization}

For problems with smooth solutions, the collocation scheme defined above is sufficient to obtain a high-order numerical approximation. However, the presence of shocks or other sharp features, which are ubiquitous in solutions of nonlinear conservation laws, will cause large oscillations and blow-up in the numerical solution. In this section we describe the stabilization mechanisms we have applied to the collocation schemes, starting with a nonlinear shock-capturing scheme using a residual-based viscosity, followed by a linear stabilization term which smooths small-scale, high frequency oscillations, similar to upwinding in the finite difference context.  

\subsection{Residual-Based Viscosity}

To stabilize our discretizations in the presence of shocks and other discontinuities, we adapt a residual-based viscosity stabilization similar to those presented in \cite{guermond2011entropy, stiernstrom2021SBP, tominec2023RBF, nazarov2017investigation} to the spline collocation setting. We start by considering an arbitrary scalar conservation law before moving on to the case of the compressible Euler equations. 

\subsubsection{Scalar Conservation Laws}

In the scalar conservation law setting, the nonlinear stabilization is added by including the Laplacian of the conserved quantity, so the statement of the semi-discretized collocation scheme becomes

\bigskip

$$
\left\{ \hspace{5pt}
\parbox{5in}{
\noindent Given $\mathbf{f} : \mathbb{R} \rightarrow \mathbb{R}^d$, $g: \partial \Omega \times [0,t_f] \rightarrow \mathbb{R}$, initial condition $\phi_0 : \Omega \rightarrow \mathbb{R}$, and final time $t_f \in \mathbb{R}$, find $\phi^h \in S^h_T$ such that:
\begin{equation}
   \frac{\partial \phi^h}{\partial t}(\mathbf{x}_i,t) + \nabla \cdot \mathbf{f} (\phi^h(\mathbf{x}_i,t)) = \nu_{art}(\mathbf{x}_i, t) \Delta \phi^h (\mathbf{x}_i, t) \quad \forall \mathbf{x}_i \in \Omega,
\end{equation}
\begin{equation}
    \phi^h(\mathbf{x}_i, t) = g(\mathbf{x}_i, t) \quad \forall \mathbf{x}_i \in \partial \Omega,
\end{equation}
\begin{equation}
    \phi^h(\mathbf{x}_i, 0) = \phi_0(\mathbf{x}_i) \quad \forall \mathbf{x}_i \in \Omega.
\end{equation}
}
\right.
$$

\bigskip

\noindent Note that we have omitted terms involving derivatives of the viscosity coefficient $ \nu_{art}(\mathbf{x}_i, t)$. We see this as analogous to choosing a viscosity which is piecewise constant in each cell, commonly assumed in other high-order shock-capturing schemes \cite{persson2006sub}. It is perhaps true that the consistency of the scheme would be improved with the inclusion of the viscosity coefficient within the viscous flux, meaning that terms involving the derivatives of $\nu_{art}$ would also be included. However, this increases the complexity of the method, as one would have to include a step such as interpolation of the $\nu_{art}$ values, as was done in \cite{aronson2023stabilized}. Moreover, we have found that the treatment of viscosity presented here has led to quality results, suggesting that these derivative terms are perhaps not necessary.  

The value of the viscosity coefficient is determined via the residual of the PDE, as it is expected that the residual will be large near the sharp features, such as shocks, which cause instability. For a scalar conservation law, we define the discrete residual as 

\begin{equation}
    R^h(\mathbf{x}, t) = \frac{\partial \phi^h}{\partial t}(\mathbf{x}, t) + \nabla \cdot \mathbf{f} (\phi^h(\mathbf{x},t)).
\end{equation}

\noindent To approximate the unsteady term in the residual, we store values of the solution at previous time steps and utilize an explicit, high-order Backward Difference Formula (BDF). We have found using the 4th-order formula sufficient for the results in this work, though in general one may want to adjust the order based on the chosen spatial and time discretization schemes. The BDF order determines how many past solutions must be stored, as well as how many initial time steps must be taken in order to generate enough history data to begin using the desired scheme. At the beginning of the simulation one can either wait to compute this residual (and thus add nonlinear stabilization) until enough history is available for the desired BDF scheme, or one can gradually increase the BDF order as data becomes available until the final order is reached. We have found that the results are not highly dependent on this choice, though the convergence rates on some smooth problems with small spatial discretization errors can be arrested by the use of low-order BDF schemes at early time steps.

\begin{figure}
\centering
\subfloat[1D grid]{\label{sfig:1D_grid}\includegraphics[width=.4\textwidth]{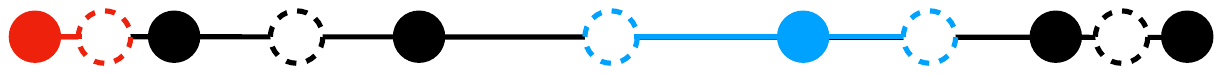}}\hfill
\subfloat[2D grid]{\label{sfig:2D_grid}\includegraphics[width=.4\textwidth]{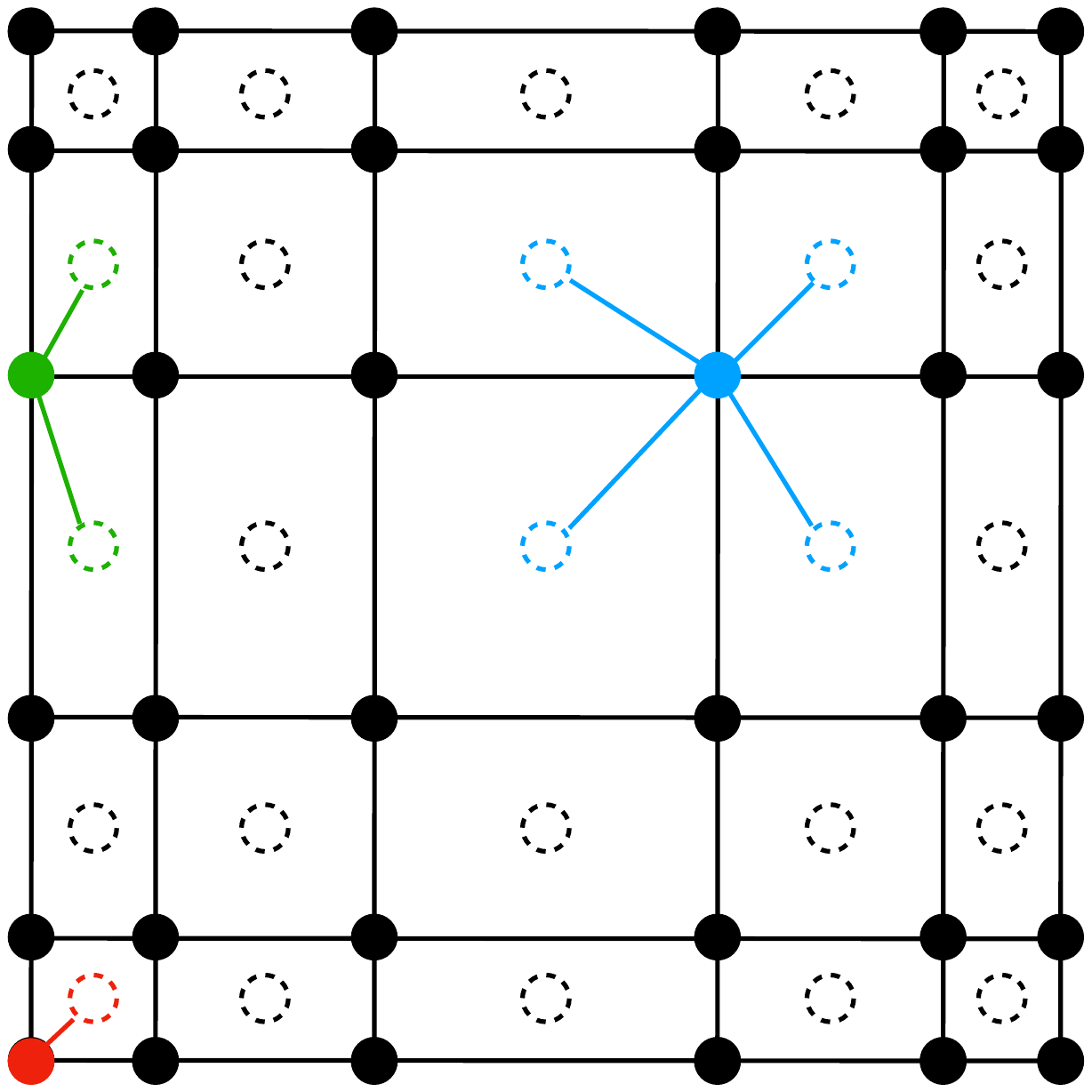}} \\
\caption{Schematic showing residual sampling strategies in 1D and 2D. Solid dots represent collocation points while dashed dots represent the possible locations at which the PDE residual is sampled. Shown in color are examples of the sampled values used to determine the viscosity at interior, edge, and corner collocation points.}
\label{fig:sampling_grids}
\end{figure}

Note that for a a steady problem, the collocation solution is defined such that the discrete residual is zero at the collocation points. Thus we choose a slightly modified residual to drive our viscosity at each collocation point. In particular, we define the residual 

\begin{equation}
    \Tilde{R}^h(\mathbf{x}_i, t) = \max_{\Tilde{\mathbf{x}}\in M}(|R^h(\Tilde{\mathbf{x}}, t)|),
    \label{eq:Rtilde}
\end{equation}

\noindent where $\mathbf{x}_i$ is a collocation point and $M$ is the set of points (which are not collocation points) at which the residual is sampled for each collocation point. These possible points are selected as the centroids of each Greville cell, and for each collocation point the set $M$ is defined as the sampling points closest to the collocation point as shown in Figure \ref{fig:sampling_grids}.

This residual defines a viscosity at each collocation point via

\begin{equation}
    \nu_{RB}(\mathbf{x}_i, t) = \frac{C_{RB} h(\mathbf{x}_i)^2 \Tilde{R}^h(\mathbf{x}_i, t)} {m(\mathbf{x}_i, t)},
    \label{eq:rb_visc}
\end{equation}

\noindent where $C_{RB}$ is a tunable constant, $h$ is the mesh size defined by the distance between collocation points, and $m$ is the normalization factor \cite{stiernstrom2021SBP}

\begin{equation}
    m(\mathbf{x}_i, t) = \Tilde{\phi} - || \phi(\mathbf{x}_i, t) - \Bar{\phi}||_{\infty}, 
\end{equation}

\noindent with $\Tilde{\phi}$ defined as the maximum value of the solution minus the minimum value of the solution over all of the collocation points directly adjacent to $\mathbf{x}_i$ at the current time, including the point itself. 

At each collocation point we also construct a first order viscosity of the form

\begin{equation}
    \nu_{FO}(\mathbf{x}_i, t) = C_{max} h(\mathbf{x}_i) c(\mathbf{x}_i, t),
\end{equation}

\noindent where $C_{max}$ is another tunable parameter and $c$ is an estimate of the wavespeed at each collocation point $|\mathbf{f}' (\phi(\mathbf{x}_i, t))|$. For each collocation point this is computed as the maximum wavespeed over the grid of width 9 collocation points in each direction centered at the current collocation point. 

The final stabilization viscosity is constructed by limiting the residual-based viscosity with the first order viscosity

\begin{equation}
    \nu_{art}(\mathbf{x}_i, t) = \min(\nu_{RB}(\mathbf{x}_i, t), \nu_{FO}(\mathbf{x}_i, t)).
    \label{eq:final_visc}
\end{equation}

Up to this point we have left all of the stabilization expressions continuous in time. Within the RK time integration schemes, we simply compute a stabilization viscosity at every collocation point using the solution information at the beginning of the time step, then use this same viscosity for all of the RK substeps. 

\subsubsection{The Euler Equations}

In the case of the Euler equations, there is slightly more freedom in the design of stabilization schemes. One option is to design stabilization terms which are directly proportional to the Laplacian of the conserved quantities, like in the scalar case above. Indeed, this is the approach used in \cite{tominec2023RBF}, and we have found that this does indeed sufficiently stabilize spline collocation solutions. In this case, the semi-discrete form can be written as

\bigskip

$$
\left\{ \hspace{5pt}
\parbox{5in}{
\noindent Given $\mathbf{g} : \partial \Omega \times [0,t_f] \rightarrow \mathbb{R}^{d+2}$, initial conditions $\rho_0 : \Omega \rightarrow \mathbb{R}$, $\rho \mathbf{u}_0 : \Omega \rightarrow \mathbb{R}^d$, and $E_0 : \Omega \rightarrow \mathbb{R}$, and final time $t_f \in \mathbb{R}$, find $\rho^h \in S^h_T$, $\rho \mathbf{u}^h \in \mathbf{S}^h_T$, and $E^h \in S^h_T$ such that:
\begin{equation}
\begin{split}
    &\frac{\partial}{\partial t} 
    \begin{bmatrix} \rho^h \\ \rho \mathbf{u}^h\\ E^h \end{bmatrix}(\mathbf{x}_i, t)
    + \nabla \cdot
    \begin{bmatrix} \rho \mathbf{u}^h \\ \rho \mathbf{u}^h \otimes \mathbf{u}^h + p^h \mathbf{I} \\ \mathbf{u}^h (E^h + p^h) \end{bmatrix}(\mathbf{x}_i,t)
    \\ &=  \nu_{art}(\mathbf{x}_i, t) \nabla \cdot
    \begin{bmatrix} \nabla \rho \\ \nabla \rho \mathbf{u} \\ \nabla E \end{bmatrix}(\mathbf{x}_i, t) \quad \forall \mathbf{x}_i \in \Omega,
    \end{split}
\end{equation}
\begin{equation}
    \begin{bmatrix} \rho^h \\ \rho \mathbf{u}^h\\ E^h \end{bmatrix}(\mathbf{x}_i, t) = \mathbf{g}(\mathbf{x}_i, t) \quad \forall \mathbf{x}_i \in \partial \Omega,
\end{equation}
\begin{equation}
    \begin{bmatrix} \rho^h \\ \rho \mathbf{u}^h\\ E^h \end{bmatrix}(\mathbf{x}_i, 0) = \begin{bmatrix} \rho_0 \\ \rho \mathbf{u}_0\\ E_0 \end{bmatrix}(\mathbf{x}_i) \quad \forall \mathbf{x}_i \in \Omega,
\end{equation}
\begin{equation}
    p^h(\mathbf{x}_i,t) = (\gamma - 1) (E^h(\mathbf{x}_i,t) - \frac{1}{2} \rho \mathbf{u}^h(\mathbf{x}_i,t) \cdot \mathbf{u}^h(\mathbf{x}_i,t)),
\end{equation}
\begin{equation}
    T^h(\mathbf{x}_i,t) = \frac{p^h(\mathbf{x}_i,t)}{\rho^h(\mathbf{x}_i,t)}.
\end{equation}
}
\right.
$$

\bigskip

\noindent For each of the conservation of mass, momentum, and energy, we define a residual at each collocation point in the same way as for a scalar conservation law, shown in Equation \eqref{eq:Rtilde}, and a viscosity based on each is constructed according to Equation \eqref{eq:rb_visc}. Then, the value of $\nu_{art}$ at each collocation point is chosen to be the maximum of each of these viscosities, limited by the first-order viscosity. We shall refer to this method as stabilization via Laplacian fluxes.

However, as mentioned in \cite{guermond2011entropy}, this approach is not necessarily consistent with the entropy inequality. Instead, we can also choose to regularize the system using the Guermond-Popov fluxes proposed in \cite{guermond2014regularization} and used in conjunction with a  residual-based viscosity in \cite{nazarov2017investigation}. Let $\kappa_{art}$ and $\mu_{art}$ be artificial viscosities to be determined, then the regularized system becomes


\bigskip

$$
\left\{ \hspace{5pt}
\parbox{5in}{
\noindent Given $\mathbf{g} : \partial \Omega \times [0,t_f] \rightarrow \mathbb{R}^{d+2}$, initial conditions $\rho_0 : \Omega \rightarrow \mathbb{R}$, $\rho \mathbf{u}_0 : \Omega \rightarrow \mathbb{R}^d$, and $E_0 : \Omega \rightarrow \mathbb{R}$, and final time $t_f \in \mathbb{R}$, find $\rho^h \in S^h_T$, $\rho \mathbf{u}^h \in \mathbf{S}^h_T$, and $E^h \in S^h_T$ such that:
\begin{equation}
\begin{split}
    &\frac{\partial}{\partial t} 
    \begin{bmatrix} \rho^h \\ \rho \mathbf{u}^h\\ E^h \end{bmatrix}(\mathbf{x}_i, t)
    + \nabla \cdot
    \begin{bmatrix} \rho \mathbf{u}^h \\ \rho \mathbf{u}^h \otimes \mathbf{u}^h + p^h \mathbf{I} \\ \mathbf{u}^h (E^h + p^h) \end{bmatrix}(\mathbf{x}_i,t)
    \\ &=  \nabla \cdot
    \begin{bmatrix} \kappa_{art} \nabla \rho \\ \mu_{art} \rho \nabla^s \mathbf{u} + \kappa_{art} \nabla \rho \otimes \mathbf{u} \\ \kappa_{art} \nabla E + \frac{1}{2} \mathbf{u} \cdot \mathbf{u} (\kappa_{art} \nabla \rho) + \mu_{art} \rho \nabla^s \mathbf{u} \cdot \mathbf{u} \end{bmatrix}(\mathbf{x}_i, t) \quad \forall \mathbf{x}_i \in \Omega,
    \end{split}
\end{equation}
\begin{equation}
    \begin{bmatrix} \rho^h \\ \rho \mathbf{u}^h\\ E^h \end{bmatrix}(\mathbf{x}_i, t) = \mathbf{g}(\mathbf{x}_i, t) \quad \forall \mathbf{x}_i \in \partial \Omega,
\end{equation}
\begin{equation}
    \begin{bmatrix} \rho^h \\ \rho \mathbf{u}^h\\ E^h \end{bmatrix}(\mathbf{x}_i, 0) = \begin{bmatrix} \rho_0 \\ \rho \mathbf{u}_0\\ E_0 \end{bmatrix}(\mathbf{x}_i) \quad \forall \mathbf{x}_i \in \Omega,
\end{equation}
\begin{equation}
    p^h(\mathbf{x}_i,t) = (\gamma - 1) (E^h(\mathbf{x}_i,t) - \frac{1}{2} \rho \mathbf{u}^h(\mathbf{x}_i,t) \cdot \mathbf{u}^h(\mathbf{x}_i,t)),
\end{equation}
\begin{equation}
    T^h(\mathbf{x}_i,t) = \frac{p^h(\mathbf{x}_i,t)}{\rho^h(\mathbf{x}_i,t)}.
\end{equation}
}
\right.
$$

\bigskip

\noindent Here we have included the viscosities $\kappa_{art}$ and $\nu_{art}$ within the derivative terms for consistency with \cite{nazarov2017investigation}, but we reiterate that terms involving derivatives of $\kappa_{art}$ and $\nu_{art}$ are ignored. The construction of these viscosities is done similarly to the previous cases. In particular we define 

\begin{equation}
    \mu_{art}(\mathbf{x}_i, t) = \min(\mu_{RB}(\mathbf{x}_i, t), \mu_{FO}(\mathbf{x}_i, t)),
\end{equation}

\noindent and 

\begin{equation}
    \kappa_{art}(\mathbf{x}_i, t) = \frac{\mathcal{P}}{C_{RB}}\mu_{art}(\mathbf{x}_i, t),
\end{equation}

\noindent similar to \cite{nazarov2017investigation}, with $\mathcal{P}$ representing an artificial Prandtl number, and the first order viscosity defined by the maximum wavespeed estimate of $c = |\mathbf{u}| + \sqrt{\gamma T}$. Unlike \cite{nazarov2017investigation}, however, we use the residuals of the governing conservation laws to determine the residual-based viscosity $\mu_{art}(\mathbf{x}_i, t)$. Like the previous case we determine the viscosity $\mu_{art}$ at every collocation point as the maximum of the viscosities computed using Equation \eqref{eq:final_visc} for each of the mass, momentum, and energy equations. 

\subsection{Linear Stabilization}

For the stability of the method, it is possible that the nonlinear viscosity alone is sufficient to prevent blow up of the solution. Indeed in the Galerkin finite element setting, residual-based viscosity stabilization results in a method which provably converges to the entropy solution of a scalar conservation law \cite{nazarov2013convergence}. However, as noted in \cite{burman2023some}, if only nonlinear stabilization is applied, the solution may still contain spurious small-scale oscillations. This was also noticed when using central finite difference stencils in \cite{stiernstrom2021SBP}. Thus it is still desirable to include some sort of linear stabilization, akin to upwinding in the finite difference and finite volume contexts. 

One option for a linear stabilization scheme is the Streamline-Upwind-Petrov-Galerkin (SUPG) method, recently adapted to the spline collocation setting in \cite{aronson2023stabilized}. In this work, however, we wish to develop a novel stabilization strategy inspired by projection stabilization for Galerkin methods \cite{braack2006lps}. Local projection stabilization adds a stabilization term to the weak formulation of the problem which depends on the difference between the discrete gradient and its projection onto a coarser finite element space, which creates a method which removes spurious oscillations. For Galerkin methods, the resulting scheme does not require the use of space-time elements for consistent time integration, and is instead more amenable to the method of lines approach pursued here. Moreover, we note that residual-based linear stabilization schemes like SUPG can exhibit instabilities when small time steps are employed \cite{hsu2010improving}, which we anticipate due to our use of explicit Runge-Kutta discretizations in time. 

Rather than project our discrete gradient onto a coarser function space, we make use of interpolation operations which we have shown previously can be performed in an efficient manner. Let the space $S^{h, k-1}$ be the B-spline space with the same number of elements as the solution space $S^h$, but with degree $k-1$ instead of $k$. One possible approach to stabilization is to interpolate the gradient of the solution using a function in $S^{h, k-1}$ to define a coarsened gradient used in the stabilization. However, we note that in 1D, the derivative of a B-spline function is itself a B-spline of one lower degree. Thus this coarsened gradient interpolation would be identical to the solution gradient, and the stabilization procedure would have no effect. 

To circumvent this issue, we start by interpolating the gradient of the solution using $S^h$ and its corresponding Greville points, and denote the result as $\widehat{\nabla \phi^h} \in S^h$. We then use the values of $\widehat{\nabla \phi^h}$ at the Greville points of $S^{h, k-1}$ to define another interpolation of the gradient, which we denote as $\Pi \widehat{\nabla \phi^h} \in S^{h, k-1}$. This function is used in a manner similar to the projected gradient in projection stabilization. Specifically, we evaluate its derivative at the Greville points of $S^h$, and another stabilization term is added to the statement of the collocation scheme of the form $\nu_{lin}(\mathbf{x}_i, t)  \nabla \cdot ( \nabla \phi^h (\mathbf{x}_i, t) - \Pi \widehat{\nabla \phi^h} (\mathbf{x}_i, t) )$, where the viscosity coefficient $\nu_{lin}$ is defined at each collocation point by 

\begin{equation}
    \nu_{lin}(\mathbf{x}_i, t) = C_{lin} h(\mathbf{x}_i) c(\mathbf{x}_i, t).
\end{equation} 

\noindent Here $C_{lin}$ is another tunable parameter and $c(\mathbf{x}_i)$ is the local wavespeed as described previously. As a final remark, we note that we have also tested other designs for this linear stabilization term, for example using a spline space with half the number of elements and the same degree as $S^h$ as the coarse space, instead of $S^{h,k-1}$. However, we have found that many of these schemes are overly dissipative on smooth solutions and thus the errors of the stabilized collocation scheme increase compared to the unstabilized scheme in these cases.

For completeness, we conclude by stating the fully stabilized semi-discrete form of a scalar conservation law, given by

$$
\left\{ \hspace{5pt}
\parbox{5in}{
\noindent Given $\mathbf{f} : \mathbb{R} \rightarrow \mathbb{R}^d$, $g: \partial \Omega \times [0,t_f] \rightarrow \mathbb{R}$, initial condition $\phi_0 : \Omega \rightarrow \mathbb{R}$, and final time $t_f \in \mathbb{R}$, find $\phi^h \in S^h_T$ such that:
\begin{equation}
\begin{split}
   \frac{\partial \phi^h}{\partial t}(\mathbf{x}_i, t) + &\nabla \cdot \mathbf{f} (\phi^h(\mathbf{x}_i,t)) = \nu_{art}(\mathbf{x}_i, t) \Delta \phi^h (\mathbf{x}_i, t) \\& + \nu_{lin}(\mathbf{x}_i, t)   \nabla \cdot ( \nabla \phi^h (\mathbf{x}_i, t) - \Pi \widehat{\nabla \phi^h} (\mathbf{x}_i, t) ) \quad \forall \mathbf{x}_i \in \Omega,
   \end{split}
\end{equation}
\begin{equation}
    \phi^h(\mathbf{x}_i, t) = g(\mathbf{x}_i, t) \quad \forall \mathbf{x}_i \in \partial \Omega,
\end{equation}
\begin{equation}
    \phi^h(\mathbf{x}_i, 0) = \phi_0(\mathbf{x}_i) \quad \forall \mathbf{x}_i \in \Omega.
\end{equation}
}
\right.
$$


For the Euler equations, the linear stabilization method is applied to the conserved quantity corresponding to each conservation law. For brevity we only list the semi-discrete form associated with the Guermond-Popov fluxes below

 $$
\left\{ \hspace{5pt}
\parbox{5in}{
\noindent Given $\mathbf{g} : \partial \Omega \times [0,t_f] \rightarrow \mathbb{R}^{d+2}$, initial conditions $\rho_0 : \Omega \rightarrow \mathbb{R}$, $\rho \mathbf{u}_0 : \Omega \rightarrow \mathbb{R}^d$, and $E_0 : \Omega \rightarrow \mathbb{R}$, and final time $t_f \in \mathbb{R}$, find $\rho^h \in S^h_T$, $\rho \mathbf{u}^h \in \mathbf{S}^h_T$, and $E^h \in S^h_T$ such that:
\begin{equation}
\begin{split}
    \frac{\partial}{\partial t} 
    \begin{bmatrix} \rho^h \\ \rho \mathbf{u}^h\\ E^h \end{bmatrix}(\mathbf{x}_i, t)
    + \nabla \cdot
    \begin{bmatrix} \rho \mathbf{u}^h \\ \rho \mathbf{u}^h \otimes \mathbf{u}^h + p^h \mathbf{I} \\ \mathbf{u}^h (E^h + p^h) \end{bmatrix}(\mathbf{x}_i,t)
    \\ =  \nabla \cdot
    \begin{bmatrix} \kappa_{art} \nabla \rho \\ \mu_{art} \rho \nabla^s \mathbf{u} + \kappa_{art} \nabla \rho \otimes \mathbf{u} \\ \kappa_{art} \nabla E + \frac{1}{2} \mathbf{u} \cdot \mathbf{u} (\kappa_{art} \nabla \rho) + \mu_{art} \rho \nabla^s \mathbf{u} \cdot \mathbf{u} \end{bmatrix}(\mathbf{x}_i, t) \\ +  \nu_{lin}(\mathbf{x}_i,t) \nabla \cdot
    \begin{bmatrix} \nabla \rho^h - \Pi\widehat{\nabla \rho^h} \\  \nabla \rho \mathbf{u}^h - \Pi\widehat{\nabla \rho \mathbf{u}^h} \\ \nabla E^h - \Pi\widehat{\nabla E^h}  \end{bmatrix}(\mathbf{x}_i,t) \quad \forall \mathbf{x}_i \in \Omega,
    \end{split}
\end{equation}
\begin{equation}
    \begin{bmatrix} \rho^h \\ \rho \mathbf{u}^h\\ E^h \end{bmatrix}(\mathbf{x}_i, t) = \mathbf{g}(\mathbf{x}_i, t) \quad \forall \mathbf{x}_i \in \partial \Omega,
\end{equation}
\begin{equation}
    \begin{bmatrix} \rho^h \\ \rho \mathbf{u}^h\\ E^h \end{bmatrix}(\mathbf{x}_i, 0) = \begin{bmatrix} \rho_0 \\ \rho \mathbf{u}_0\\ E_0 \end{bmatrix}(\mathbf{x}_i) \quad \forall \mathbf{x}_i \in \Omega,
\end{equation}
\begin{equation}
    p^h(\mathbf{x}_i,t) = (\gamma - 1) (E^h(\mathbf{x}_i,t) - \frac{1}{2} \rho \mathbf{u}^h(\mathbf{x}_i,t) \cdot \mathbf{u}^h(\mathbf{x}_i,t)),
\end{equation}
\begin{equation}
    T^h(\mathbf{x}_i,t) = \frac{p^h(\mathbf{x}_i,t)}{\rho^h(\mathbf{x}_i,t)}.
\end{equation}
}
\right.
$$

\section{Results}

With a full description of the stabilized collocation schemes complete, we now move on to consider the results of numerical tests on a variety of conservation laws. We start by considering the linear advection and Burgers equations, which demonstrate the effectiveness of the stabilization strategies in the presence of shocks while maintaining high-order accuracy on problems with smooth solutions. We then move on to consider solutions of the Buckley-Leverett equations, which show that accurate solutions can be obtained even for scalar conservation laws defined by non-convex fluxes. Finally, we turn to the compressible Euler equations to show the promise of the method for high-order CFD applications. In all cases we perform time integration using the standard RK4 method. In addition, results are obtained with both linear and nonlinear stabilization terms included, unless otherwise stated in the text. 

\subsection{Scalar Advection}

The simplest scalar conservation law is the linear advection equation, given by

\begin{equation}
    \frac{\partial \phi}{\partial t} + \nabla \cdot (\mathbf{a} \phi) = 0.
    \label{eq:scalar_adv}
\end{equation}

\noindent Here $\mathbf{a}$ is the advection velocity field and $\phi$ is the conserved quantity. In this case the solution is simply transported according to the velocity $\mathbf{a}$, meaning that, while discontinuities may be present in the initial condition, shocks will not form from smooth solutions in finite time. For all of the cases in this section we define $C_{RB} = 4$, $C_{max} = 0.5$, and $C_{lin} = 0.25$. 

\subsubsection{Smooth Solution}

To start, we consider a scenario with a smooth solution in order to assess the effect of the stabilization strategies on errors produced by the collocation schemes. We define the initial condition as $\phi(\mathbf{x}, 0) = \sin(2 \pi x) \sin(2 \pi y)$ over the domain $\mathbf{x} \in [0,1]^2$, a constant advection velocity of $\mathbf{a} = (1,1)$, and use periodic boundary conditions in each direction. From this setup it is clear that the exact solution at $t_f = 1$ is the same as the initial condition, $\phi(\mathbf{x}, 0) = \phi(\mathbf{x}, 1)$, and we compute the error in the numerical solution at this time compared to the exact solution. 

\begin{figure}
\centering
\subfloat[Linear stabilization only]{\label{sfig:Adv2D_smooth_upwind}\includegraphics[width=.5\textwidth]{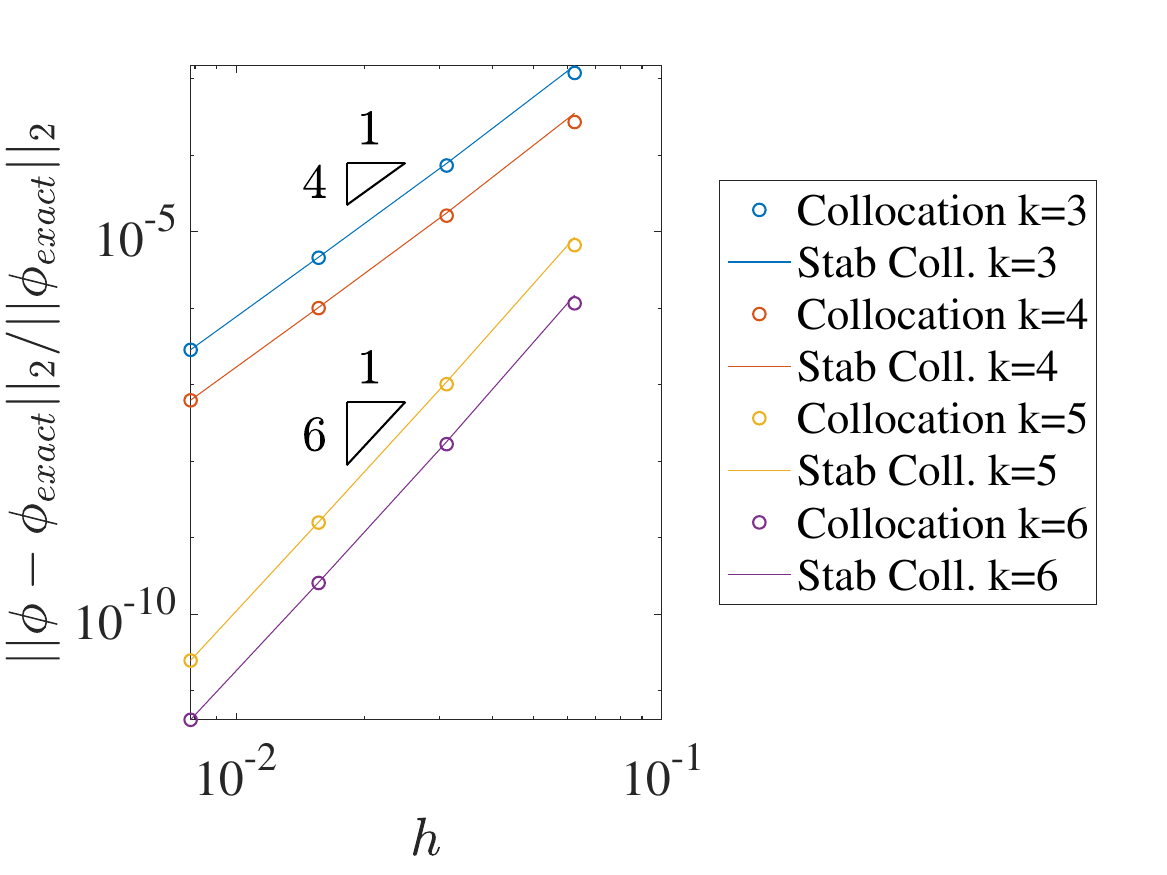}}\hfill
\subfloat[Linear and nonlinear stabilization]{\label{sfig:Adv2D_smooth_both}\includegraphics[width=.5\textwidth]{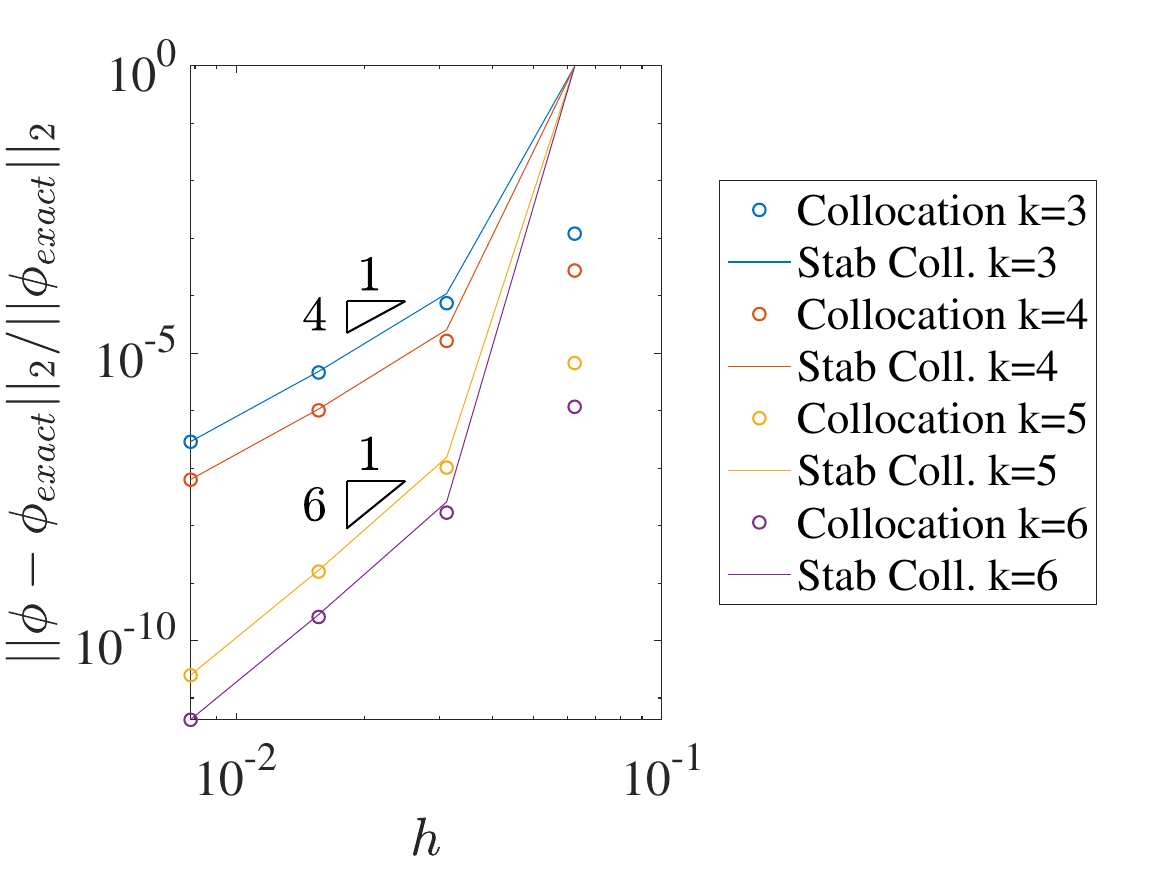}} \\
\caption{Convergence rates to smooth solution of advection equation}
\label{fig:Adv2D_smooth_conv}
\end{figure}

Figure \ref{fig:Adv2D_smooth_conv} shows the results of this study using a time step of $1 \times 10^{-4}$, selected experimentally to be small enough such that the time discretization errors are small compared to the spatial discretization errors for all of the meshes tested. In Figure \ref{sfig:Adv2D_smooth_upwind} we show the $L^2$ errors obtained with standard, unstabilized collocation as well as with only linear stabilization included. Clearly the expected convergence rates of $k+1$ for odd $k$ and $k$ for even $k$ are recovered using both schemes and the inclusion of linear stabilization has not noticeably increased the magnitude of the errors in this smooth case.

Figure \ref{sfig:Adv2D_smooth_both} shows the $L^2$ errors obtained with linear and nonlinear stabilization active compared to unstabilized collocation. As the mesh is refined, the errors obtained with the stabilized scheme do match the unstabilized results and the proper convergence rates are obtained. However, on the coarsest mesh tested we see that the inclusion of nonlinear stabilization has dramatically increased the error. On this coarse mesh the solution is treated as underresolved by the stabilization and thus the nonlinear viscosity is active everywhere and behaves like a first-order viscosity term. After time integrating for one full period the resulting damping has destroyed the structure of the numerical solution, resulting in the large errors. 

\subsubsection{Non-Smooth Solution}

Next we consider the transport of a non-smooth solution profile. While the previous test case demonstrated the ability of the stabilization to properly vanish in smooth regions where the solution is well resolved, this test will demonstrate the effectiveness of the stabilization near shocks. The setup is the same as in the previous case, only now we define the initial condition as 

\begin{equation}
    \phi(\mathbf{x}, 0) = 
    \begin{cases}
    1 & \textup{if } \mathbf{x} \in (0.3, 0.7)^2 \\
    0 & \textup{Otherwise}.
    \end{cases}
\end{equation}

\begin{figure}
\centering
\subfloat[$L^2$ errors]{\label{sfig:Adv2D_nonsmooth_l2}\includegraphics[width=.5\textwidth]{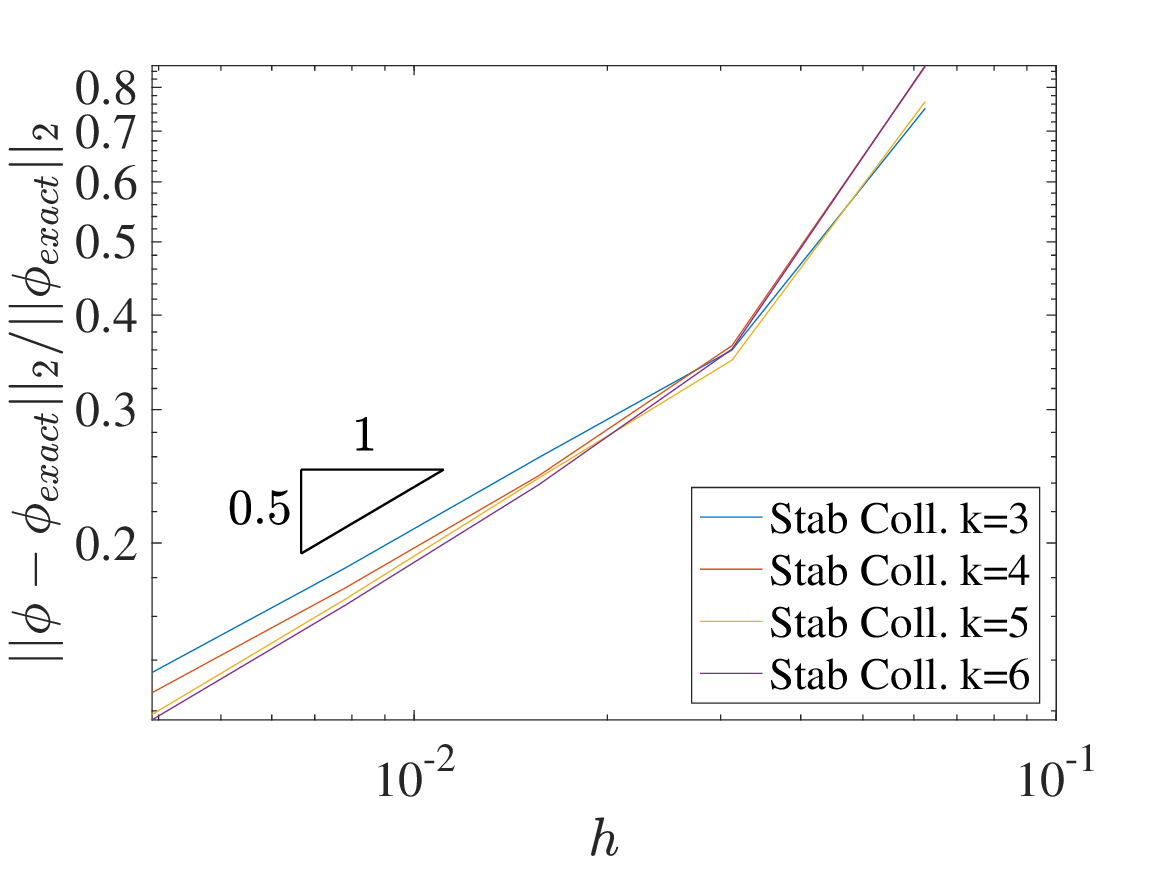}}\hfill
\subfloat[$L^1$ errors]{\label{sfig:Adv2D_nonsmooth_l1}\includegraphics[width=.5\textwidth]{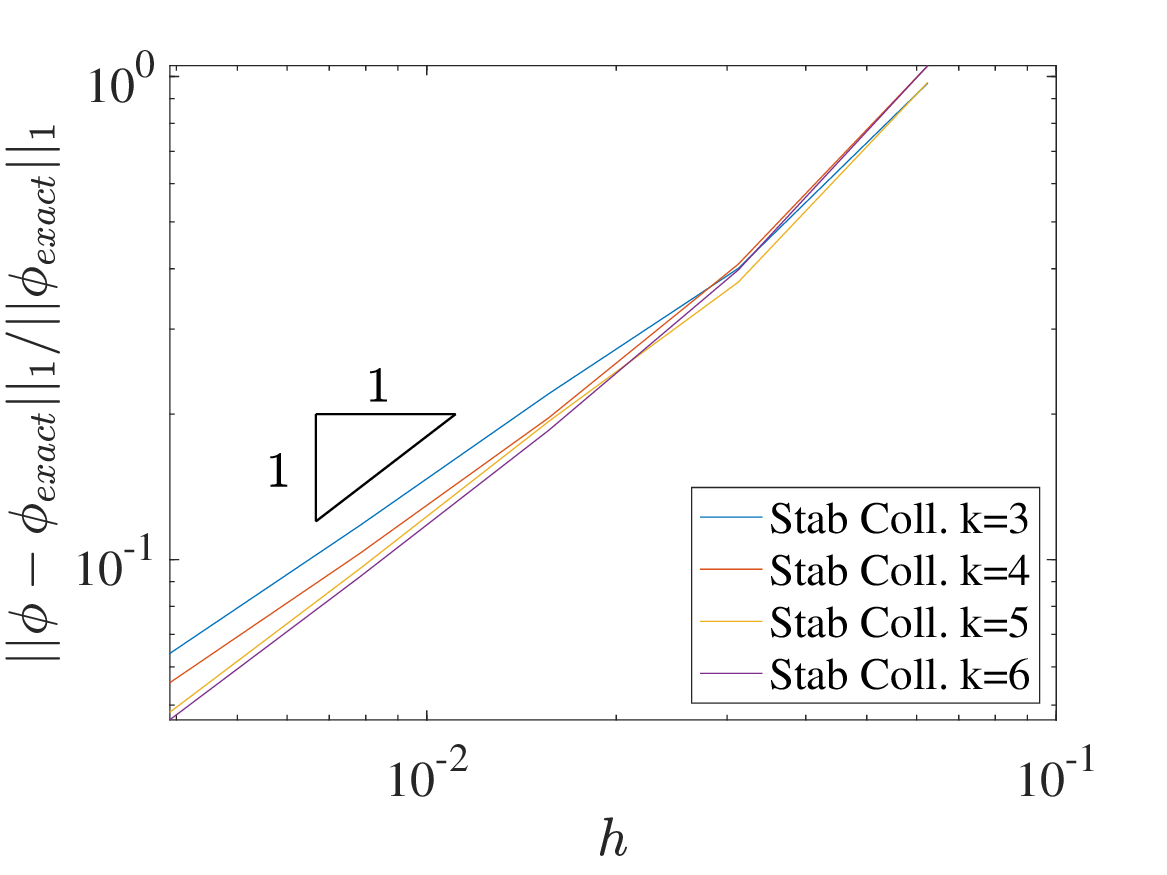}} \\
\caption{Convergence rates to non-smooth solution of advection equation}
\label{fig:Adv2D_nonsmooth_conv}
\end{figure}

Figure \ref{fig:Adv2D_nonsmooth_conv} shows the $L^1$ and $L^2$ errors in the solution after advancing in time for one period. As the solution in this case contains a shock, the optimal theoretical convergence rates are 1 in the $L^1$ norm and 1/2 in the $L^2$ norm. Our stabilized collocation scheme recovers close to these optimal rates with all degrees $k$, with a small increase in rate and decrease in error magnitude as $k$ increases. 

\begin{figure}
\centering
\subfloat[First-order viscosity]{\label{sfig:Adv2D_nonsmooth_FO}\includegraphics[width=.5\textwidth]{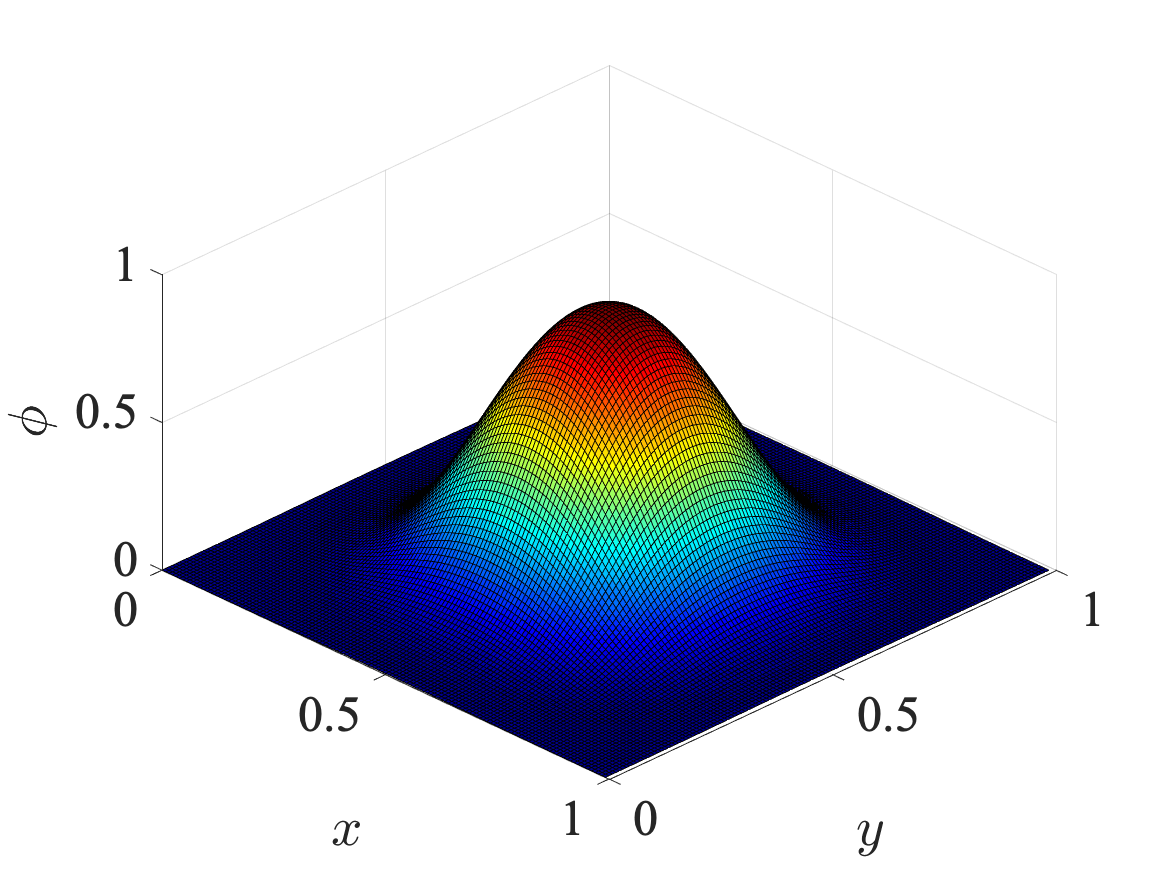}}\hfill
\subfloat[Residual-based viscosity]{\label{sfig:Adv2D_nonsmooth_RB}\includegraphics[width=.5\textwidth]{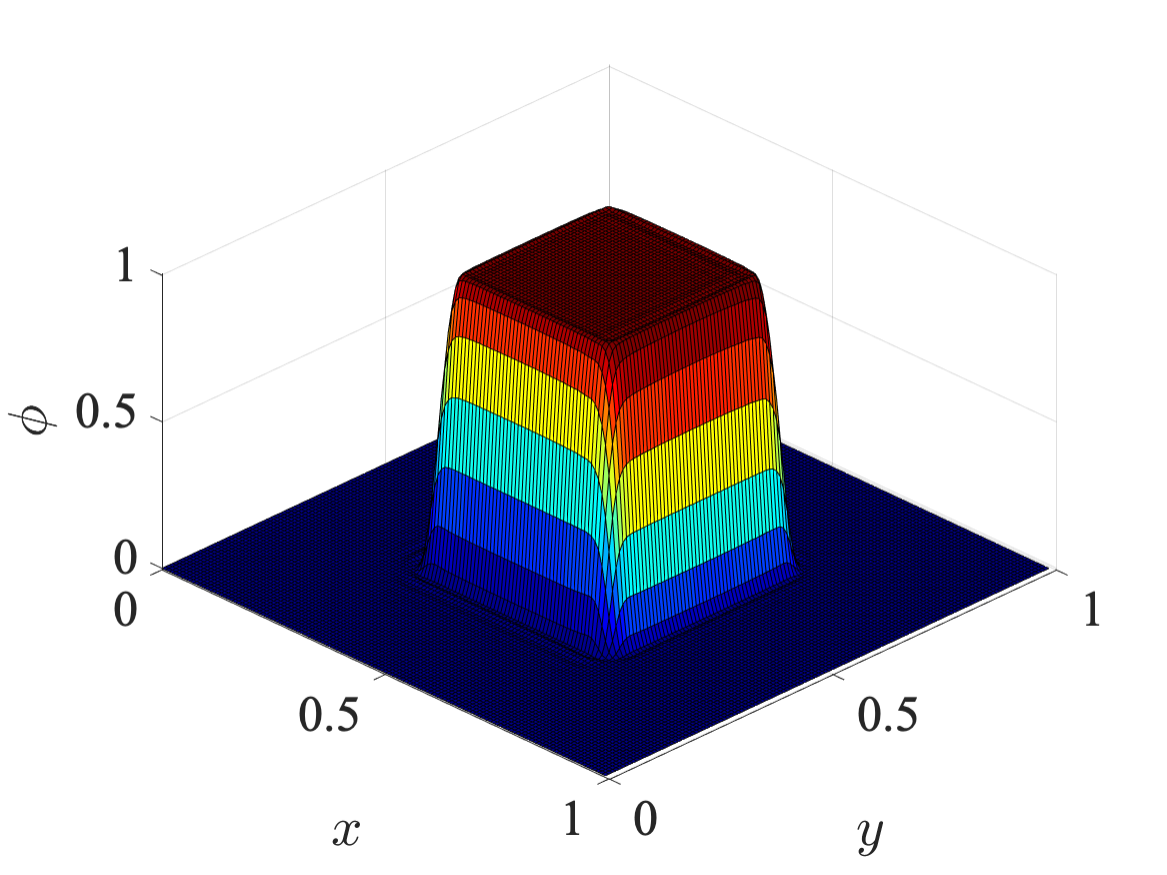}} \\
\caption{Comparison of nonlinear stabilization techniques on non-smooth solution of advection equation}
\label{fig:Adv2D_nonsmooth_FOcomp}
\end{figure}

Moreover, we can also demonstrate the superior performance of the residual-based viscosity when compared to a first-order viscosity using this example. Figure \ref{fig:Adv2D_nonsmooth_FOcomp} shows the final solution using $128^2$ B\'{e}zier elements and $k = 5$ with a first-order viscosity (obtained by overwriting the residual-based viscosity with the limiting first-order value at all locations) and our residual-based stabilization scheme. As expected, the residual-based scheme recovers a much more accurate and sharp solution profile than the first-order viscosity scheme.

\begin{figure}
\centering
\subfloat[Solution contours: Nonlinear stabilization only]{\label{sfig:Adv2D_nonsmooth_RBonly_contours}\includegraphics[width=.5\textwidth]{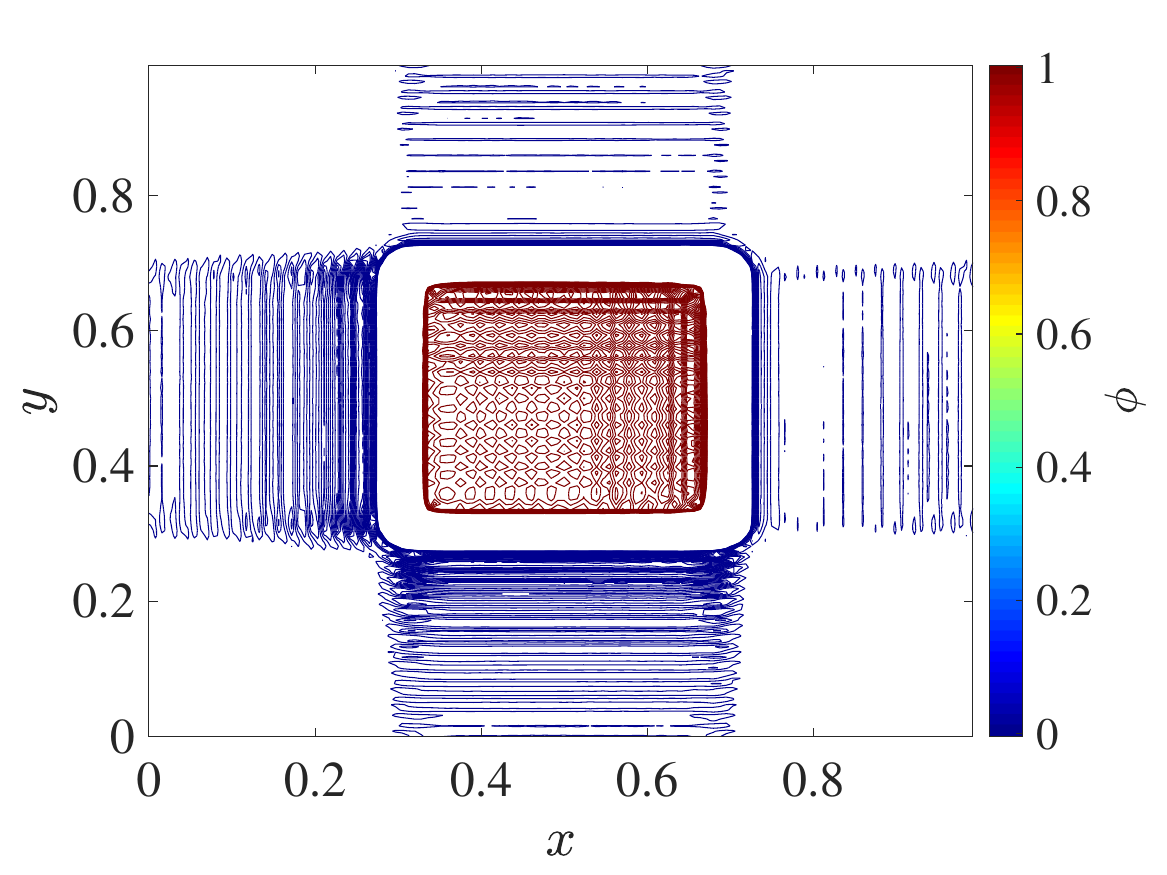}}\hfill
\subfloat[Viscosity field: Nonlinear stabilization only]{\label{sfig:Adv2D_nonsmooth_RBonly_nu}\includegraphics[width=.5\textwidth]{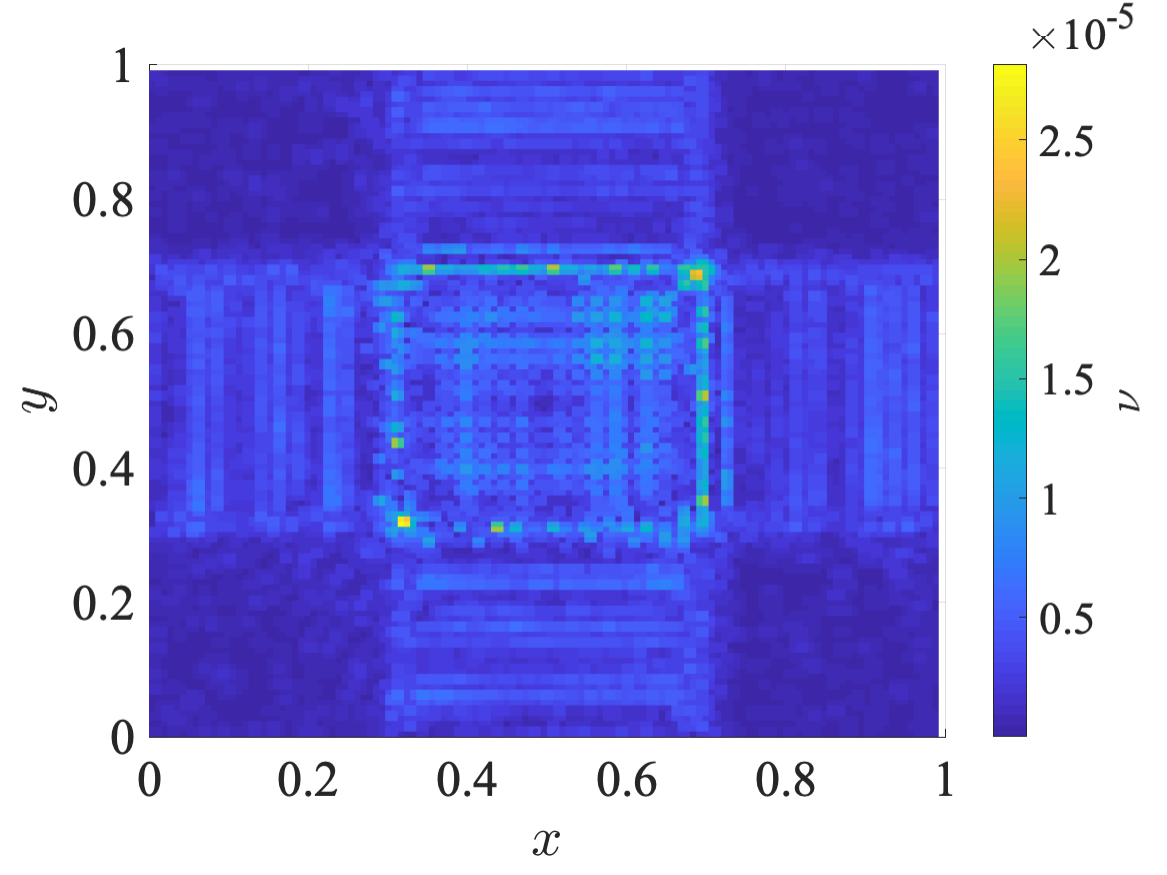}} \\
\subfloat[Solution contours: Linear and nonlinear stabilization]{\label{sfig:Adv2D_nonsmooth_both_contours}\includegraphics[width=.5\textwidth]{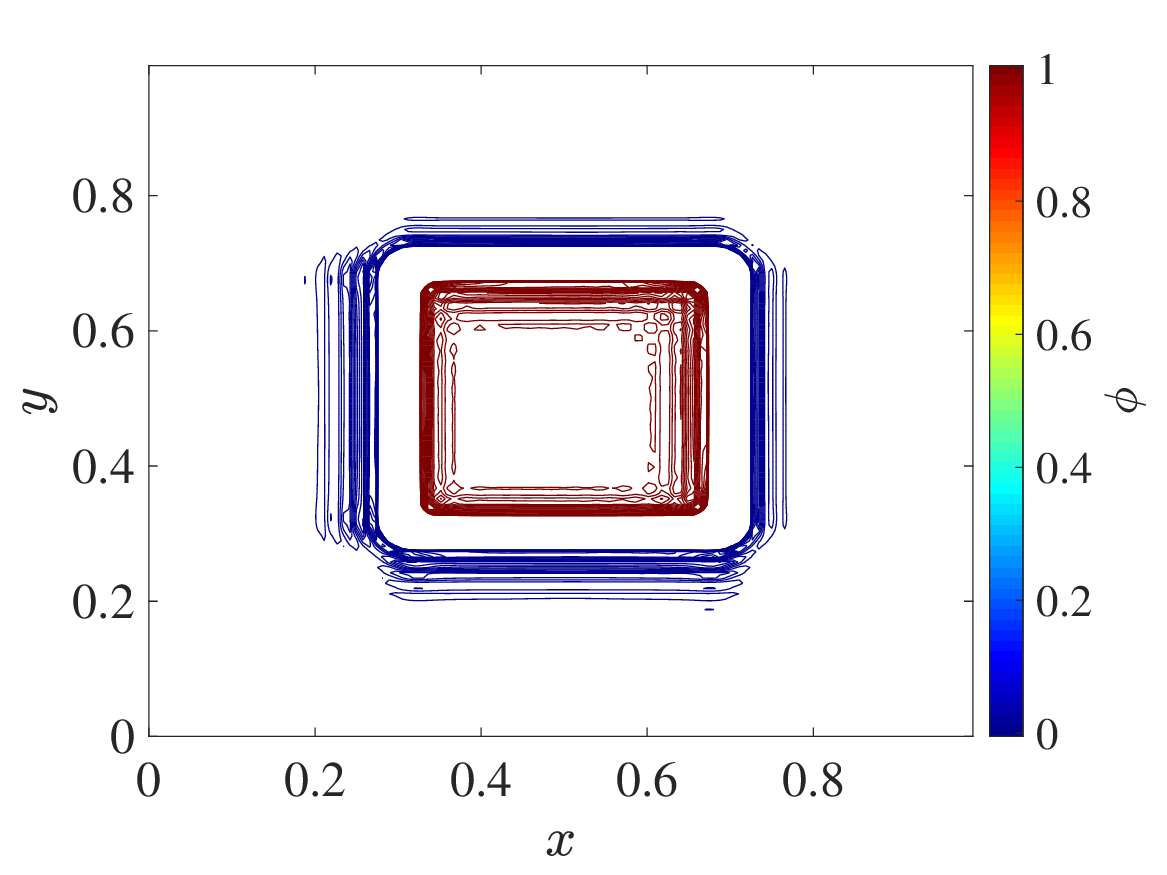}}\hfill
\subfloat[Viscosity field: Linear and nonlinear stabilization]{\label{sfig:Adv2D_nonsmooth_both_nu}\includegraphics[width=.5\textwidth]{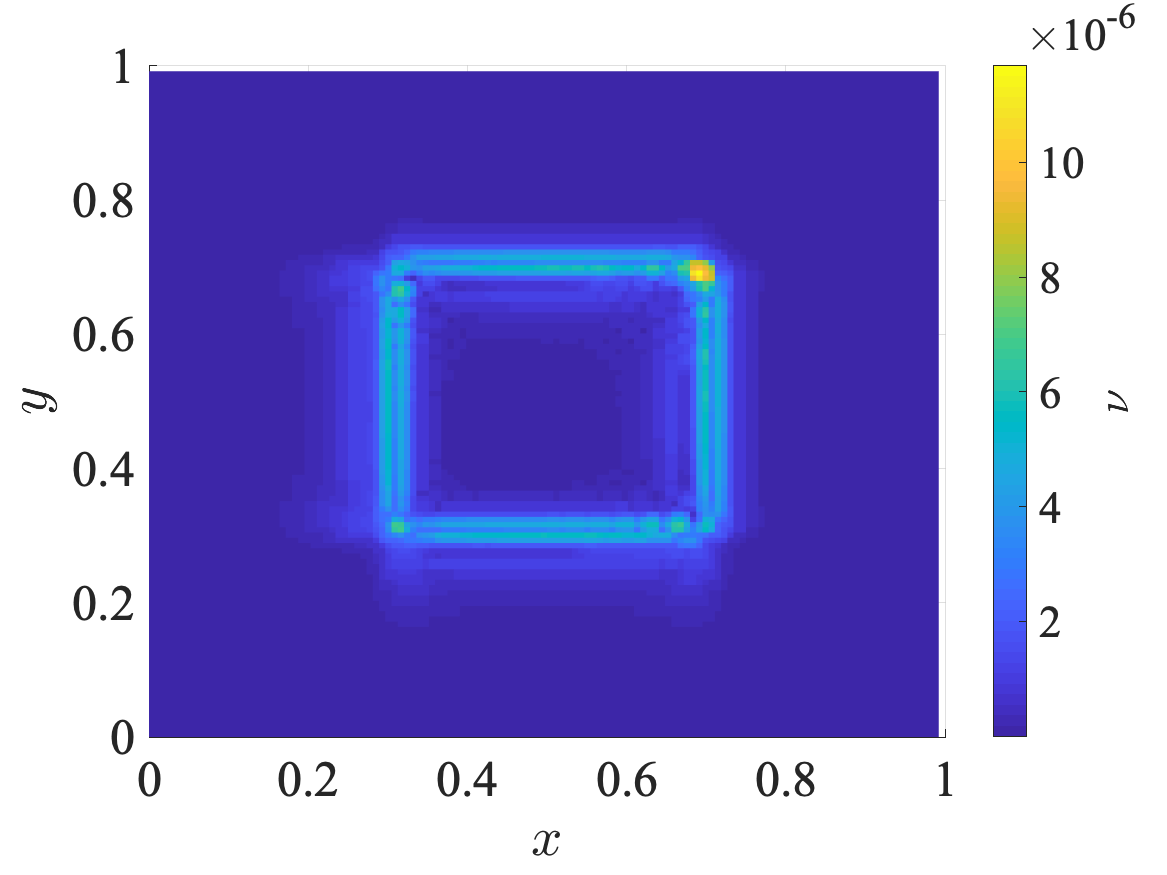}} \\
\caption{Effect of linear stabilization on non-smooth solution of advection equation}
\label{fig:Adv2D_nonsmooth_lincomp}
\end{figure}

Finally, this same test case also demonstrates the effect of the linear stabilization technique when coupled with the residual-based viscosity. Figure \ref{fig:Adv2D_nonsmooth_lincomp} shows contour plots of the solutions obtained with and without linear stabilization using 128 elements in each direction and $k=5$, as well as the corresponding residual-based viscosity fields. These plots of the solution include 30 equispaced contours in the range $[-0.01, 0.01]$ and 30 equispaced contours in $[0.99, 1.01]$. While the collocation scheme is certainly stable without linear stabilization, we see that there are many small oscillations that travel away from the shock locations. Moreover we see that the residual-based viscosity becomes active in the oscillatory regions, but it is not effective at removing them. With linear stabilization included, however, we see that these spurious oscillations are effectively removed and the residual-based viscosity is much more focused in the shock regions. 

\subsection{Burgers Equation}

Next we consider the inviscid Burgers equation, written in conservative form as 
\begin{equation}
    \frac{\partial \phi}{\partial t} + \nabla \cdot \left(\frac{\phi^2}{2} \mathbf{v} \right) = 0,
\end{equation}

\noindent with $\mathbf{v} = 1$ in 1D and $\mathbf{v} = (1,1)$ in 2D. Here the unknown $\phi$ represents the advection speed in the $\mathbf{v}$ direction, making the problem nonlinear and creating the possibility of shock formation in finite time. Similar to linear advection, we define the constants $C_{RB} = 4$, $C_{max} = 0.5$, and $C_{lin} = 0.25$ for all of the examples in this section. 

\subsubsection{Smooth Solution}

As in the linear advection case we start by considering a problem with a smooth analytical solution to verify convergence of the stabilized methods. We consider a 1D problem on the interval $x \in [0,1]$ with initial condition given by $\phi(x, 0) = e^x - 1$. We specify the exact solution as a Dirichlet boundary condition at $x = 0$ and do not specify anything at the outflow, $x = 1$. We advance in time using a time step of size $5 \times 10^{-5}$ until a final time of $t = 0.01$, where the exact solution is computed via the method of characteristics. 

\begin{figure}
\centering
\subfloat[Linear stabilization only]{\label{sfig:Burgers_smooth_upwind}\includegraphics[width=.5\textwidth]{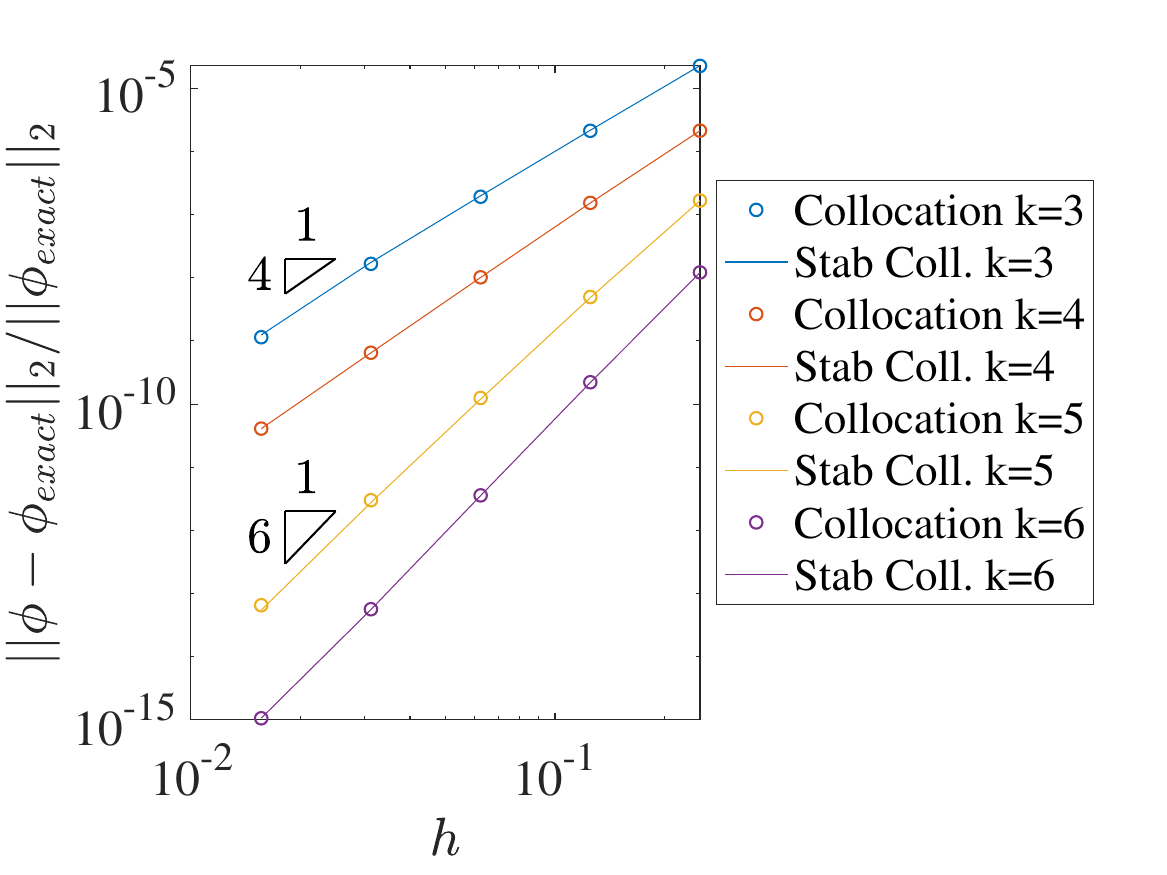}}\hfill
\subfloat[Linear and nonlinear stabilization]{\label{sfig:Burgers_smooth_both}\includegraphics[width=.5\textwidth]{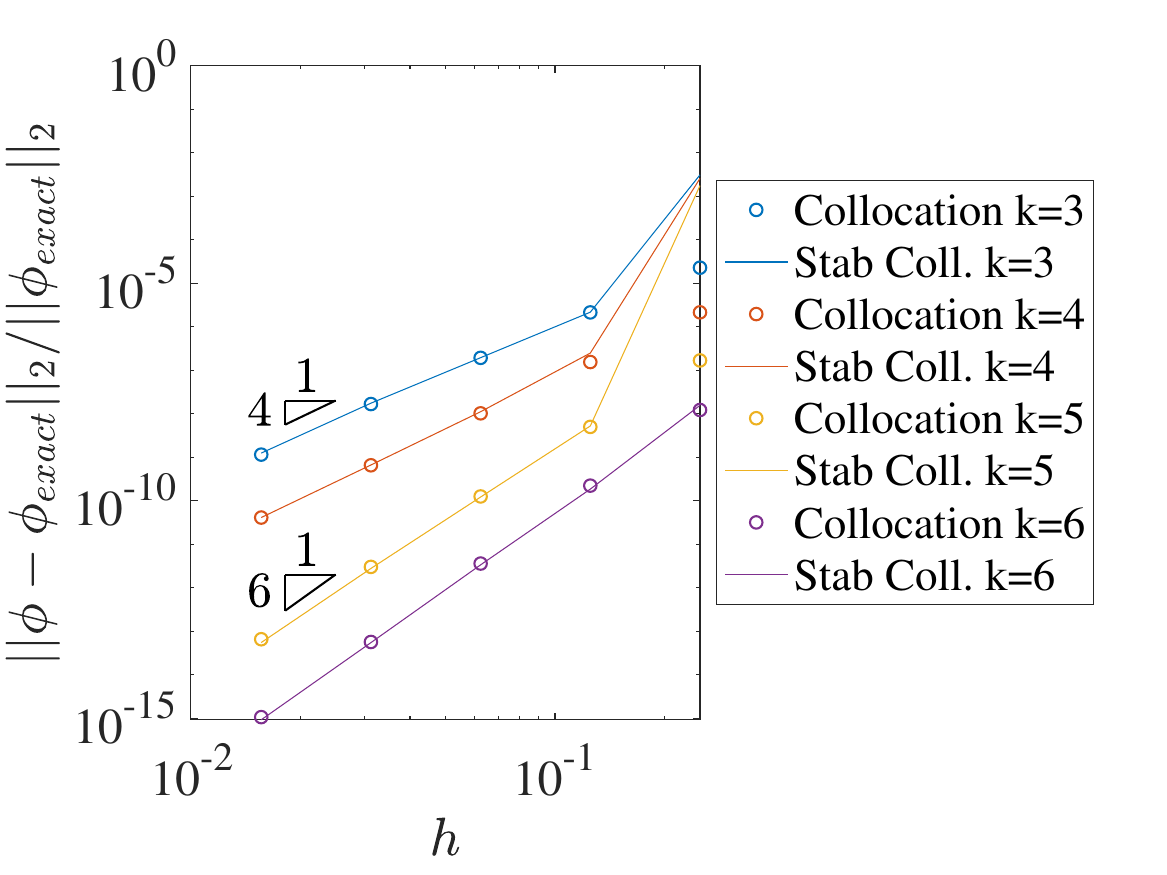}} \\
\caption{Convergence rates to smooth solution of Burgers equation}
\label{fig:Burgers_smooth_conv}
\end{figure}

Figure \ref{fig:Burgers_smooth_conv} details the $L^2$ errors obtained with unstabilized collocation, collocation with only linear stabilization, and collocation with linear and nonlinear stabilization. Like the scalar advection case we recover the expected convergence rates with the stabilized schemes and the errors are the same as those obtained without stabilization. We also see the same coarse mesh effect where the nonlinear stabilization becomes active and increases the error, but again this is limited to very coarse resolution simulations. 

\subsubsection{1D Riemann Problem}

We now move on to solving a 1D Riemann problem, namely a moving shock. The initial condition is defined over the interval $x \in [0,1]$ as 

\begin{equation}
    \phi(x, 0) = 
    \begin{cases}
    1 & \textup{if } x < 1/3 \\
    0 & \textup{Otherwise}.
    \end{cases}
\end{equation}

\noindent This shock travels to the right at a constant speed, and at the final time of $t_f = 0.2$ is located at $x = 13/30$. For this example we select a time step of $1 \times 10^{-5}$.

\begin{figure}
\centering
\subfloat[$L^2$ Error]{\label{sfig:Burgers_step_l2}\includegraphics[width=.5\textwidth]{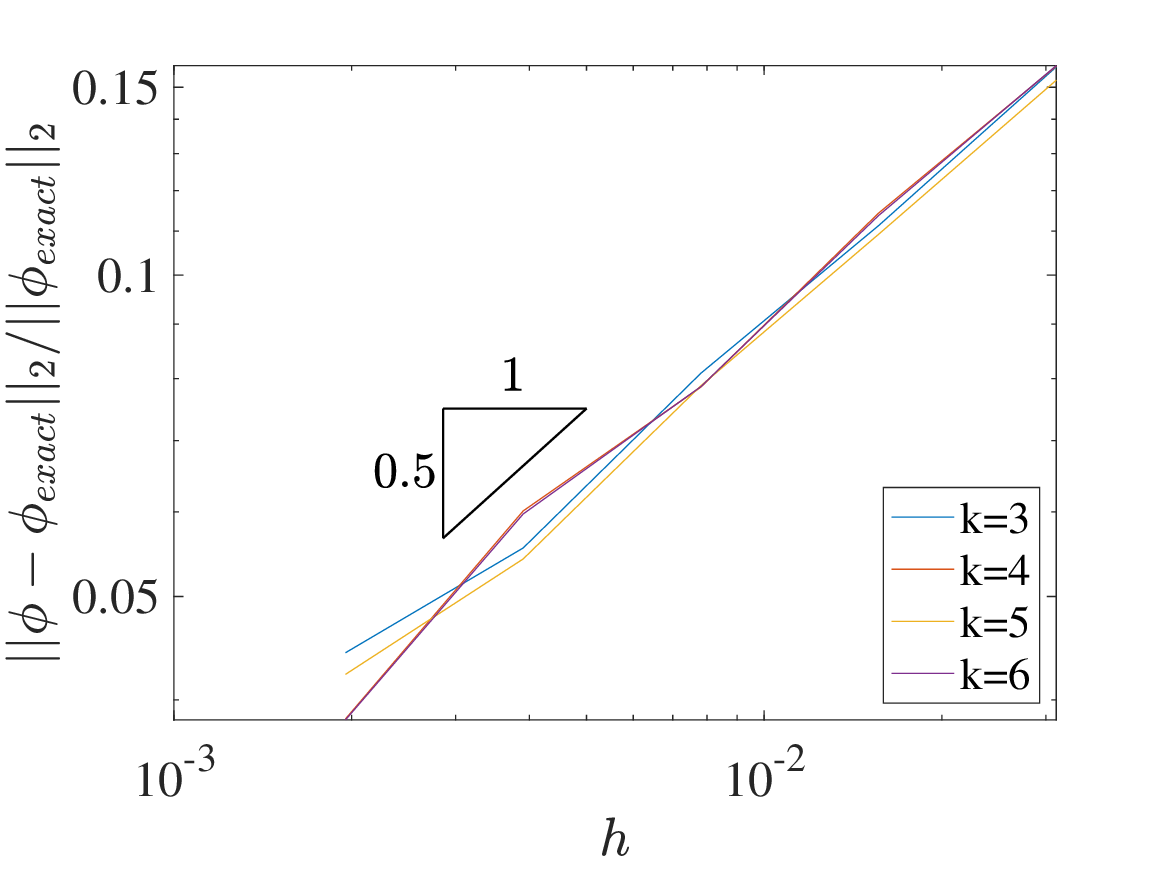}}\hfill
\subfloat[$L^1$ Error]{\label{sfig:Burgers_step_l1}\includegraphics[width=.5\textwidth]{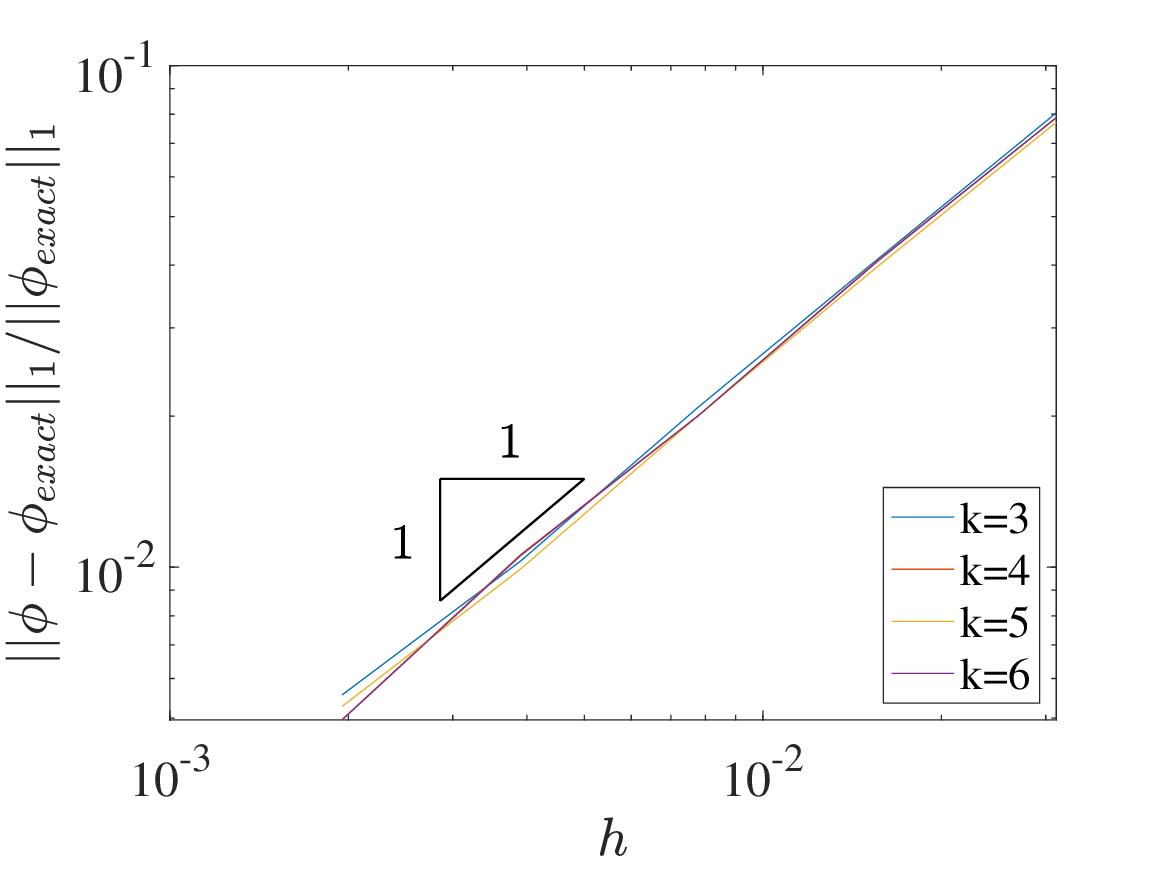}}\\
\caption{Convergence rates of 1D Burgers equation Riemann problem}
\label{fig:Burgers_step_conv}
\end{figure}

Figure \ref{fig:Burgers_step_conv} shows the convergence in the $L^2$ and $L^1$ norms for the stabilized collocation schemes. Once again, they recover close to the optimal rates of 0.5 and 1, respectively. We also see some interesting behavior where, especially in the $L^2$ error, as the mesh is refined the rate of convergence alternates between a higher and lower value. Moreover, it seems that the odd $k$ results are similar to one another, while the even $k$ are also similar but shifted by one refinement level compared to the odd results. The errors are comparable in magnitude to those obtained by residual-based viscosity SBP finite difference methods \cite{stiernstrom2021SBP}.

\begin{figure}
\centering
\subfloat[Solution]{\label{sfig:Burgers_elemComp_phi}\includegraphics[width=.5\textwidth]{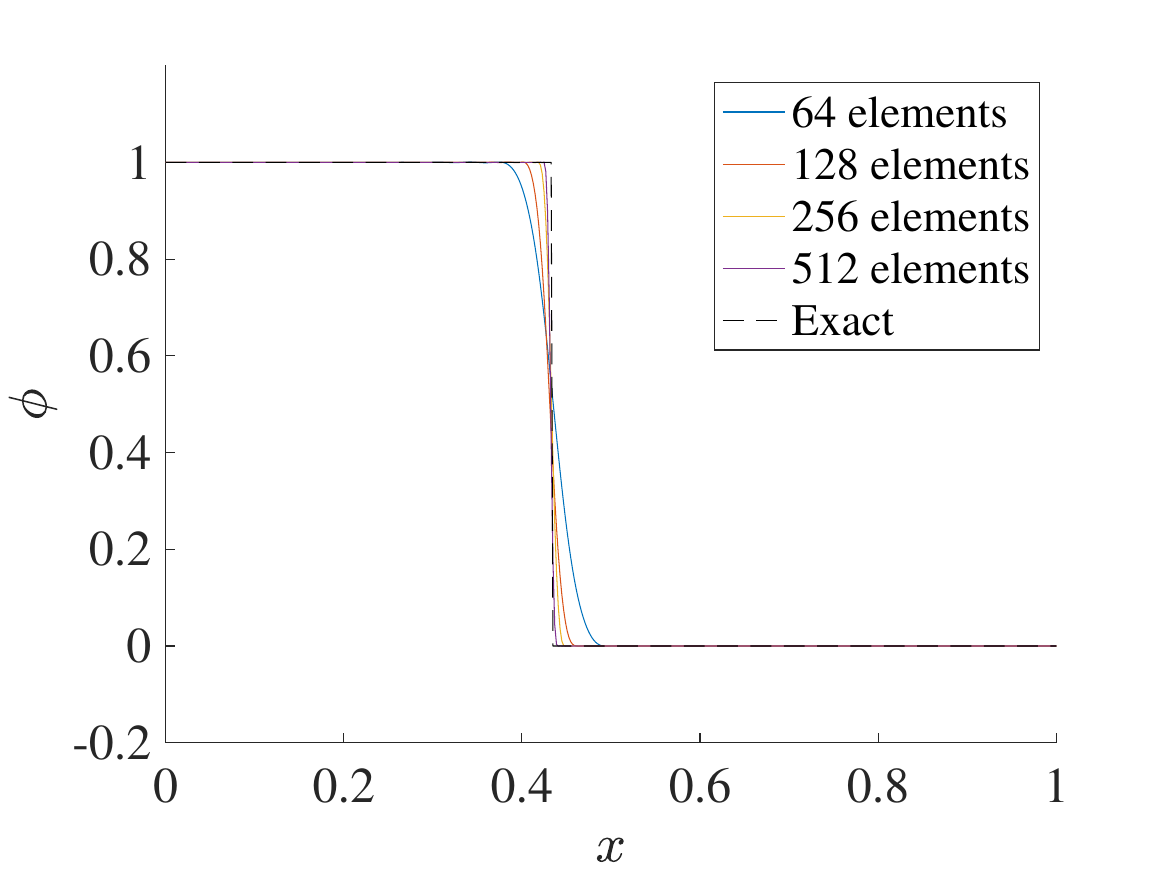}}\hfill
\subfloat[Residual-based viscosity]{\label{sfig:Burgers_elemComp_nu}\includegraphics[width=.5\textwidth]{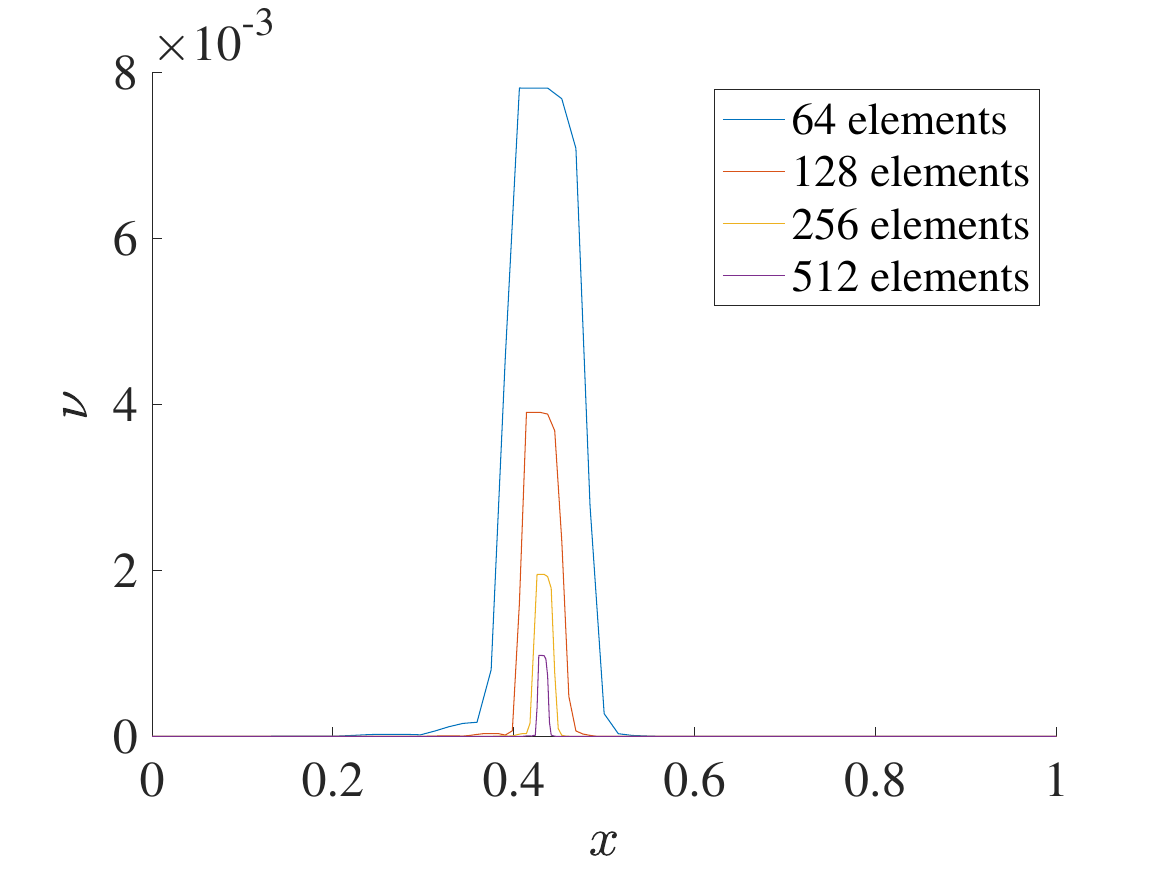}}\\
\caption{1D Burgers equation Riemann problem solution with varying number of elements}
\label{fig:Burgers_elemComp}
\end{figure}

\begin{figure}
\centering
\subfloat[Solution]{\label{sfig:Burgers_pComp_phi}\includegraphics[width=.5\textwidth]{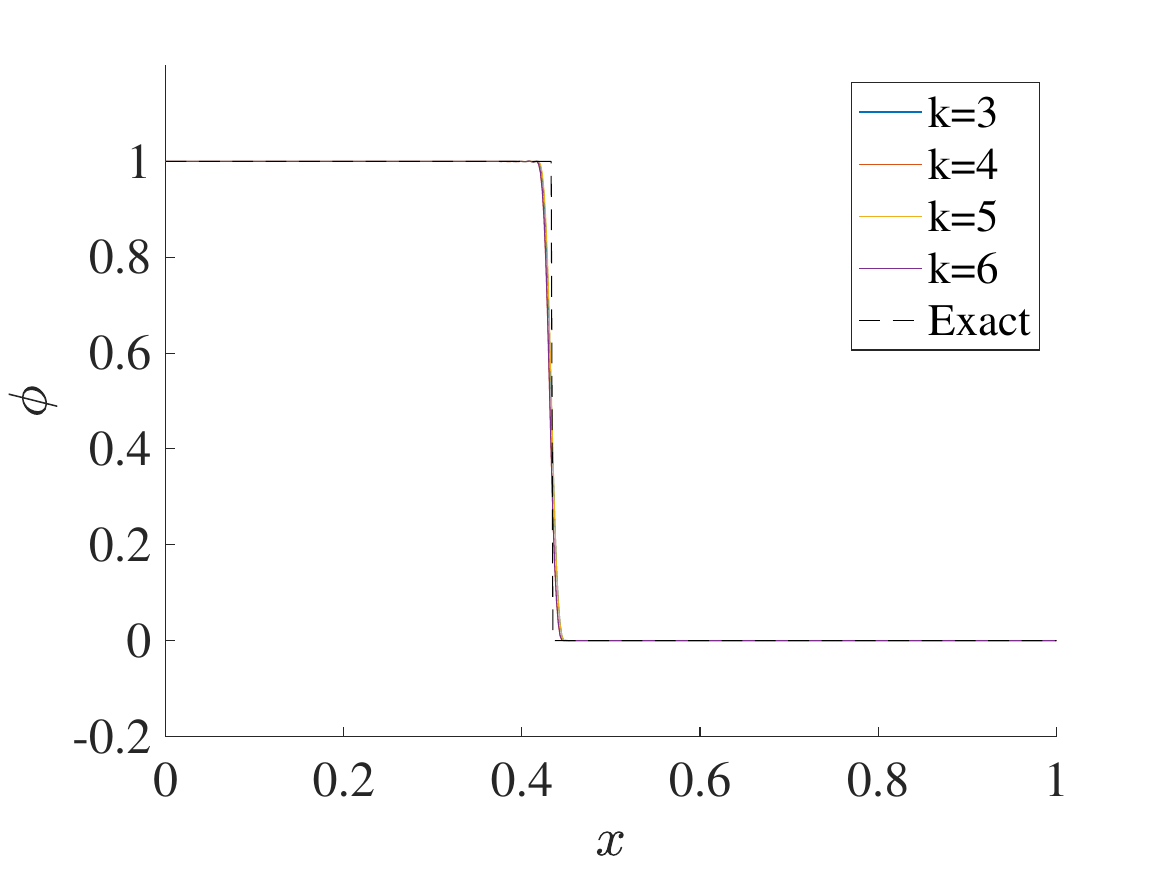}}\hfill
\subfloat[Residual-based viscosity]{\label{sfig:Burgers_pComp_nu}\includegraphics[width=.5\textwidth]{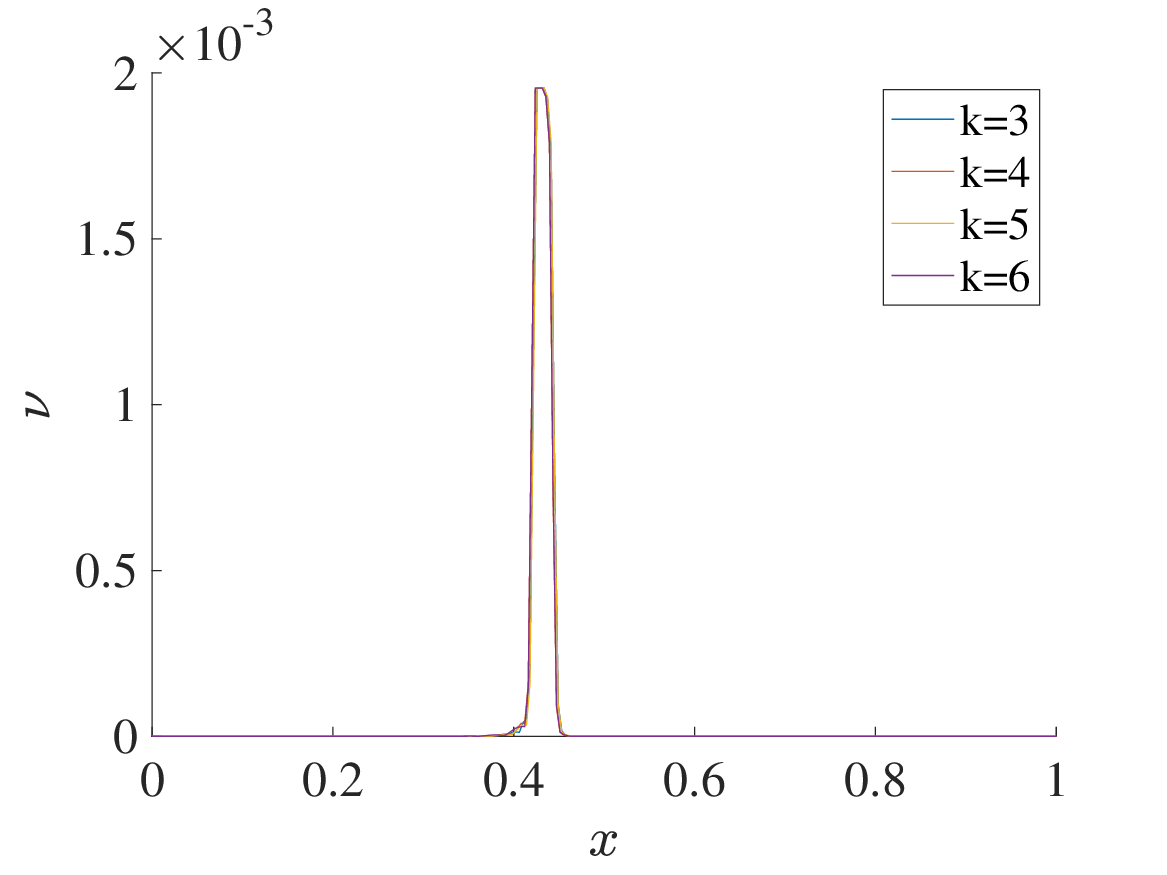}}\\
\caption{1D Burgers equation Riemann problem solution with varying polynomial degree}
\label{fig:Burgers_pComp}
\end{figure}

\begin{figure}
\centering
\subfloat[Mesh refinement]{\label{sfig:Burgers_elemComp_phi_zoom}\includegraphics[width=.5\textwidth]{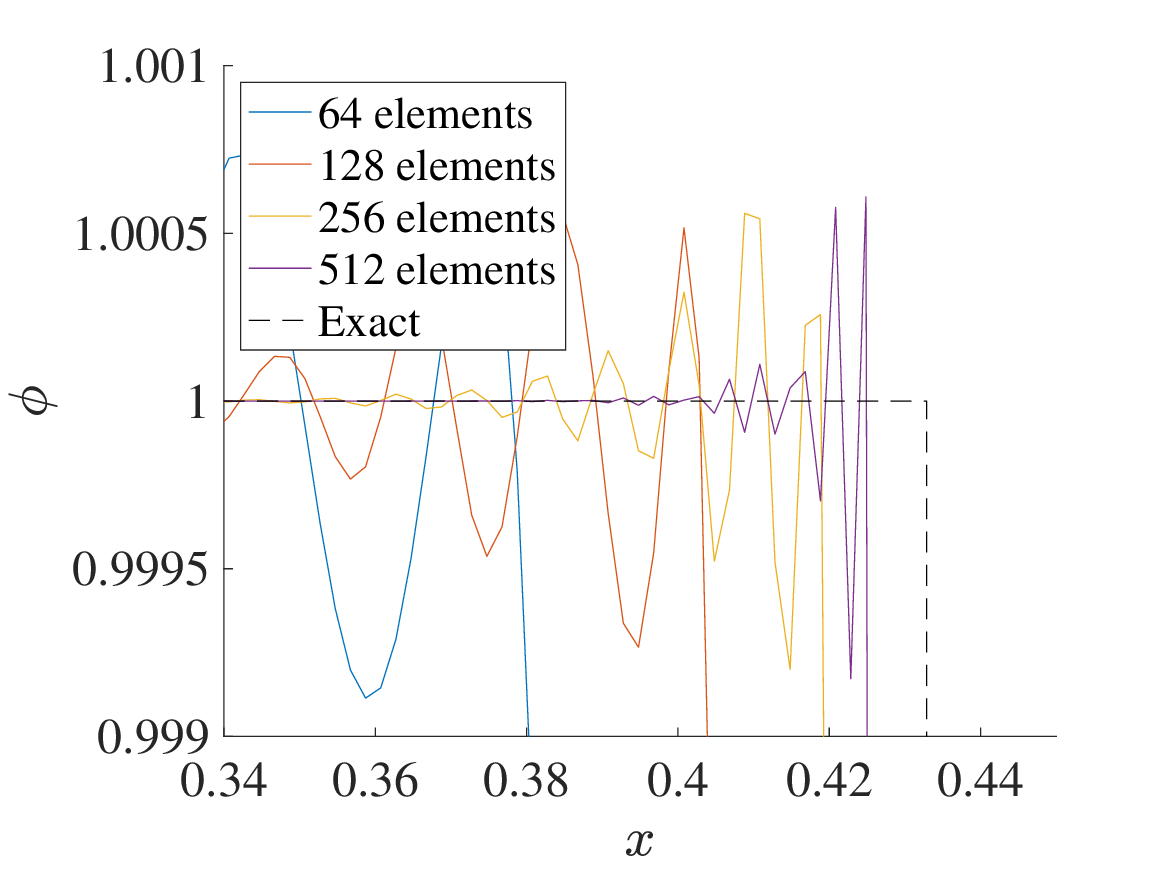}}\hfill
\subfloat[Degree elevation]{\label{sfig:Burgers_pComp_phi_zoom}\includegraphics[width=.5\textwidth]{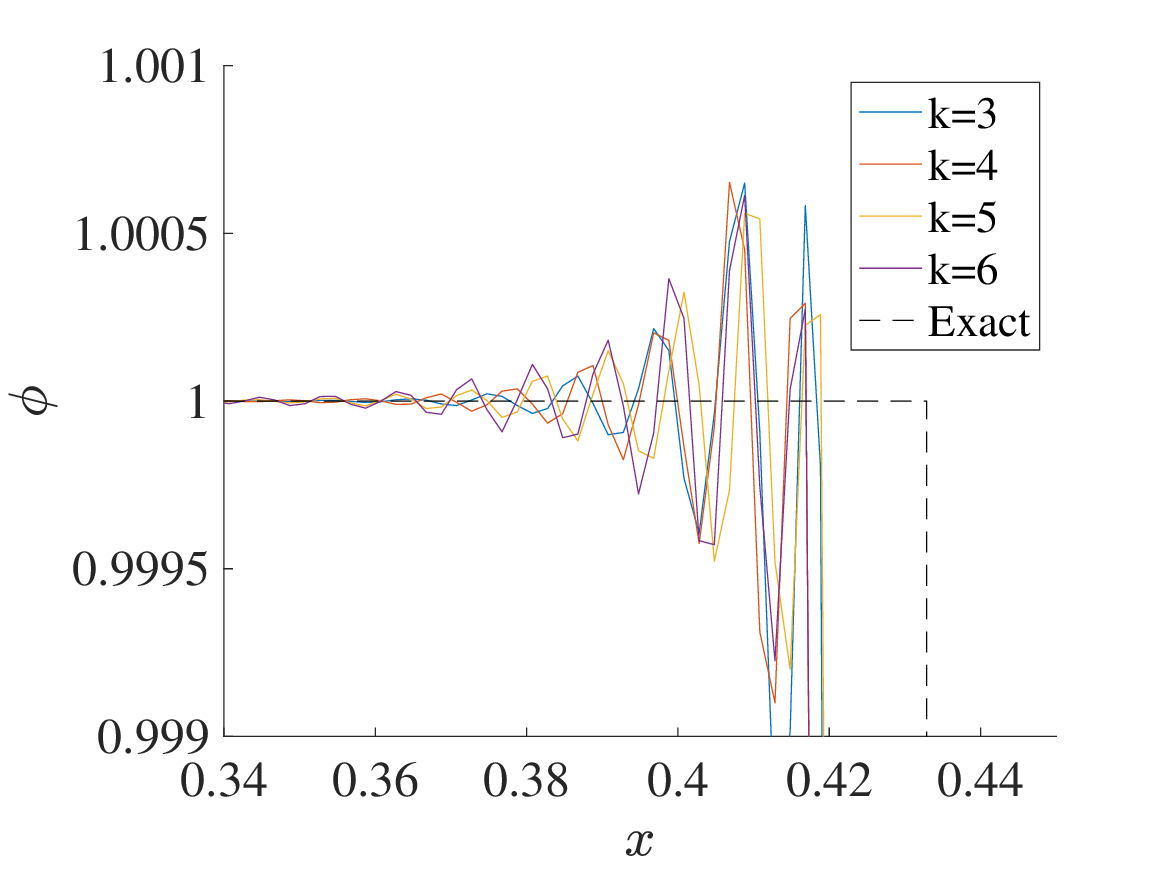}}\\
\caption{Zoom view of 1D Burgers equation Riemann problem solutions}
\label{fig:Burgers_zoom}
\end{figure}

For a more qualitative analysis we can also plot the solutions and viscosity fields for a variety of meshes and polynomial degrees. Figure \ref{fig:Burgers_elemComp} shows the solution and residual-based viscosity fields obtained using $k=5$ and a various numbers of elements, while Figure \ref{fig:Burgers_pComp} details the results obtained using 256 elements and a variety of polynomial degrees. The solutions behave as expected, with the shocks in simulations run with fewer elements being spread over a larger region, and variations in polynomial order having little effect. Figure \ref{fig:Burgers_zoom} shows a zoomed view of these solutions at the top of the shock location, illustrating the level of monotonicity achieved by the solutions. Increasing the number of elements results in the oscillations being confined to a region closer to the shock, while degree elevation again produces little change in the solution. The location and magnitude of the pre-shock oscillations are very similar to those seen in \cite{stiernstrom2021SBP} using upwind SBP finite difference schemes for the same problem.

\subsubsection{2D Riemann Problem}

Another common test case for Burgers equation is the two-dimensional Riemann problem defined on $\mathbf{x} \in [0,1]^2$ with the initial condition

\begin{equation}
    \phi(\mathbf{x}, 0) = 
    \begin{cases}
    0.5 & \textup{if } x < 1/2 \textup{ and } y < 1/2 \\
    -0.2 & \textup{if } x < 1/2 \textup{ and } y > 1/2 \\
    0.8 & \textup{if } x > 1/2 \textup{ and } y < 1/2 \\
    -1 & \textup{if } x > 1/2 \textup{ and } y > 1/2 .
    \end{cases}
\end{equation}

\noindent For our simulations we extend the domain to $\mathbf{x} \in [0,2]^2$, symmetrically extend the initial data about $x = 1$ and $y = 1$, and use periodic boundary conditions. The simulations are run using a time step of $2 \times 10^{-4}$ to a final time of $t_f = 0.5$, where the solution is compared against the exact solution found in \cite{guermond2011entropy}. 

\begin{figure}
\centering
\subfloat[Solution]{\label{sfig:Burgers_128p6_soln}\includegraphics[width=.5\textwidth]{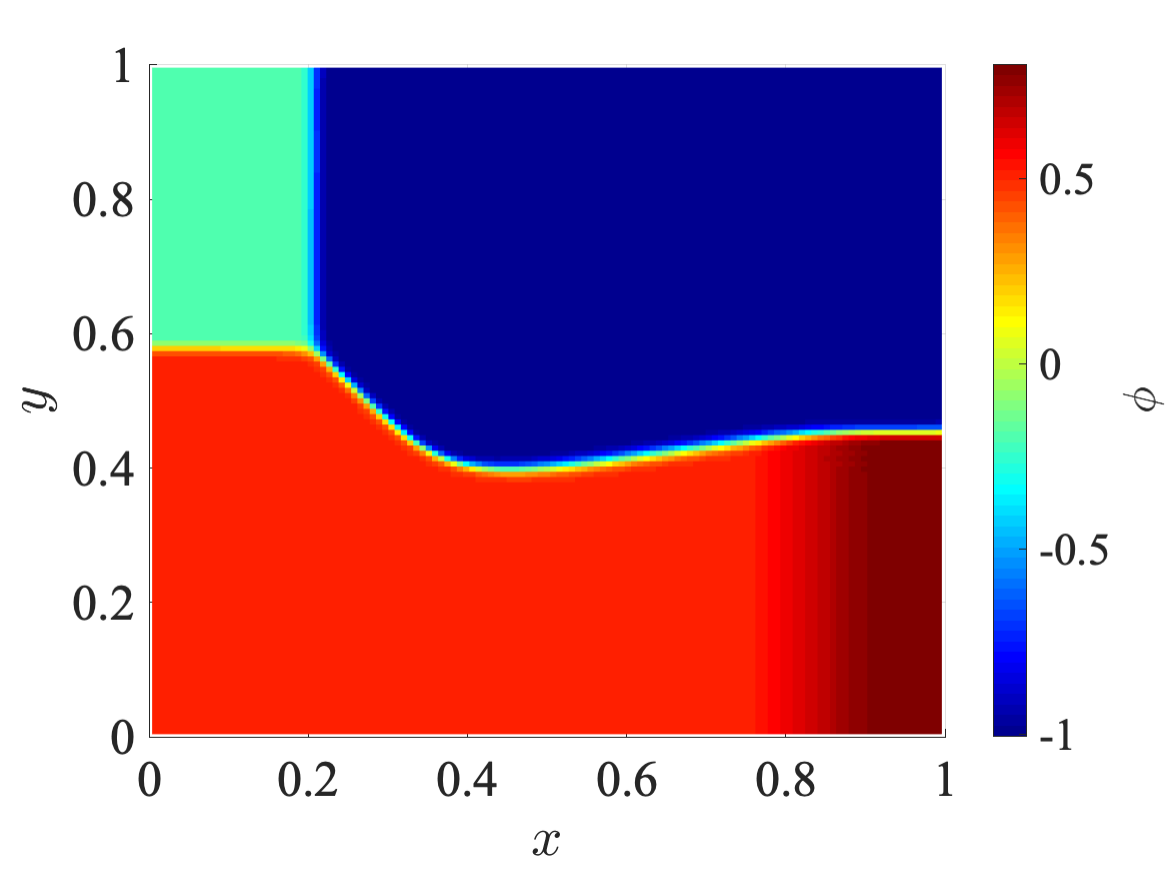}}\hfill
\subfloat[Residual-based viscosity]{\label{sfig:Burgers_128p6_nu}\includegraphics[width=.5\textwidth]{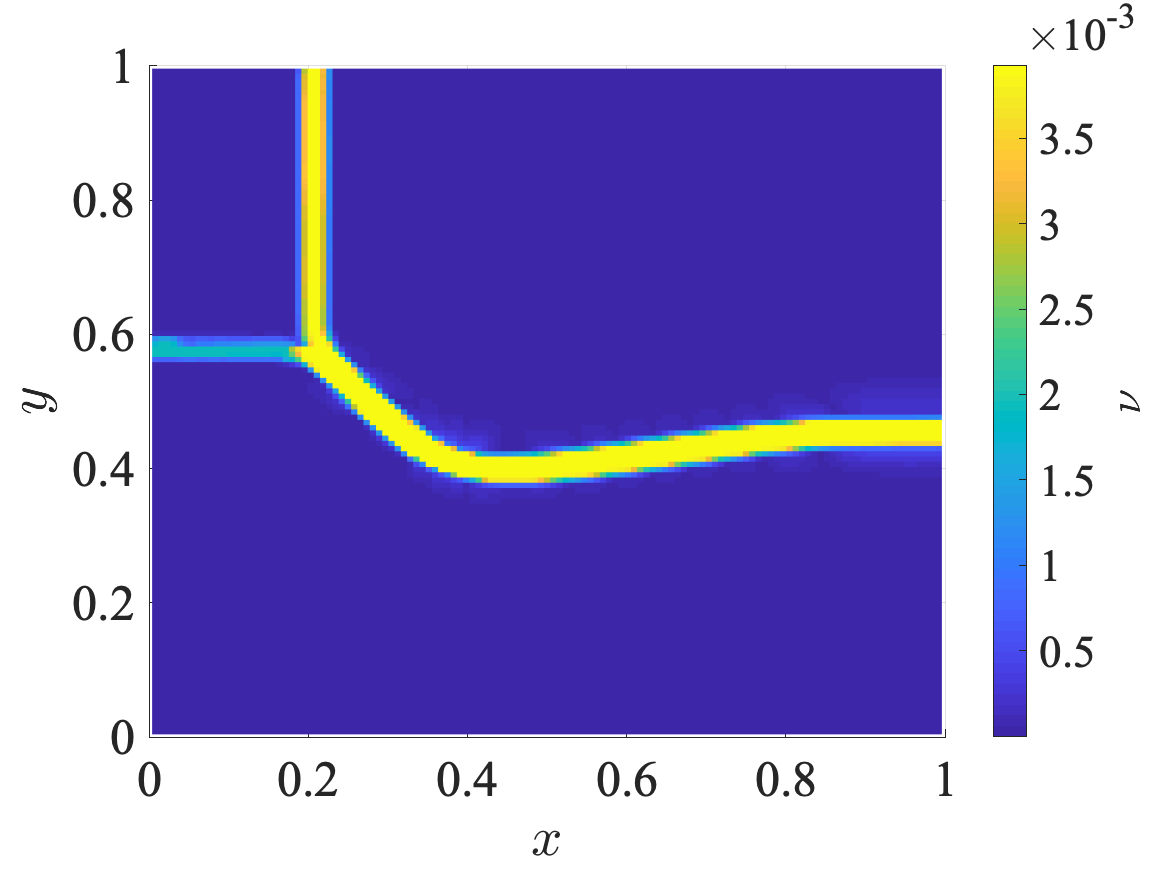}} \\
\caption{Solution to 2D Burgers Riemann problem with $128^2$ elements and $k=6$}
\label{fig:Burgers_128p6}
\end{figure}

\begin{figure}
\centering
\subfloat[$L^2$ Error]{\label{sfig:Burgers_2D_l2}\includegraphics[width=.5\textwidth]{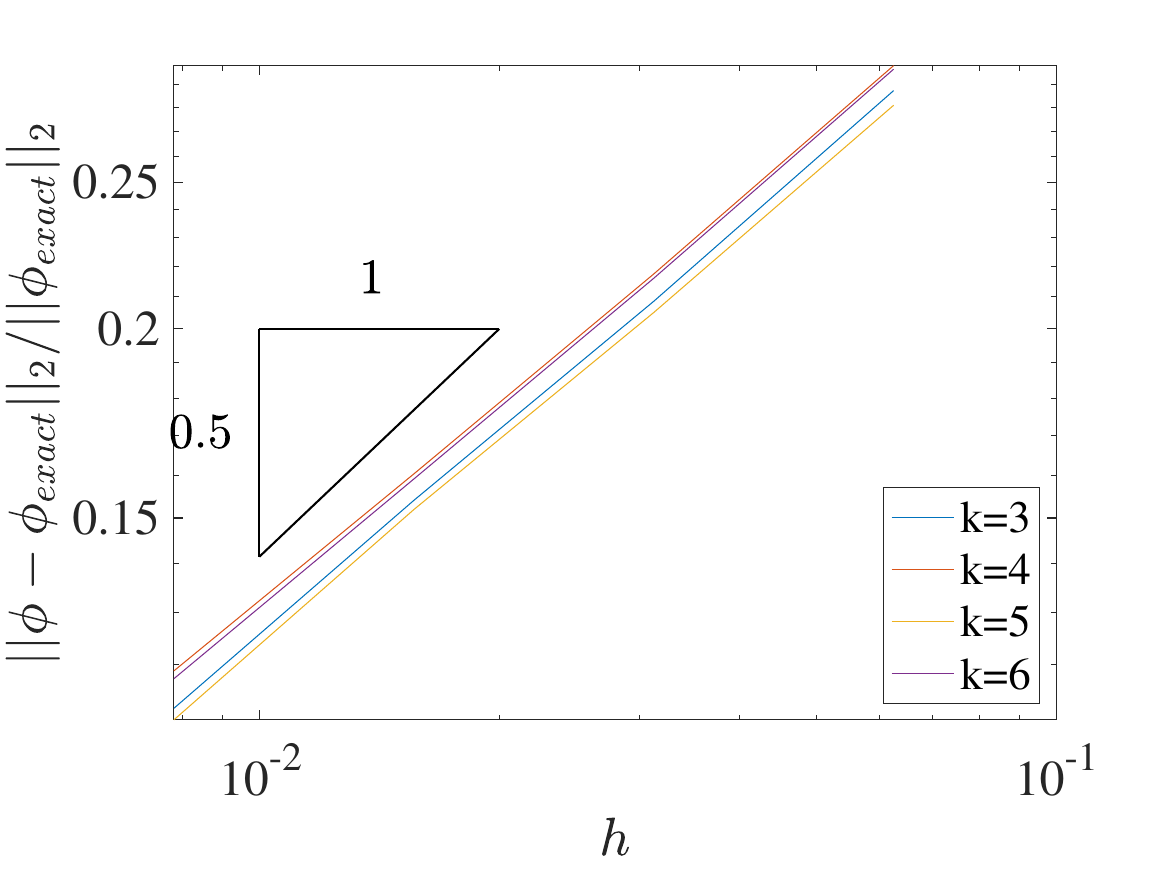}}\hfill
\subfloat[$L^1$ Error]{\label{sfig:Burgers_2D_l1}\includegraphics[width=.5\textwidth]{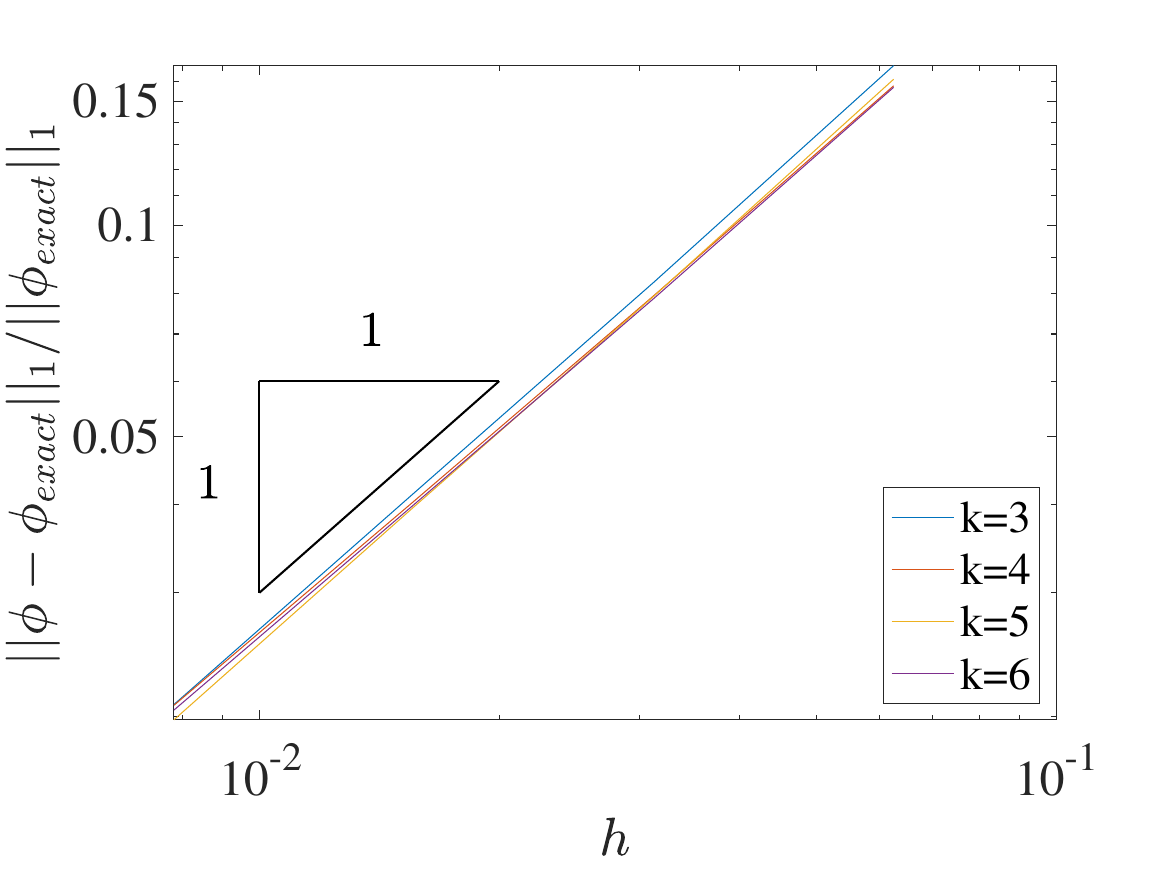}}\\
\caption{Convergence rates of 2D Burgers equation Riemann problem}
\label{fig:Burgers_2D_conv}
\end{figure}

Figure \ref{fig:Burgers_128p6} shows a sample solution at the final time computed using $128^2$ elements and $k = 6$ in the domain $\mathbf{x} \in [0,1]^2$ and the corresponding residual-based viscosity field. Figure \ref{fig:Burgers_2D_conv} shows the convergence of the $L^2$ and $L^1$ errors as the mesh is refined. The numerical solution and resulting errors match well with those produced in \cite{guermond2008entropySpectral, stiernstrom2021SBP, tominec2023RBF} using a variety of residual-based viscosity stabilized discretizations, and the residual-based viscosity is properly focused in the shock region.

\subsection{Buckley-Leverett Equations}

In contrast to the previous examples, conservation laws defined by non-convex flux functions can challenge numerical schemes due to the composite wave nature of their solutions, resulting in connected shocks and rarefactions. Many standard finite volume techniques can fail to converge to the correct entropy solution in these cases \cite{kurganov2007adaptive}. The Buckely-Leverett equations can lead to one such type of conservation law. For a full description of these equations, we refer to \cite{gerritsen2005modeling, christov2008FV}. These equations describe the the flow of 2 immiscible fluids (historically oil and water) within a porous media, and are given by an elliptic equation describing the pore pressure, Darcy's law relating pressure to fluid velocity, and a hyperbolic equation relating the saturation of each phase to the fluid velocity. 

A common technique to solve these equations is known as Implicit Pressure Explicit Saturation (IMPES) \cite{coats2000impes}, where the pressure (and thus velocity) is computed via an implicit solver and the saturation is updated via an explicit scheme. If we assume there are no sources in the domain and the relevant material properties are all set to 1 as in \cite{christov2008FV}, the 1D conservation law to be solved for saturation is given by

\begin{equation}
    \frac{\partial \phi}{\partial t} + \frac{\partial}{ \partial x} \left(\frac{\phi^2}{\phi^2 + (1-\phi)^2} \right) = 0,
\end{equation}

\noindent where $\phi$ is the water saturation field (or the saturation field of the equivalent phase).

An even more complex conservation law results when gravitational forces are included in the formulation. For a 2D problem with gravity acting in the negative $y$ direction and the same set of parameters as \cite{christov2008FV, guermond2011entropy}, the conservation law defining saturation is given by 

\begin{equation}
    \frac{\partial \phi}{\partial t} + \nabla \cdot \mathbf{f} = 0,
\end{equation}

\noindent where $\mathbf{f} = \left( \frac{\phi^2}{\phi^2 + (1-\phi)^2}, \frac{\phi^2 (1-5(1-\phi)^2)}{\phi^2 + (1-\phi)^2} \right)$. For all of the examples in this section we set $C_{RB} = 4$, $C_{max} = 0.25$, and $C_{lin} = 0.25$. While the results were stable with $C_{max} = 0.5$, we found a reduction necessary to achieve optimal convergence rates for the 1D Riemann problem described next.

\subsubsection{1D Riemann Problem}

For this conservation law we focus on results for non-smooth problems only, starting with a standard 1D Riemann problem. The domain is defined as $x \in [-1, 1]$ with initial data 

\begin{equation}
    \phi(x, 0) = 
    \begin{cases}
    1 & \textup{if } x < 0 \\
    0 & \textup{Otherwise}.
    \end{cases}
\end{equation}

\noindent The solution in this case exhibits a compound structure of both a rarefaction and shock moving to the right. We simulate until a final time of $t_f = 0.25$ using a time step of $5 \times 10^{-5}$ and compare the results with the exact solution computed using the method of Osher \cite{osher1984riemann}. 

\begin{figure}
\centering
\subfloat[$L^2$ Error]{\label{sfig:BL_1D_l2}\includegraphics[width=.5\textwidth]{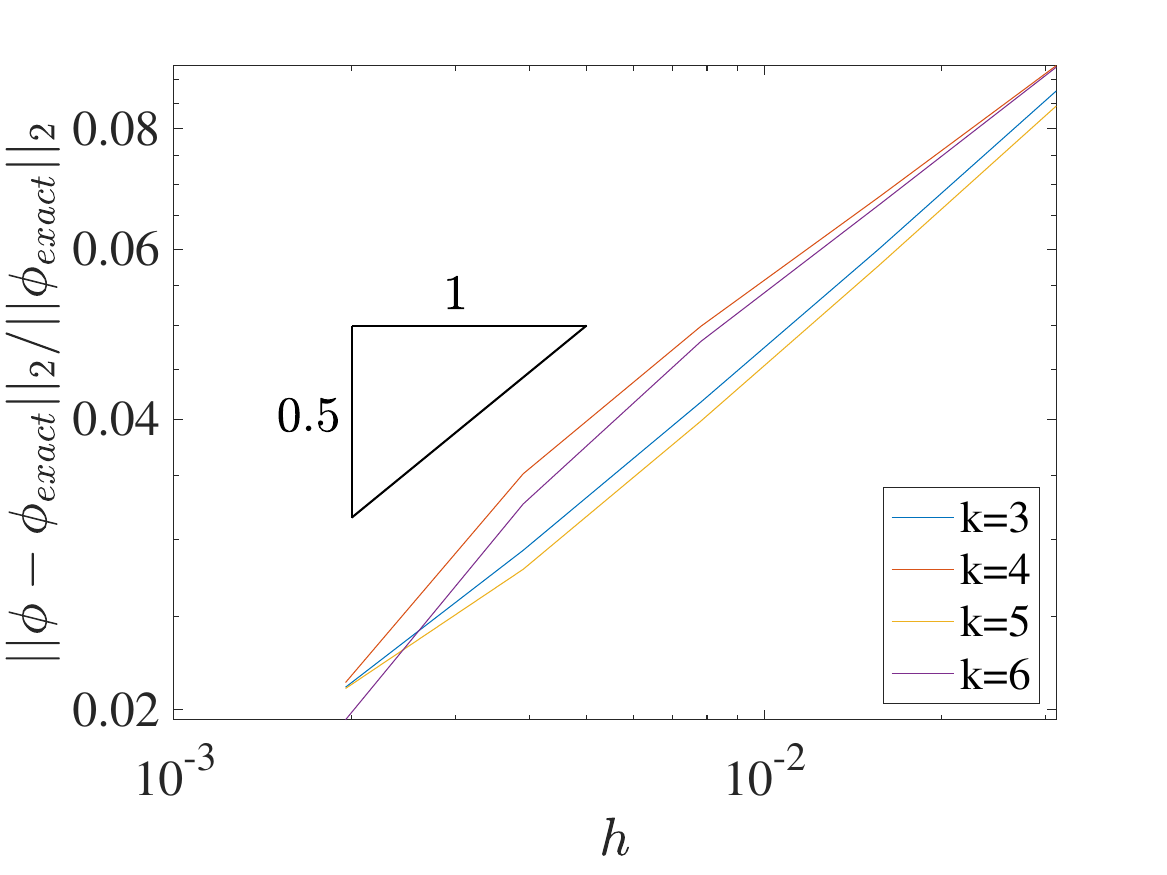}}\hfill
\subfloat[$L^1$ Error]{\label{sfig:BL_1D_l1}\includegraphics[width=.5\textwidth]{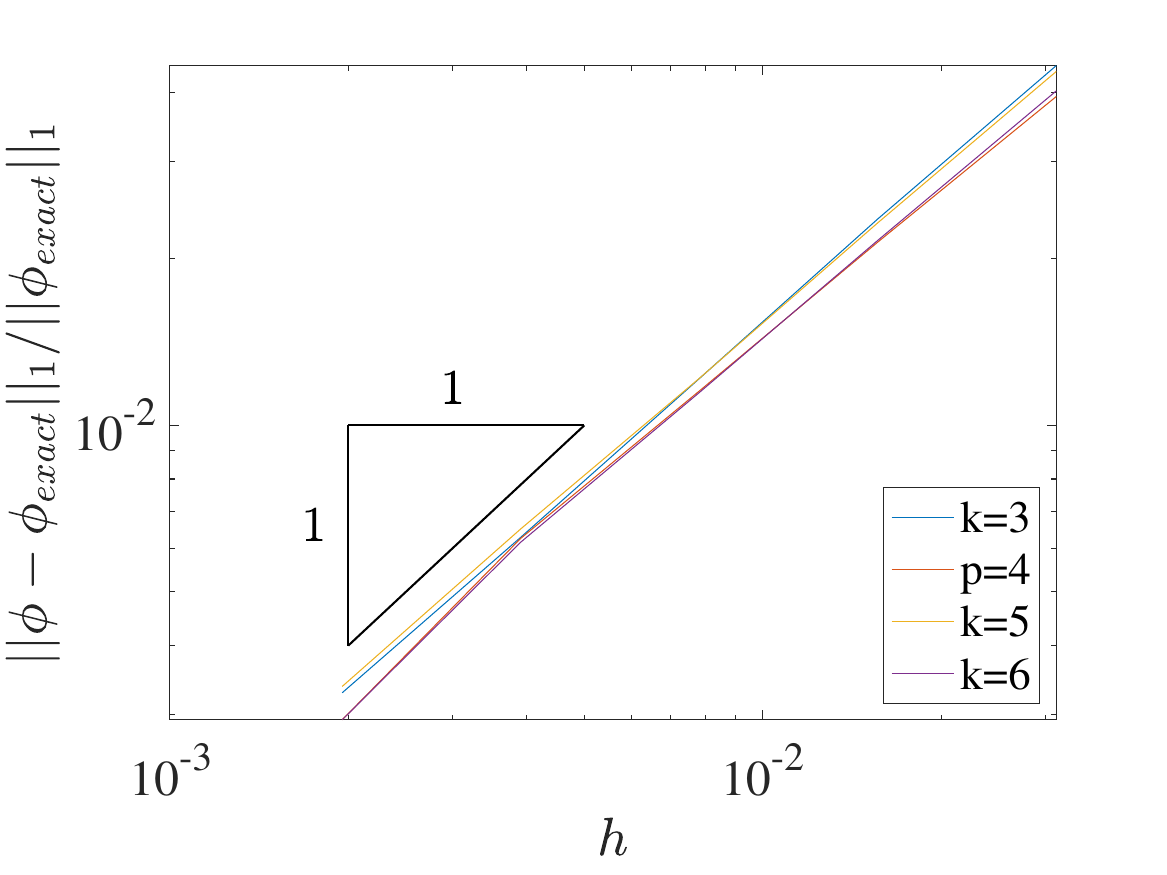}}\\
\caption{Convergence rates of 1D Buckley-Leverett equation Riemann problem}
\label{fig:BL_1D_conv}
\end{figure}

\begin{figure}
\centering
\subfloat[Solution]{\label{sfig:BL_elemComp_phi}\includegraphics[width=.5\textwidth]{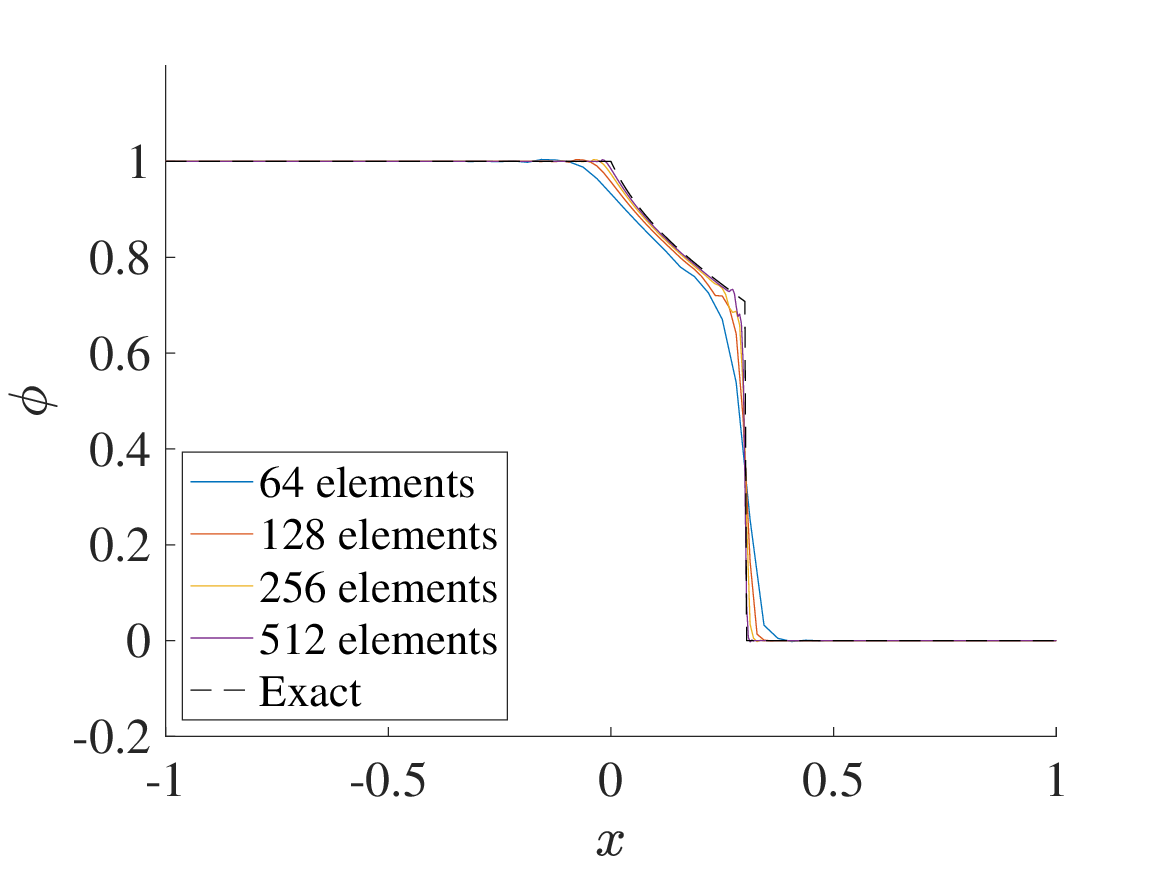}}\hfill
\subfloat[Residual-based viscosity]{\label{sfig:BL_elemComp_nu}\includegraphics[width=.5\textwidth]{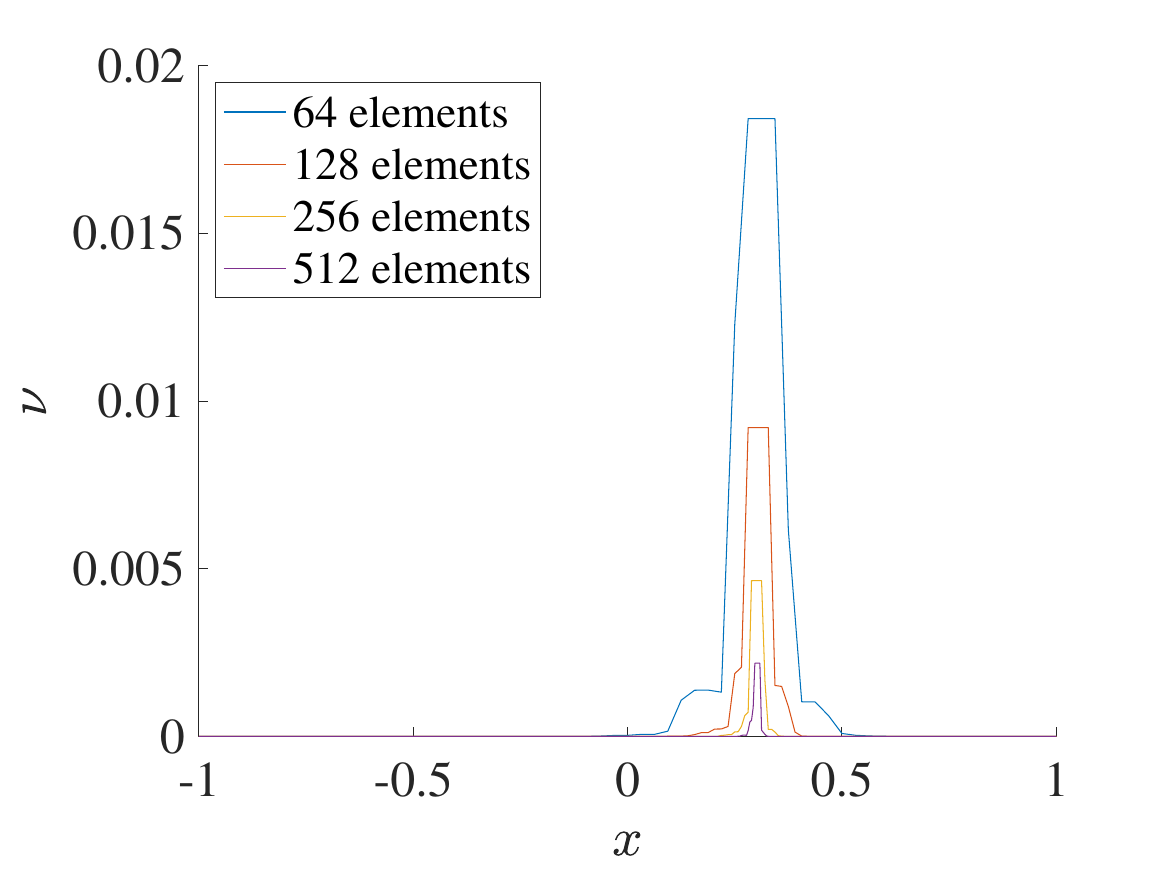}}\\
\caption{1D Buckley-Leverett equation Riemann problem solution with varying number of elements}
\label{fig:BL_elemComp}
\end{figure}

Figure \ref{fig:BL_1D_conv} details the convergence rates of the $L^2$ and $L^1$ errors of the numerical solution at the final time. Once again rates that are close to optimal are recovered in both norms. We also begin to observe some split behavior for the odd and even polynomial degrees as the mesh is refined, like in the previous tests. Figure \ref{fig:BL_elemComp} shows the solutions obtained using $k = 5$ and varying numbers of elements, as well as the corresponding residual-based viscosities. The rarefaction-shock structure is appropriately captured and the viscosity is focused on the shock, though there are some small oscillations that appear where the shock and rarefaction meet. 

\subsubsection{2D Problem with Gravity}

\begin{figure}
\centering
\subfloat[Solution]{\label{sfig:BL_2D_soln}\includegraphics[width=.5\textwidth]{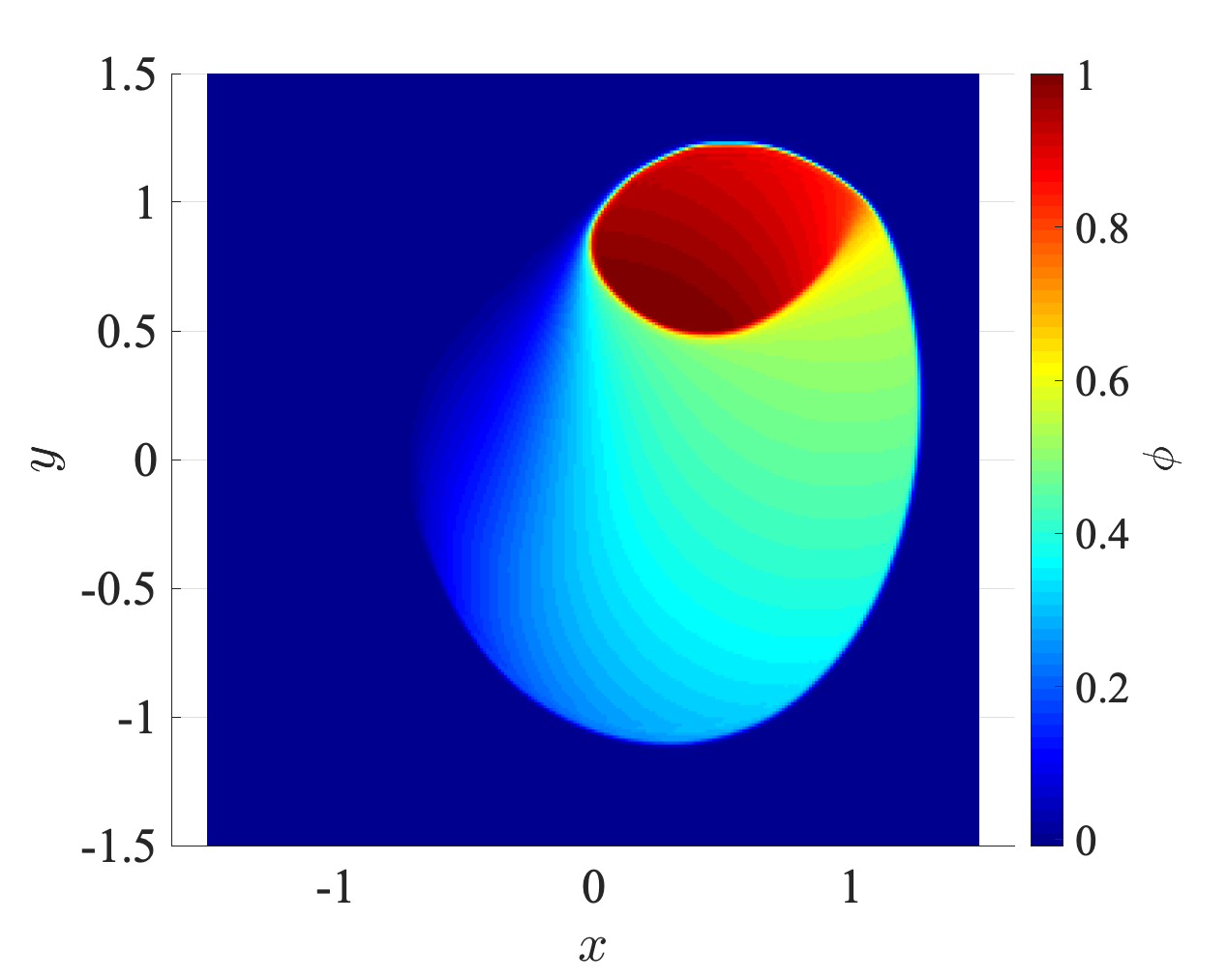}}\hfill
\subfloat[Residual-based viscosity]{\label{sfig:BL_2D_nu}\includegraphics[width=.5\textwidth]{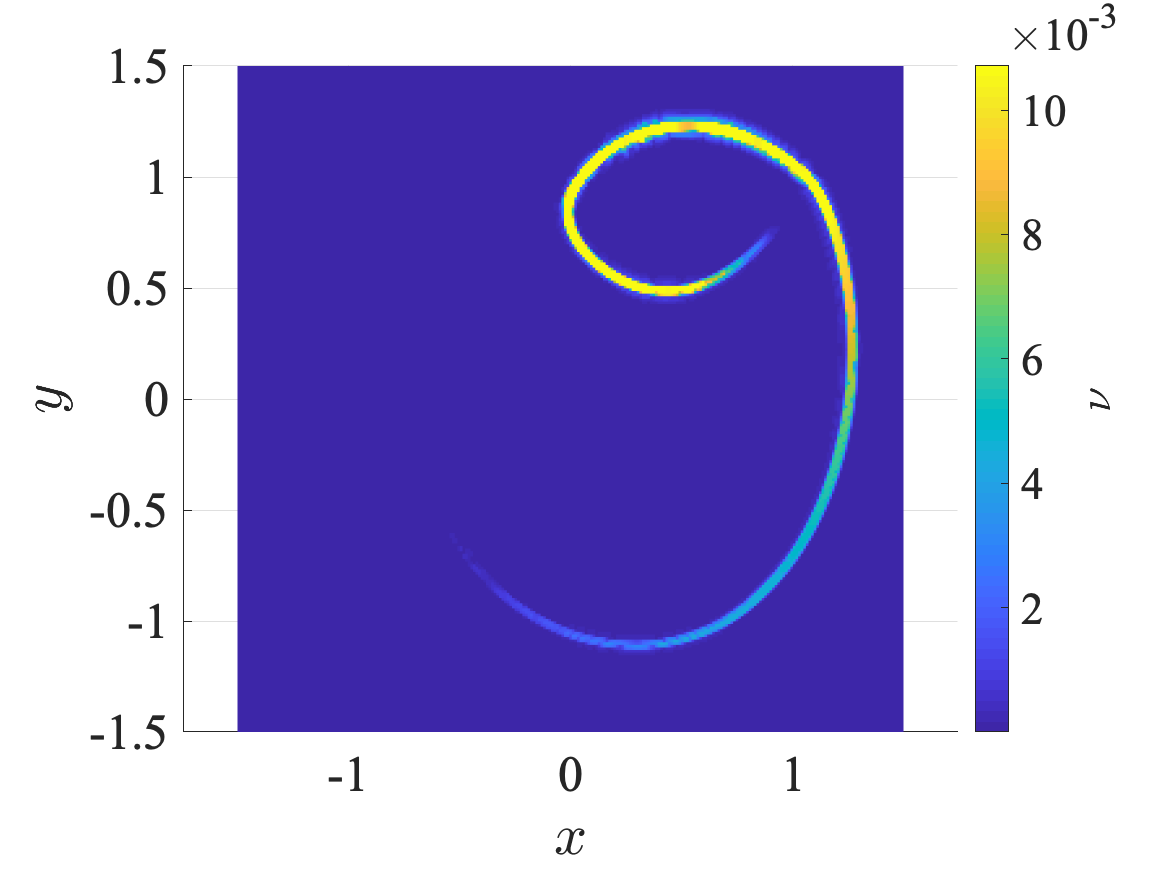}} \\
\caption{Solution to 2D Buckley-Leverett problem with $256^2$ elements and k=4}
\label{fig:BL_2D}
\end{figure}

Next we consider the 2D problem with gravity with initial condition defined by 

\begin{equation}
    \phi(\mathbf{x}, 0) = 
    \begin{cases}
    1 & \textup{if } x^2 + y^2 < 1/2\\
    0 & \textup{Otherwise}.
    \end{cases}
\end{equation}

\noindent on the domain $\mathbf{x} \in [-1.5, 1.5]^2$. There is no exact solution to this problem, but we show a sample numerical solution computed at $t_f = 0.5$ with $256^2$ elements, $k=5$, and a time step of $1 \times 10^{-4}$ in Figure \ref{fig:BL_2D}. The solution agrees well with the ones presented in \cite{christov2008FV, guermond2011entropy} using a limited finite volume technique and entropy viscosity stabilized finite elements. The 2D composite wave structure is well captured by the residual-based viscosity and correctly depicts the variation in strength of the shock. We do note, however, that there is still some oscillatory content in the solution, particularly near the strongest shock location at the top of the high $\phi$ region. This can be removed by increasing $C_{RB}$ or $C_{max}$ to further smooth the solution.

\subsection{Euler Equations}

Our final set of results is obtained for the Euler equations, in which we also show the different effects produced by the Laplacian and Guermond-Popov flux regularizations. For the Laplacian flues, we set $C_{RB} = 4$, $C_{max} = 0.1$, and $C_{lin} = 0.25$. When using the Guermond-Popov fluxes, we have found a little more variation necessary to achieve optimal results. For the 1D Euler equation examples we set $C_{RB} = 4$, $C_{max} = 0.2$, and $\mathcal{P} = 0.5$ while for the 2D case we have found setting $C_{RB} = 4$, $C_{max} = 0.1$, and $\mathcal{P} = 1$ produces a slightly higher quality result. A similar change was also made in \cite{nazarov2017investigation} for this 2D test case. In both one and two dimensions we use $C_{lin} = 0.25$.

\subsubsection{Isentropic Flow}

We begin by considering a problem with a smooth solution for one final verification that the stabilization terms do not adversely affect the accuracy of the scheme. The computational domain is $x \in [-1, 1]$ with initial conditions

\begin{equation}
    \rho(x, 0) = 1 + 0.9\sin(\pi x)
\end{equation}

\begin{equation}
    \rho u(x, 0) = 0
\end{equation}

\begin{equation}
    E(x, 0) = \frac{ \rho^\gamma}{\gamma-1}.
\end{equation}

\noindent We set $\gamma = 3$ and enforce periodic boundary conditions. We simulate to a final time of $t_f = 0.1$ using a time step of size $5 \times 10^{-5}$, and the resulting flow should be isentropic. Thus we should recover the expected convergence rates for smooth problems with the stabilized collocation schemes.  

\begin{figure}
\centering
\subfloat[Density with linear stabilization]{\label{sfig:isen_r_lin}\includegraphics[width=.5\textwidth]{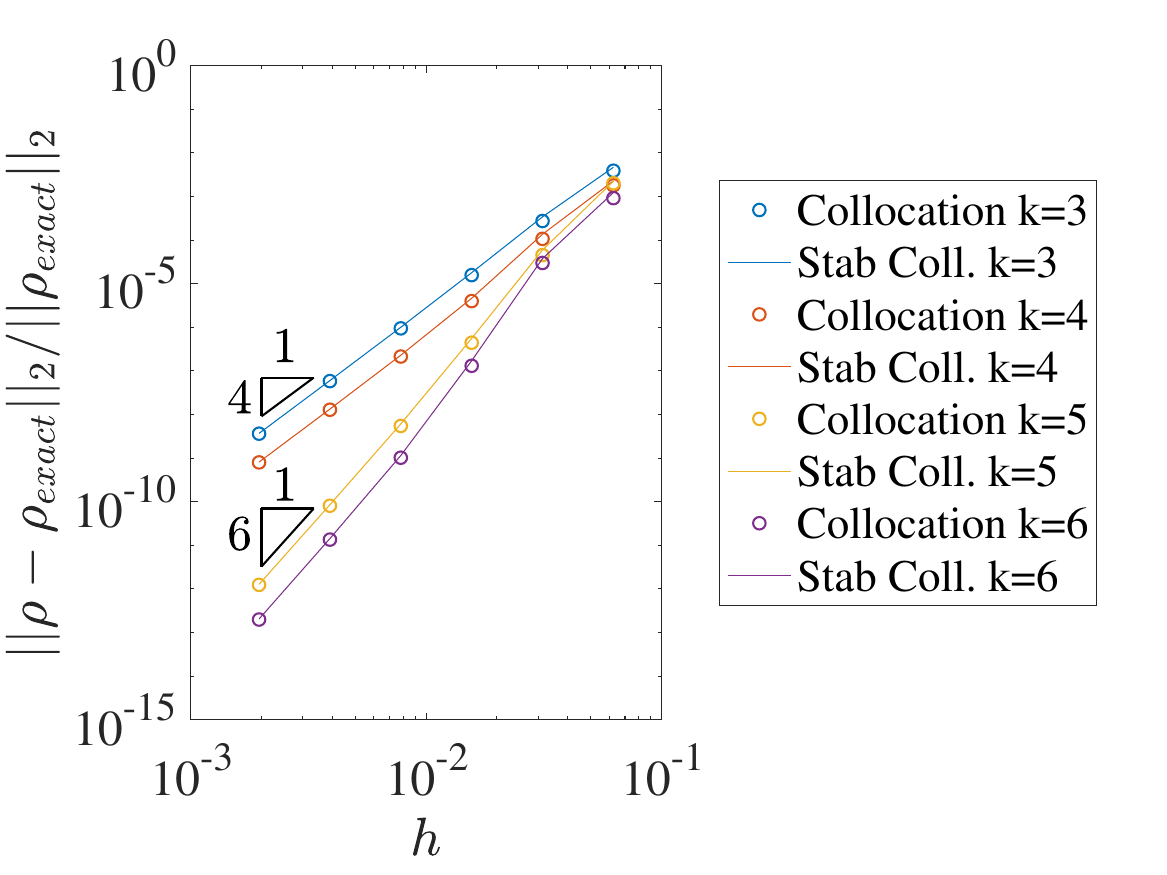}}\hfill
\subfloat[Density with linear and nonlinear stabilization]{\label{sfig:isen_r_stab}\includegraphics[width=.5\textwidth]{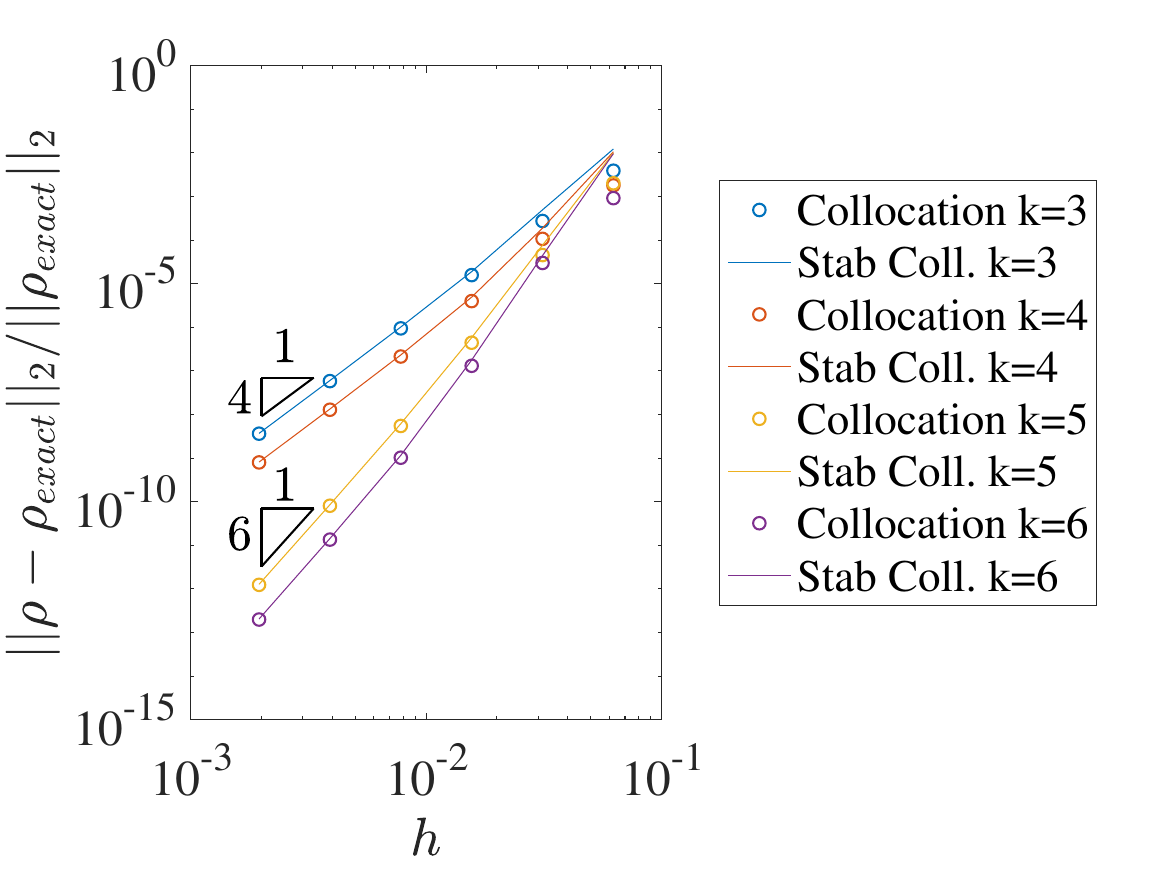}}\\
\subfloat[Momentum with linear stabilization]
{\label{sfig:isen_ru_lin}\includegraphics[width=.5\textwidth]{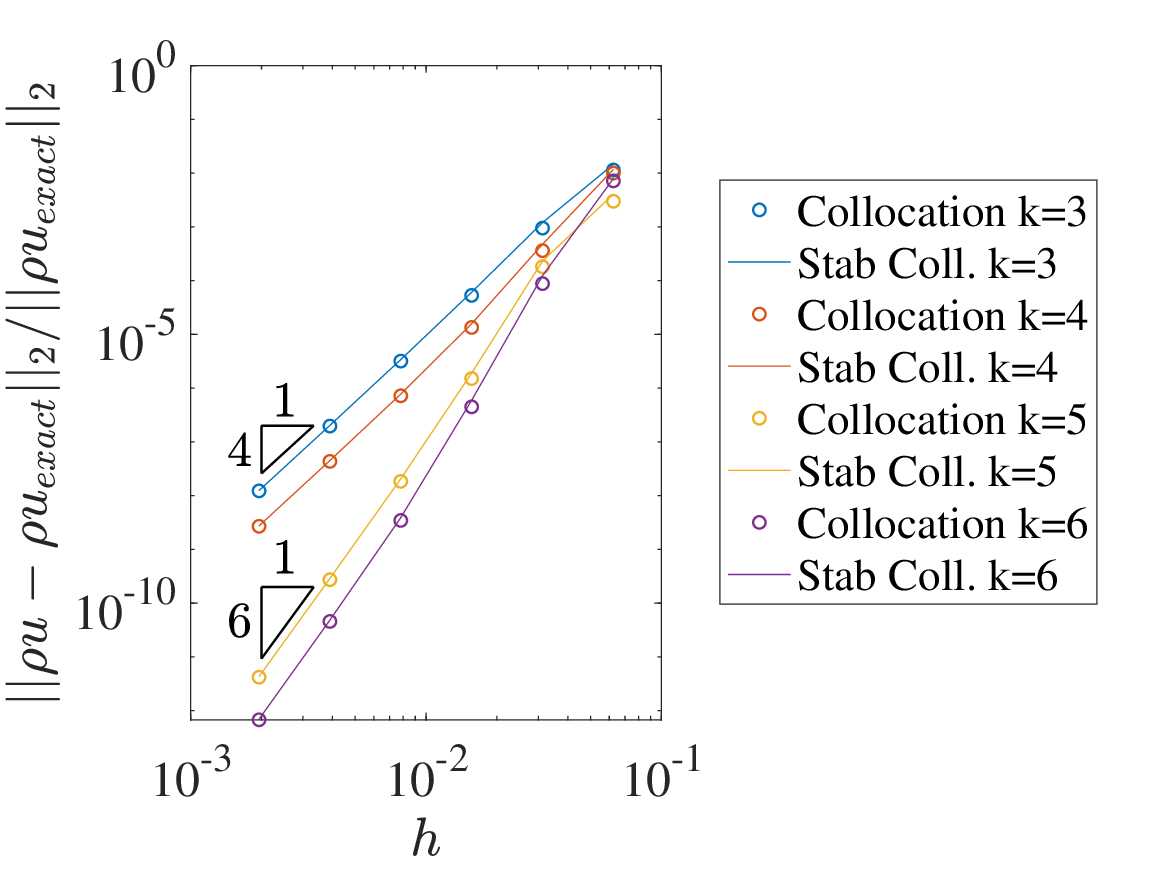}}\hfill
\subfloat[Momentum with linear and nonlinear stabilization]{\label{sfig:isen_ru_stab}\includegraphics[width=.5\textwidth]{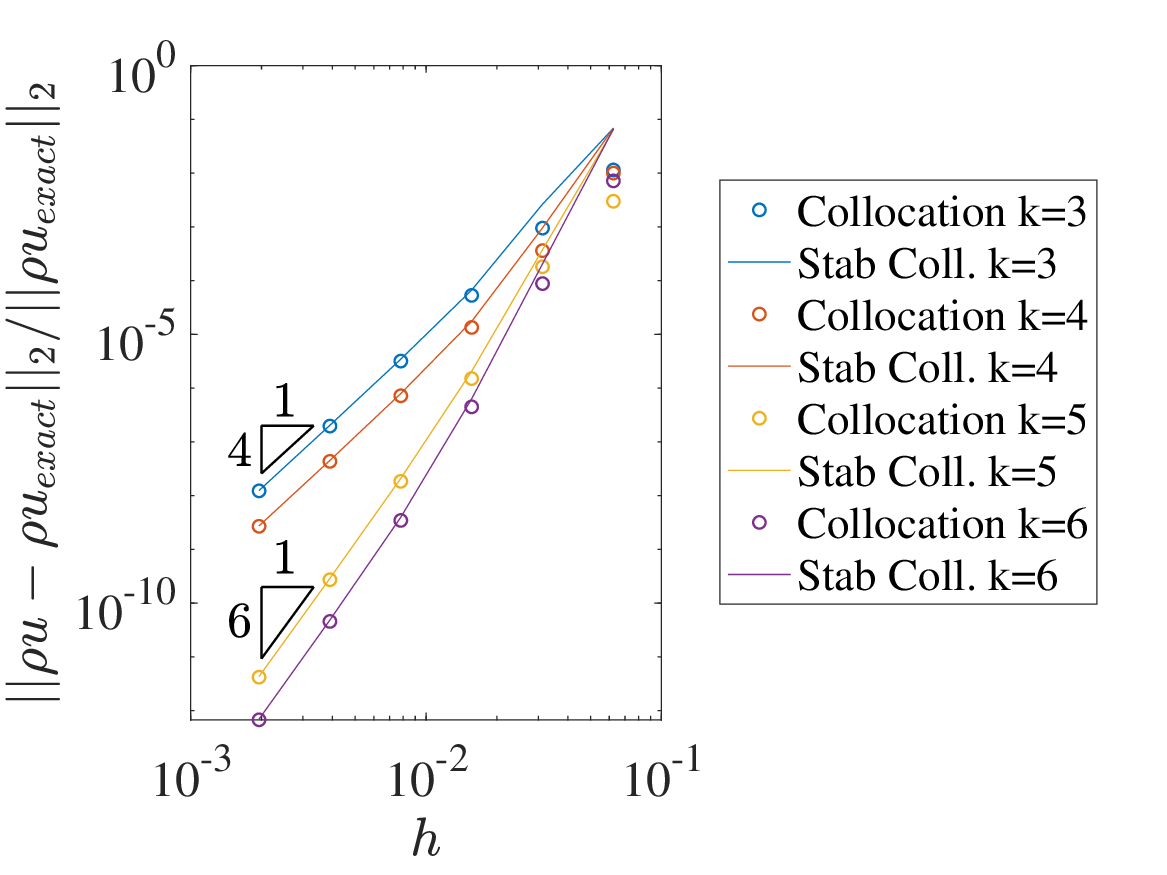}}\\
\subfloat[Energy with linear stabilization]
{\label{sfig:isen_E_lin}\includegraphics[width=.5\textwidth]{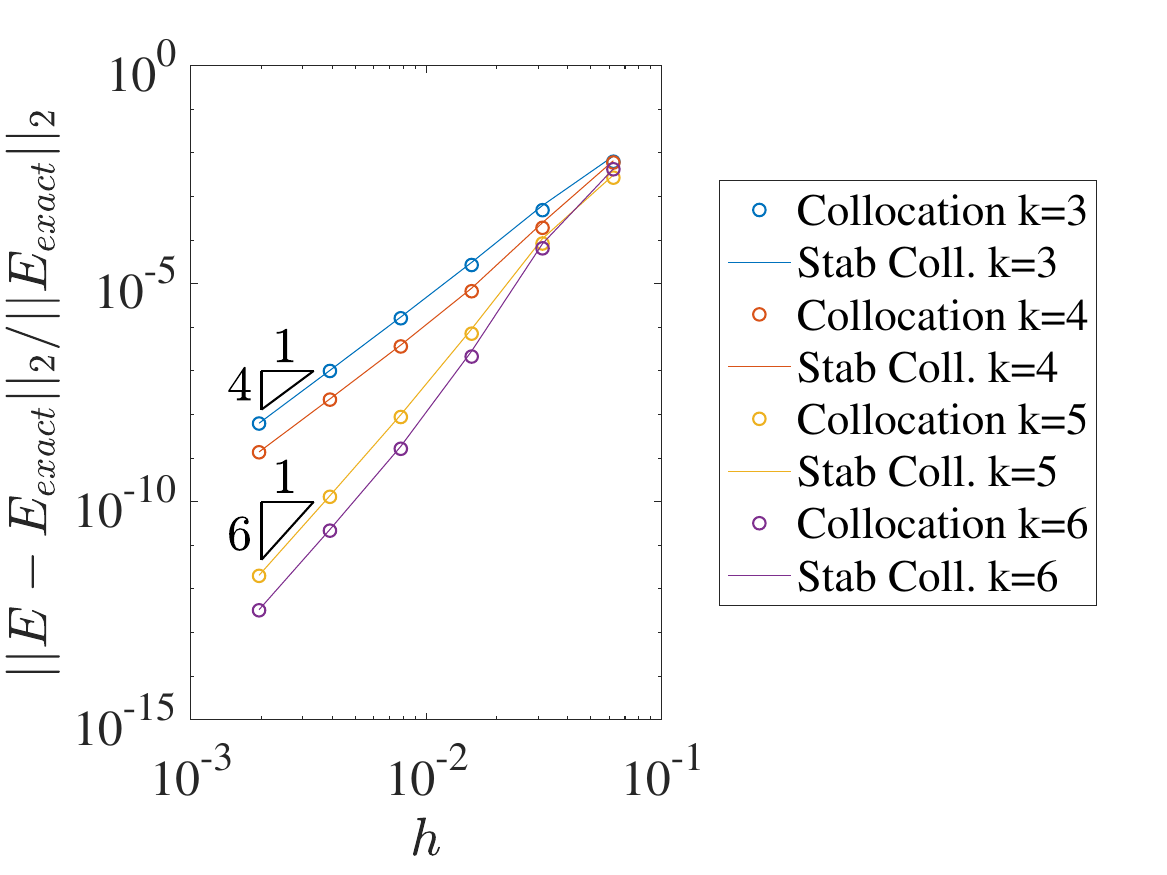}}\hfill
\subfloat[Energy with linear and nonlinear stabilization]{\label{sfig:isen_E_stab}\includegraphics[width=.5\textwidth]{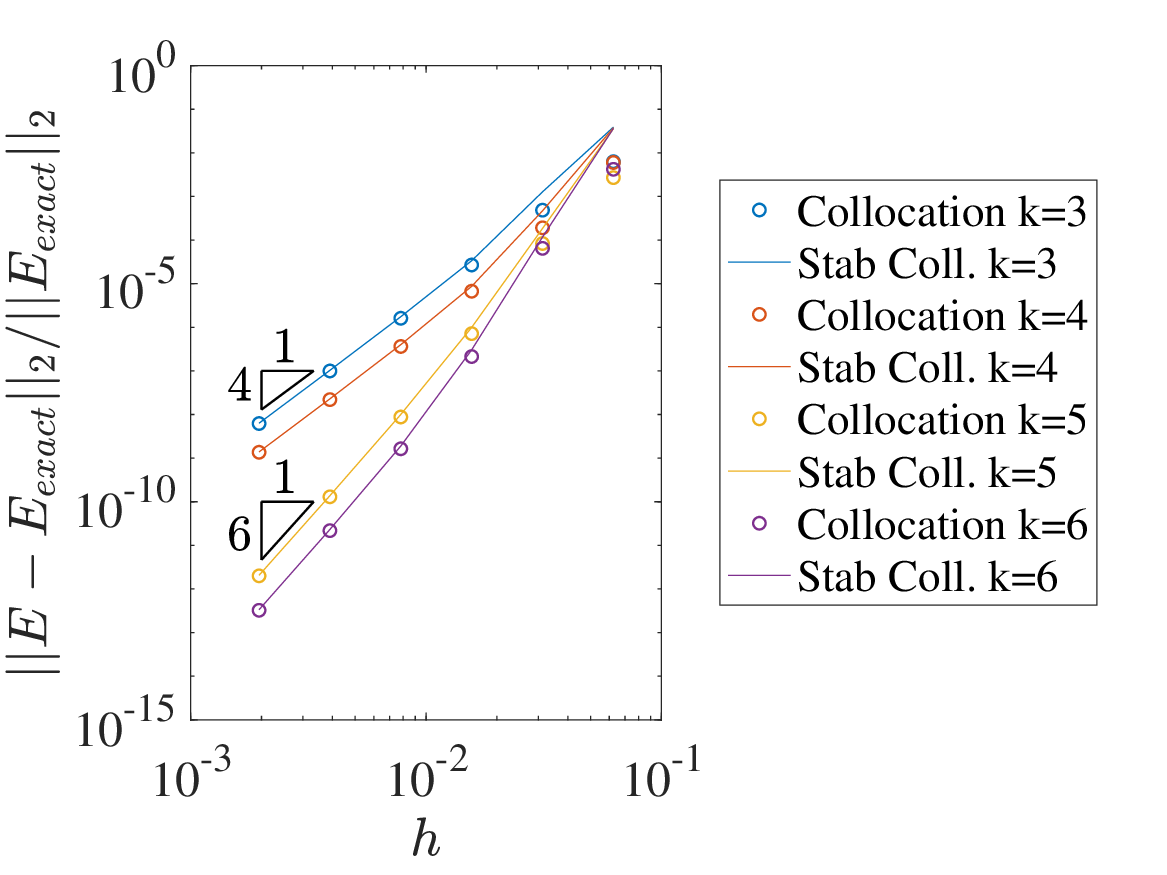}}\\
\caption{Convergence rates of isentropic flow with Guermond-Popov fluxes}
\label{fig:isen_conv}
\end{figure}

For brevity, we only present results obtained using Guermond-Popov fluxes, the results obtained with the Laplacian regularization are very similar. Figure \ref{fig:isen_conv} details the $L^2$ errors obtained for each of the conserved quantities with unstabilized collocation, collocation with only linear stabilization, and collocation with linear and nonlinear stabilization. These results behave much like the corresponding linear advection and Burgers results. For all three conserved quantities the expected rates of $k$ and $k+1$ for even and odd degrees are recovered, and only on the coarsest meshes does the residual-based viscosity noticeably increase the magnitude of the errors. 


\subsubsection{Sod Shock Tube}

Next we consider the standard Sod shock tube Riemann problem, with initial conditions given by

\begin{equation}
    \rho(x, 0) = 
    \begin{cases}
    1 & \textup{if } x < 1/2\\
    0.125 & \textup{Otherwise}.
    \end{cases}
\end{equation}

\begin{equation}
    \rho u(x, 0) = 0
\end{equation}

\begin{equation}
    E(x, 0) = 
    \begin{cases}
    \frac{1}{\gamma -1} & \textup{if } x < 1/2\\
    \frac{0.1}{\gamma - 1}& \textup{Otherwise}.
    \end{cases}
\end{equation}

\noindent over the domain $x \in [0,1]$ with $\gamma = 1.4$. We integrate in time until $t_f = 0.25$ using a time step of size $1 \times 10^{-4}$, during which the initial discontinuity separates into a rarefaction, contact wave, and shock according to the characteristics of the system. 

\begin{figure}
\centering
\subfloat[Density $L^2$ Error]{\label{sfig:sod_r_l2_lap}\includegraphics[width=.5\textwidth]{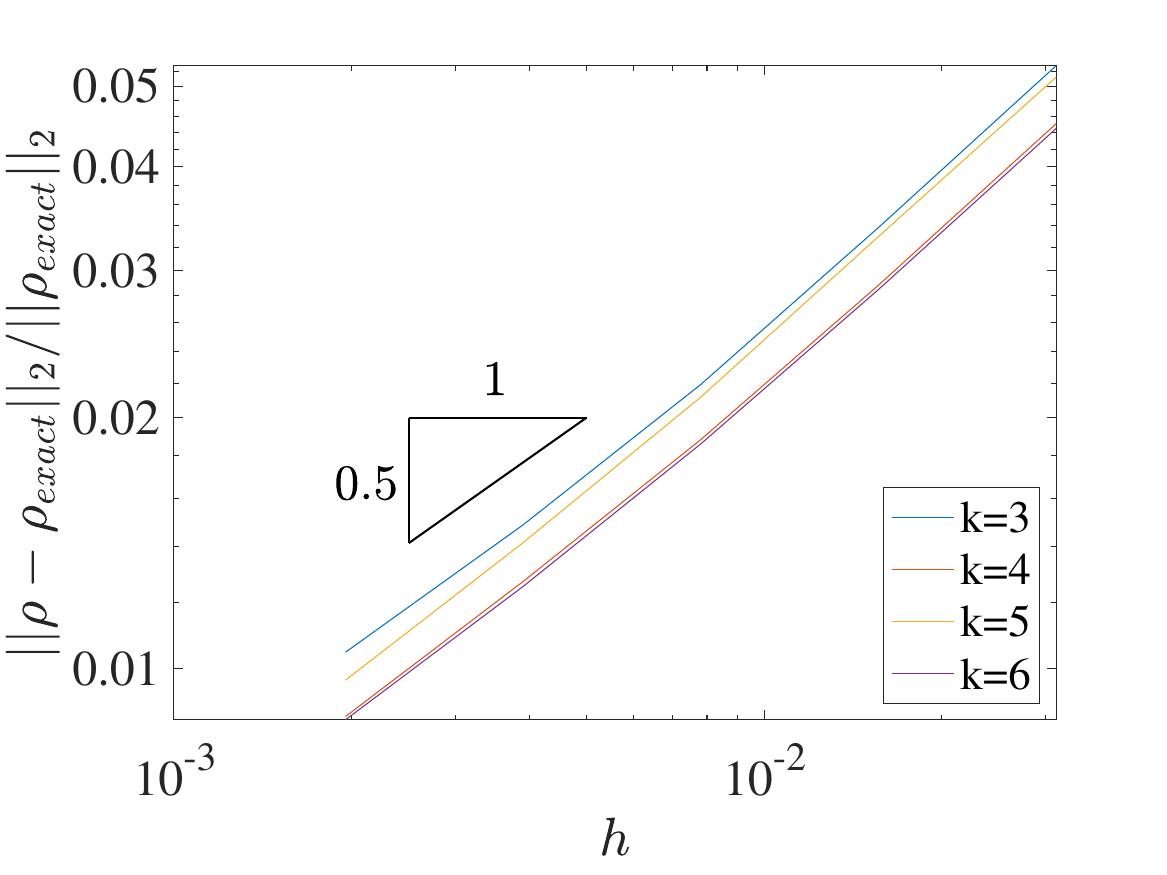}}\hfill
\subfloat[Density $L^1$ Error]{\label{sfig:sod_r_l1_lap}\includegraphics[width=.5\textwidth]{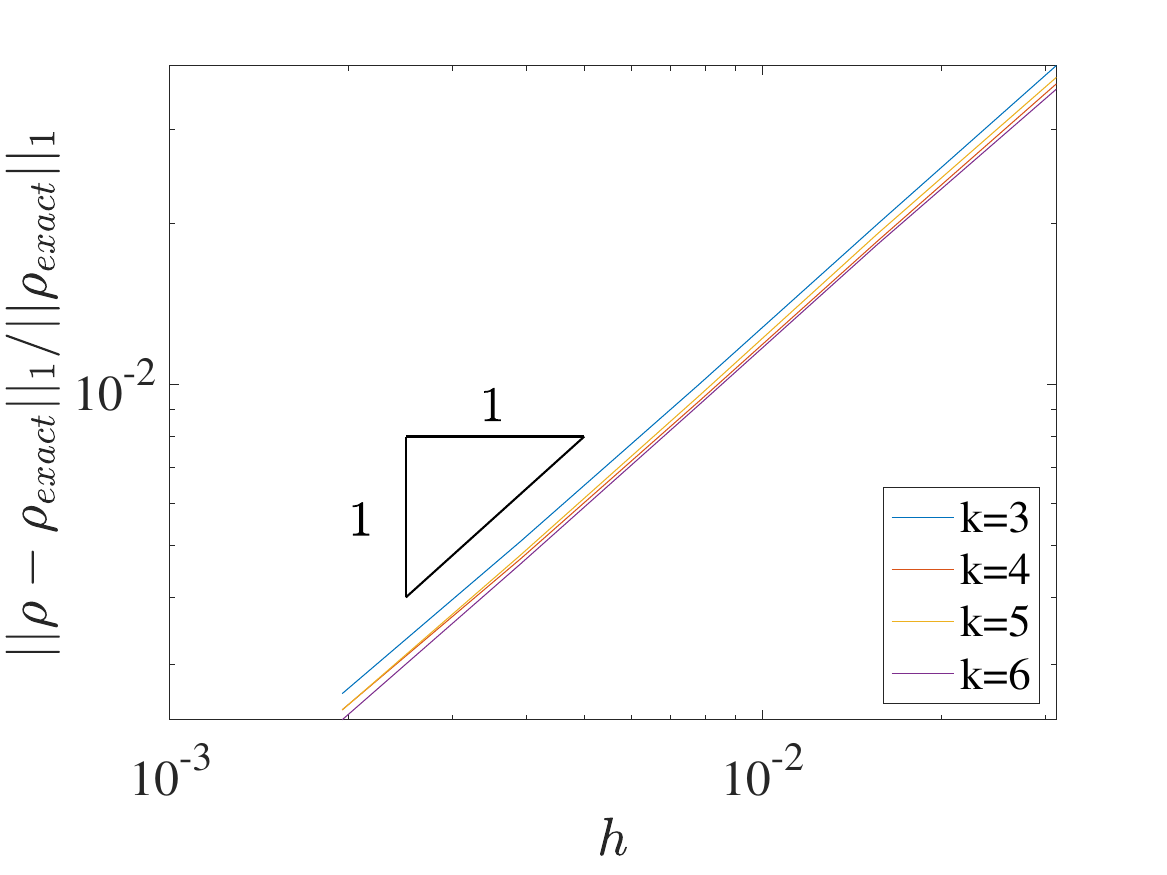}}\\
\subfloat[Momentum $L^2$ Error]
{\label{sfig:sod_ru_l2_lap}\includegraphics[width=.5\textwidth]{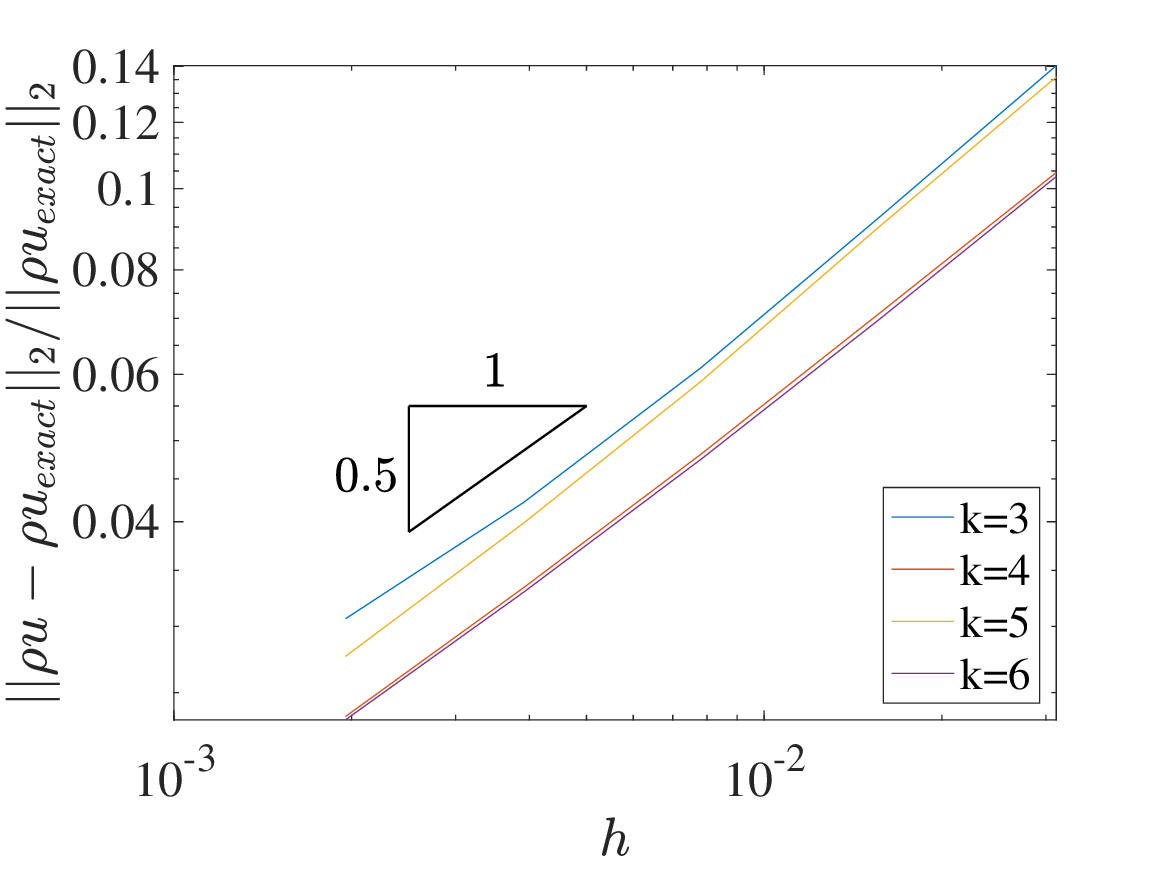}}\hfill
\subfloat[Momentum $L^1$ Error]{\label{sfig:sod_ru_l1_lap}\includegraphics[width=.5\textwidth]{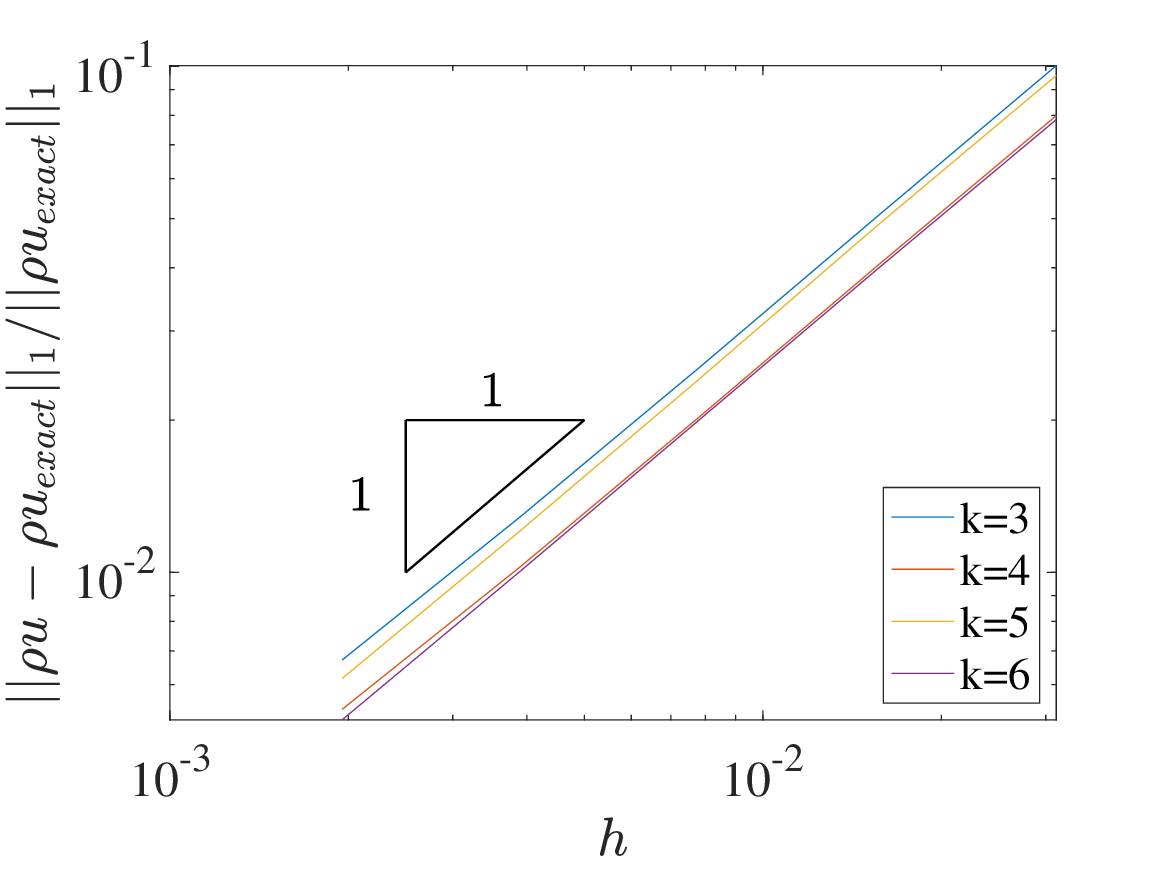}}\\
\subfloat[Energy $L^2$ Error]
{\label{sfig:sod_E_l2_lap}\includegraphics[width=.5\textwidth]{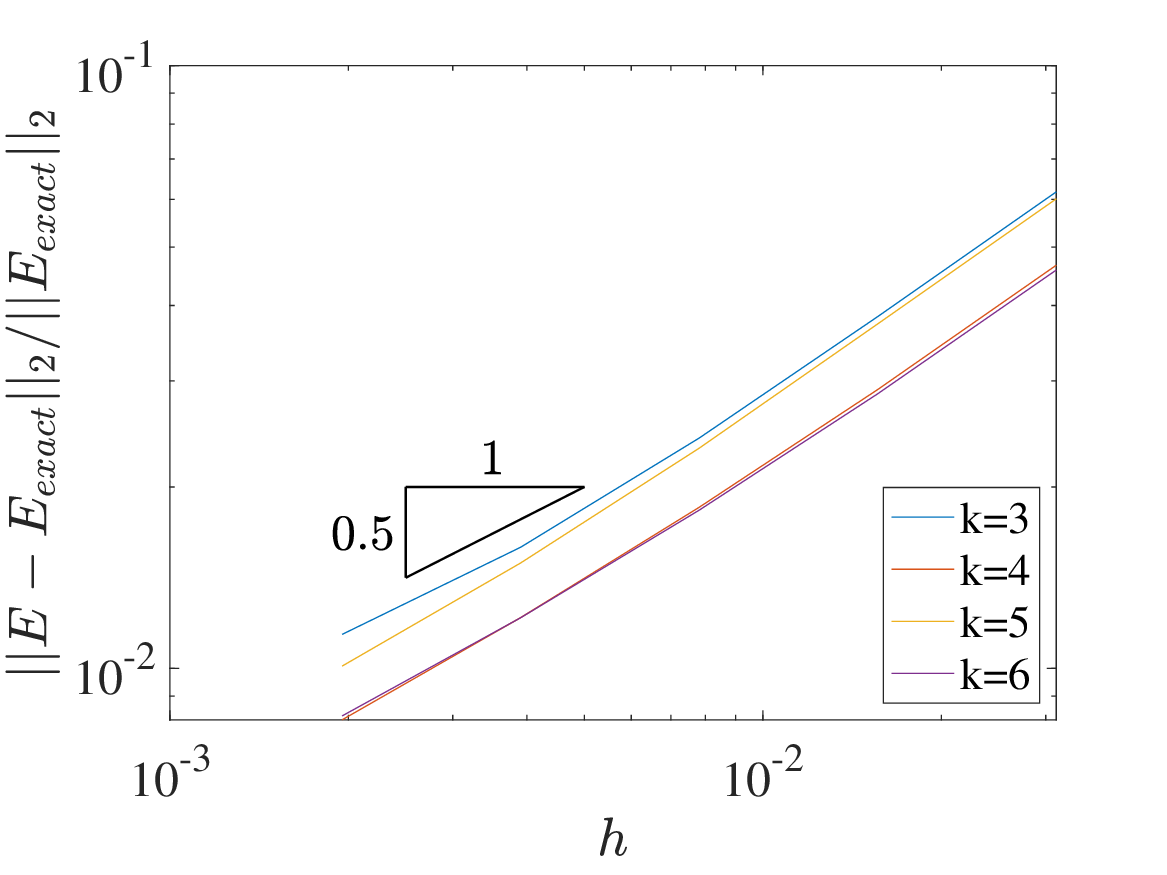}}\hfill
\subfloat[Energy $L^1$ Error]{\label{sfig:sod_E_l1_lap}\includegraphics[width=.5\textwidth]{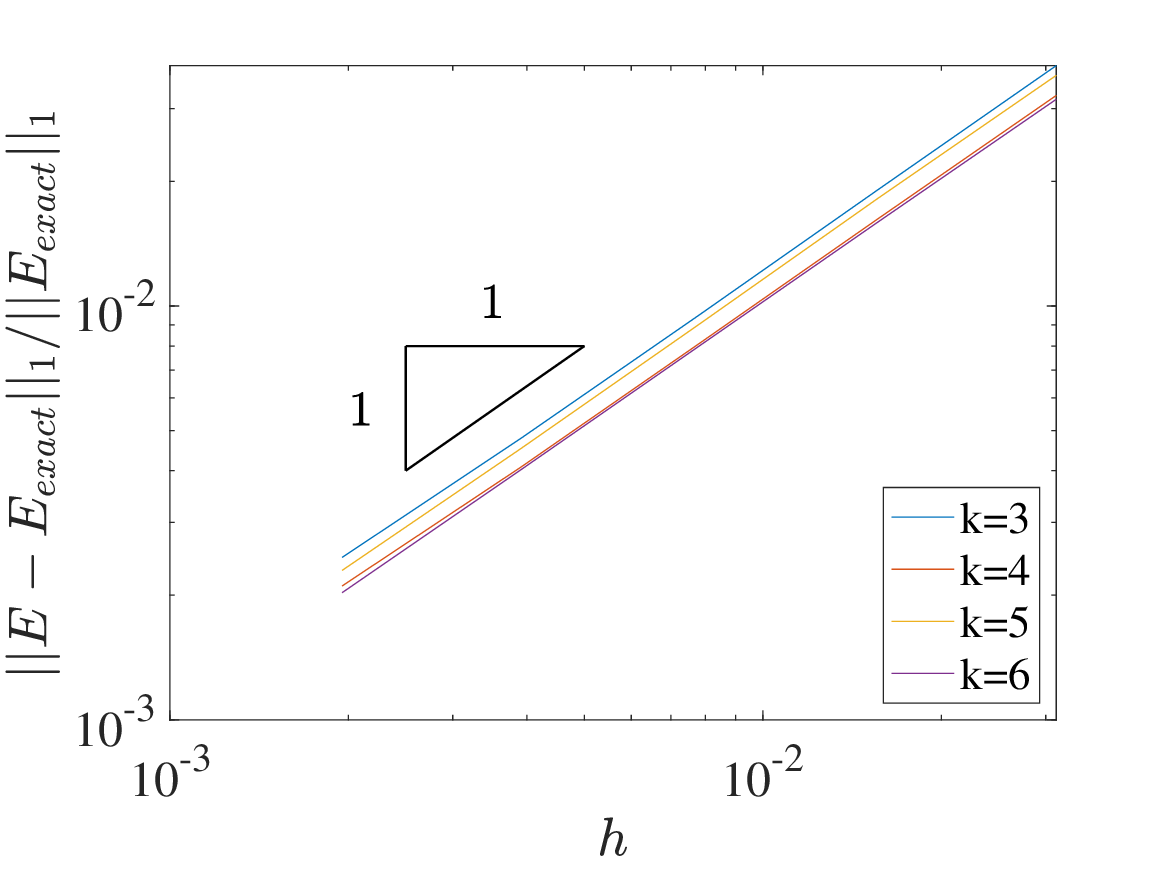}}\\
\caption{Convergence rates of Sod shock tube with Laplacian fluxes}
\label{fig:sod_conv_lap}
\end{figure}

\begin{figure}
\centering
\subfloat[Density $L^2$ Error]{\label{sfig:sod_r_l2_GP}\includegraphics[width=.5\textwidth]{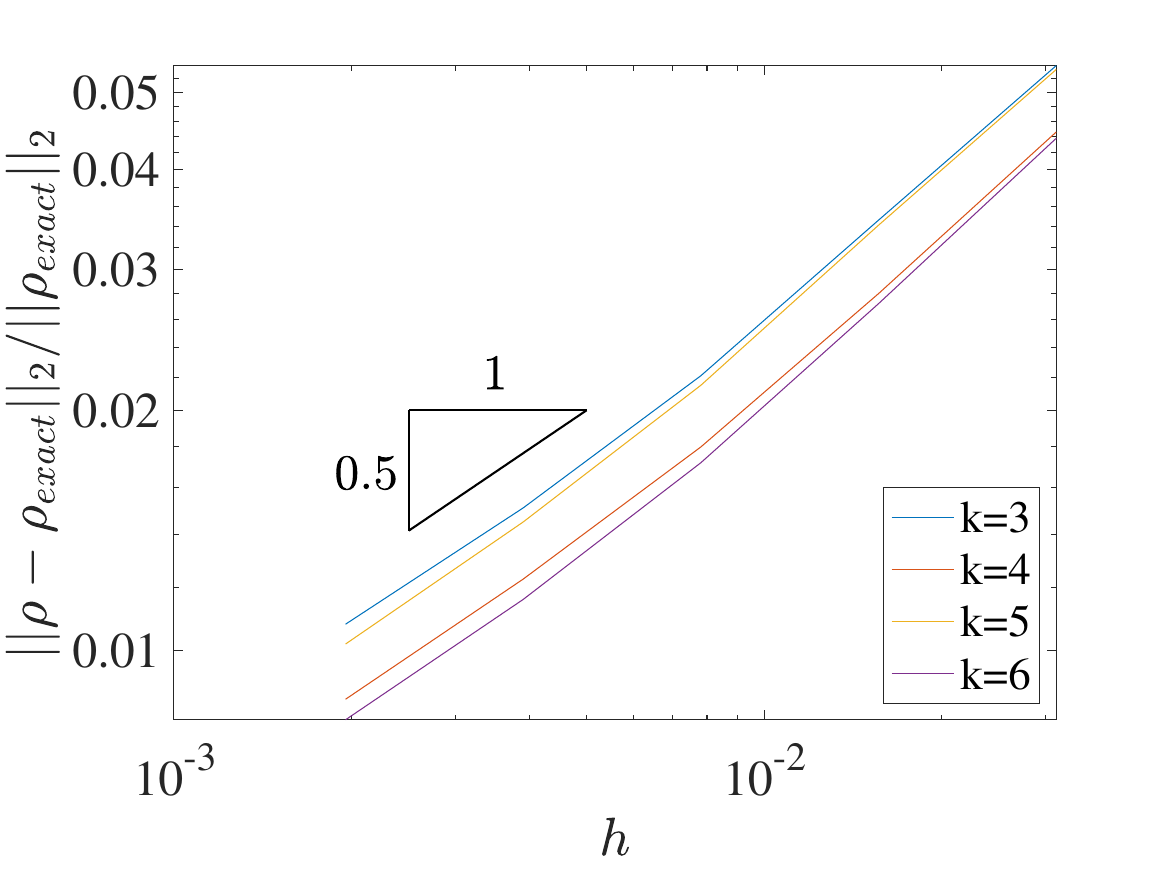}}\hfill
\subfloat[Density $L^1$ Error]{\label{sfig:sod_r_l1_GP}\includegraphics[width=.5\textwidth]{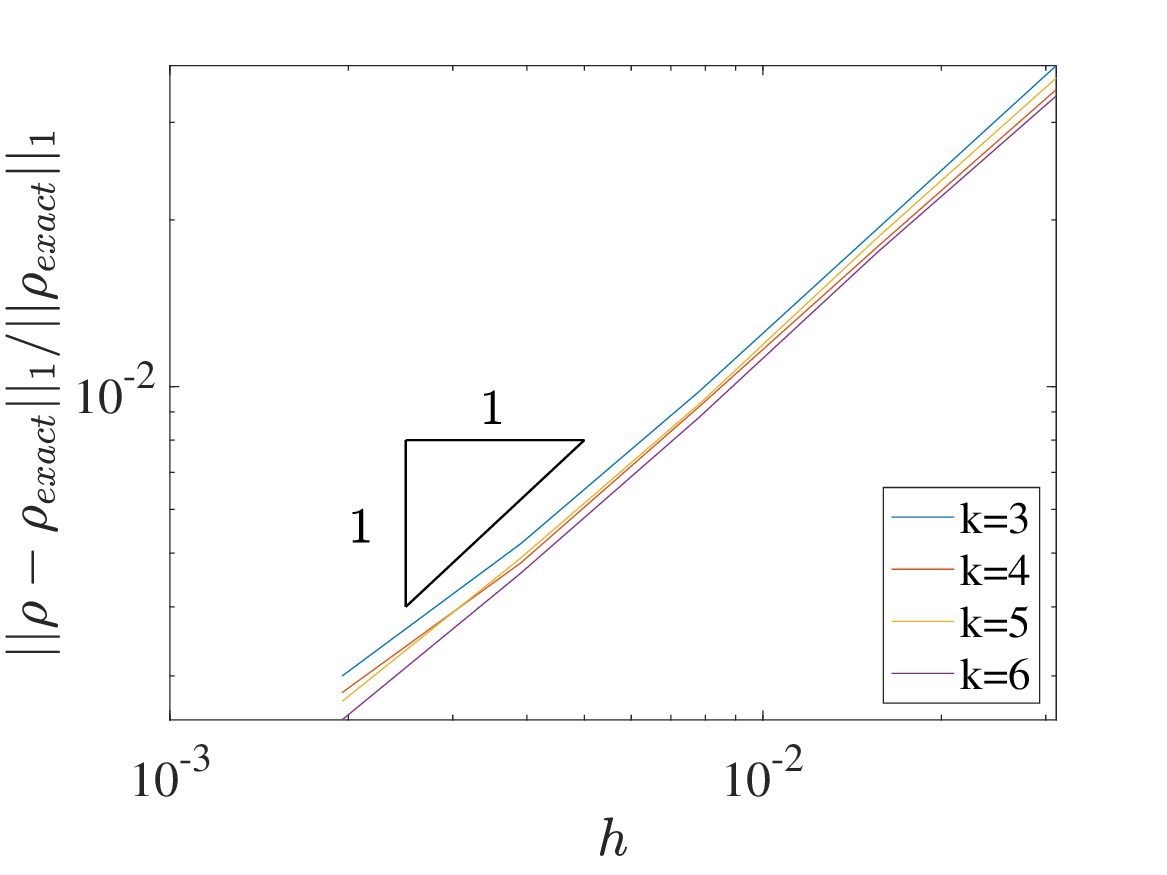}}\\
\subfloat[Momentum $L^2$ Error]
{\label{sfig:sod_ru_l2_GP}\includegraphics[width=.5\textwidth]{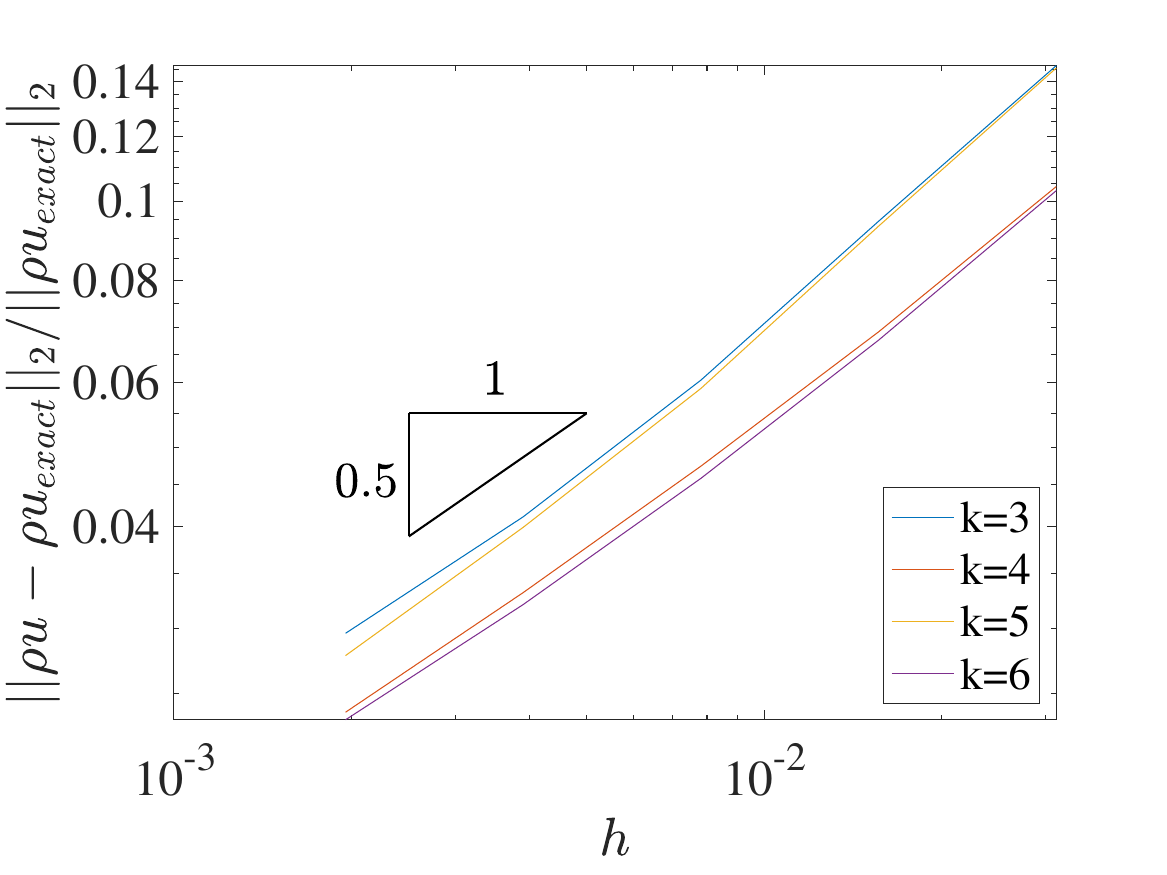}}\hfill
\subfloat[Momentum $L^1$ Error]{\label{sfig:sod_ru_l1_GP}\includegraphics[width=.5\textwidth]{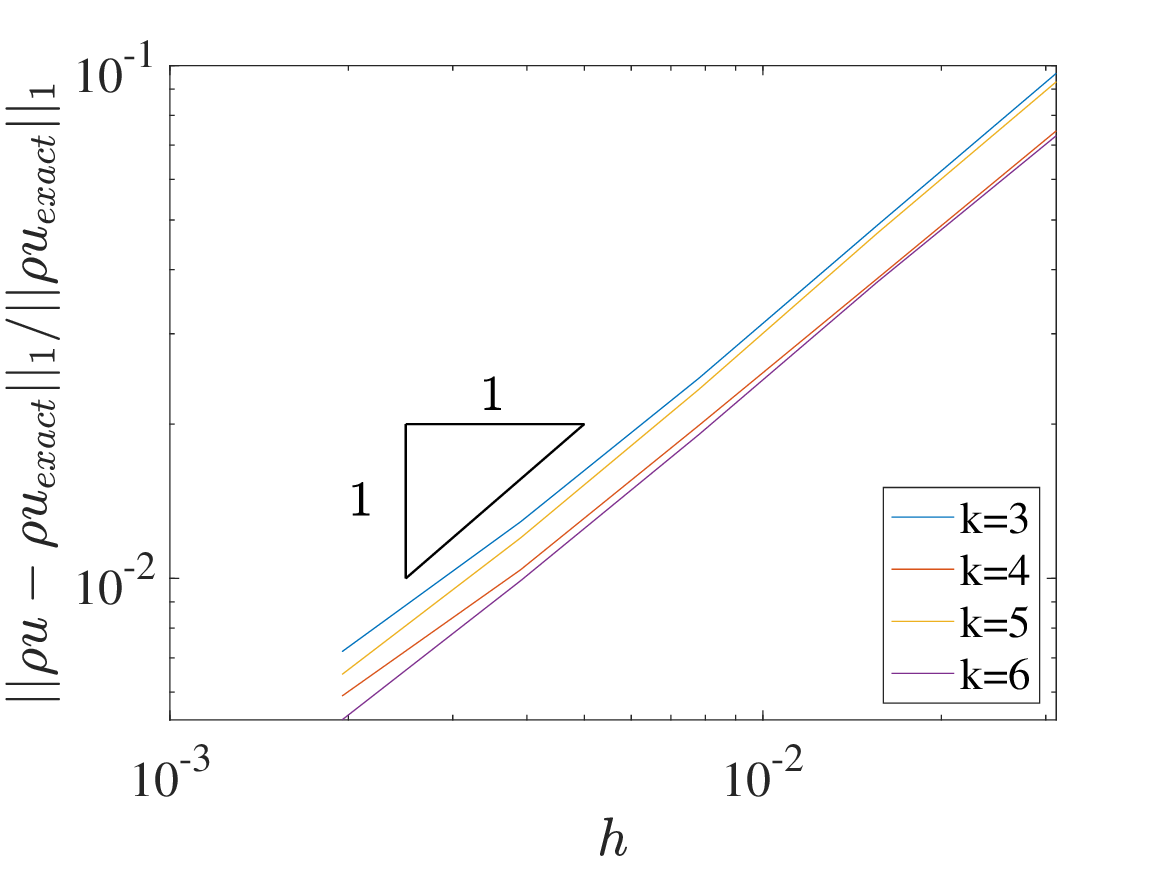}}\\
\subfloat[Energy $L^2$ Error]
{\label{sfig:sod_E_l2_GP}\includegraphics[width=.5\textwidth]{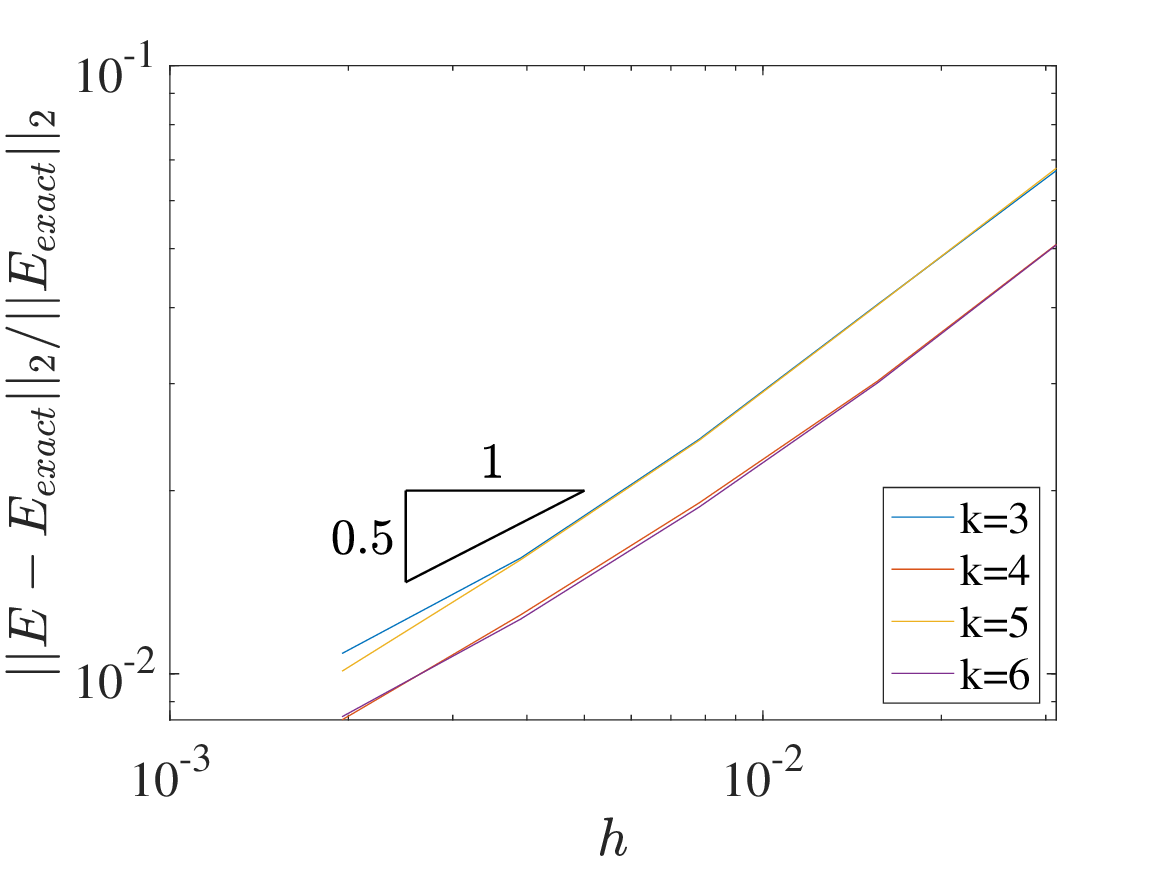}}\hfill
\subfloat[Energy $L^1$ Error]{\label{sfig:sod_E_l1_GP}\includegraphics[width=.5\textwidth]{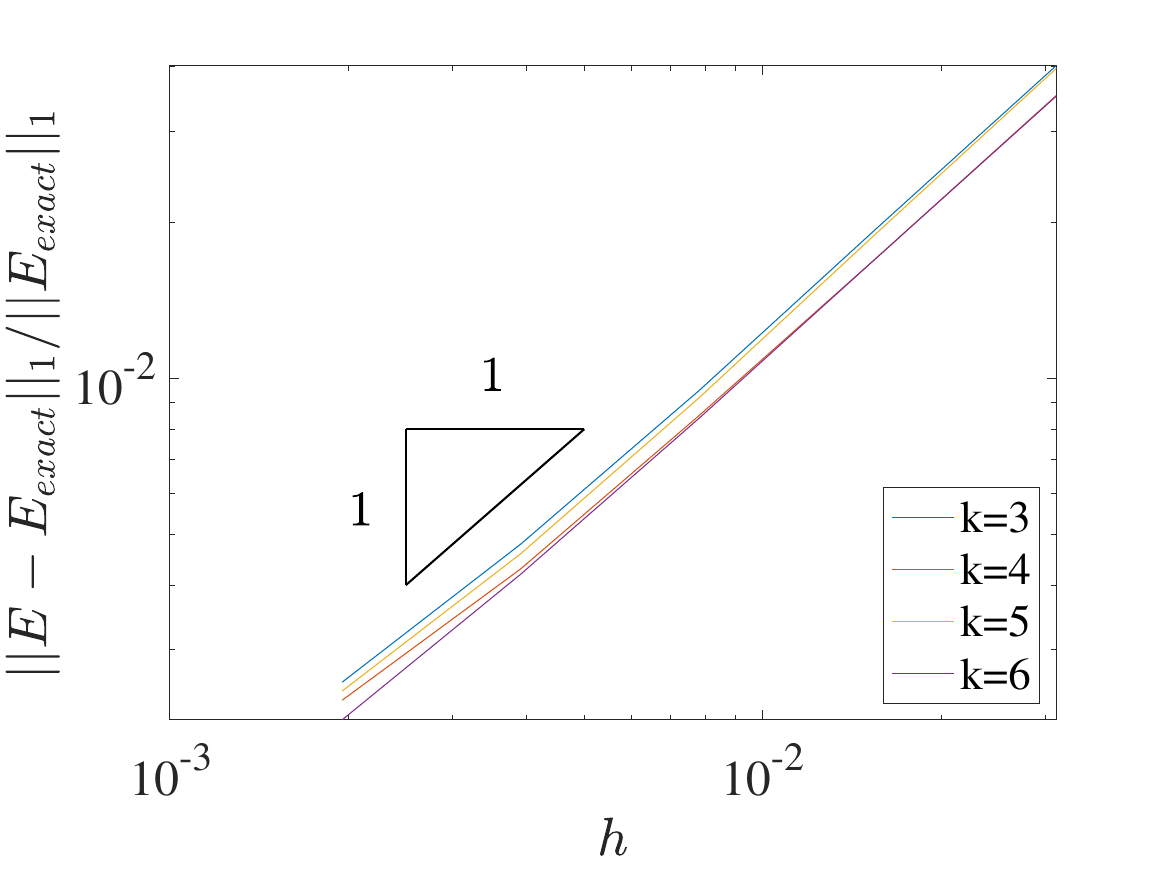}}\\
\caption{Convergence rates of Sod shock tube with Guermond-Popov fluxes}
\label{fig:sod_conv_GP}
\end{figure}

The $L^2$ and $L^1$ errors in the conserved quantities compared to the exact solution at the final time are shown in Figure \ref{fig:sod_conv_lap} for Laplacian flux regularization and Figure \ref{fig:sod_conv_GP} for Guermond-Popov regularization. In all cases rates close the optimal values of 0.5 and 1 are obtained. We also see lower errors for even values of $k$ in all of the cases, especially in the $L^2$ errors of the conserved quantities. Finally, the magnitude of the errors obtained match with those obtained in \cite{nazarov2017investigation} using finite elements stabilized using an entropy viscosity and the Guermond-Popov fluxes as well as the results obtained in \cite{tominec2023RBF} using RBF-FD methods stabilized with a residual-based viscosity and the Laplacian fluxes. 

\begin{figure}
\centering
\subfloat[Density]{\label{sfig:sod_r_elemComp_lap}\includegraphics[width=.5\textwidth]{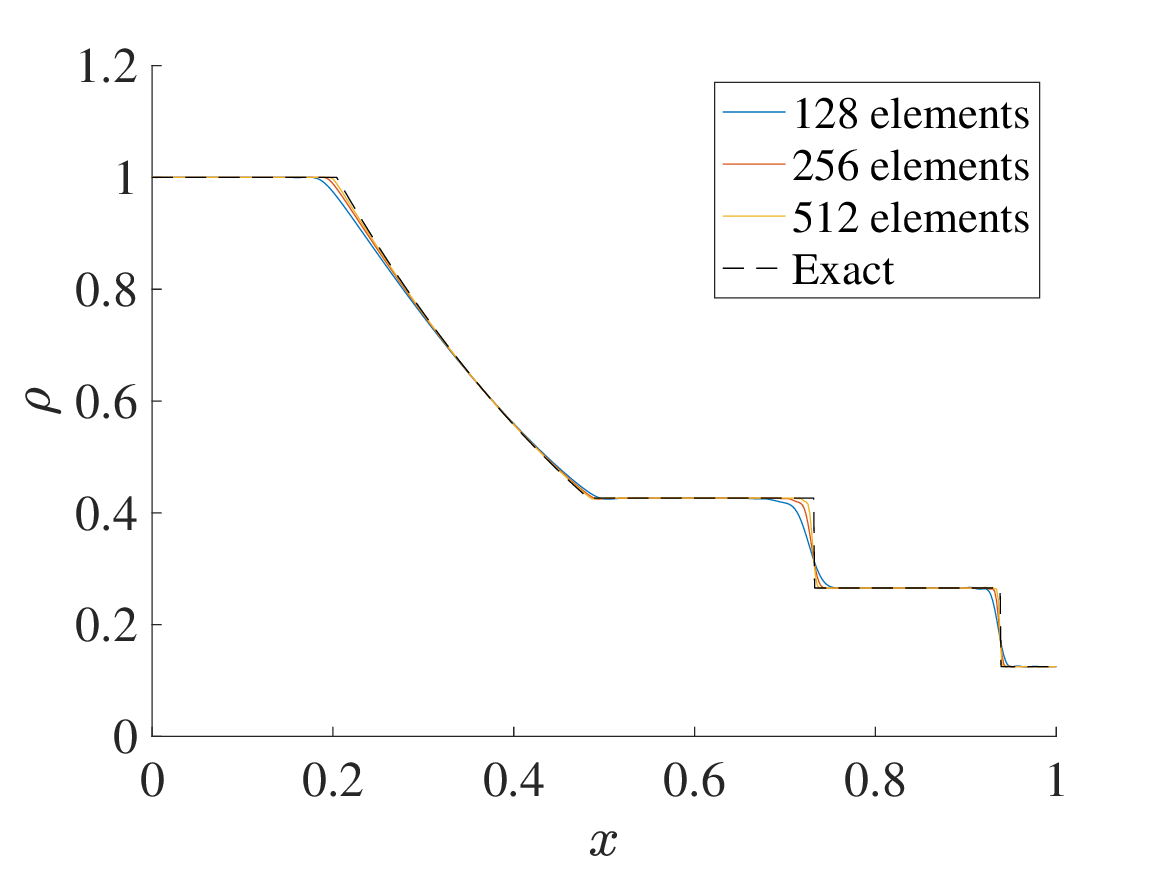}}\hfill
\subfloat[Velocity]{\label{sfig:sod_u_elemComp_lap}\includegraphics[width=.5\textwidth]{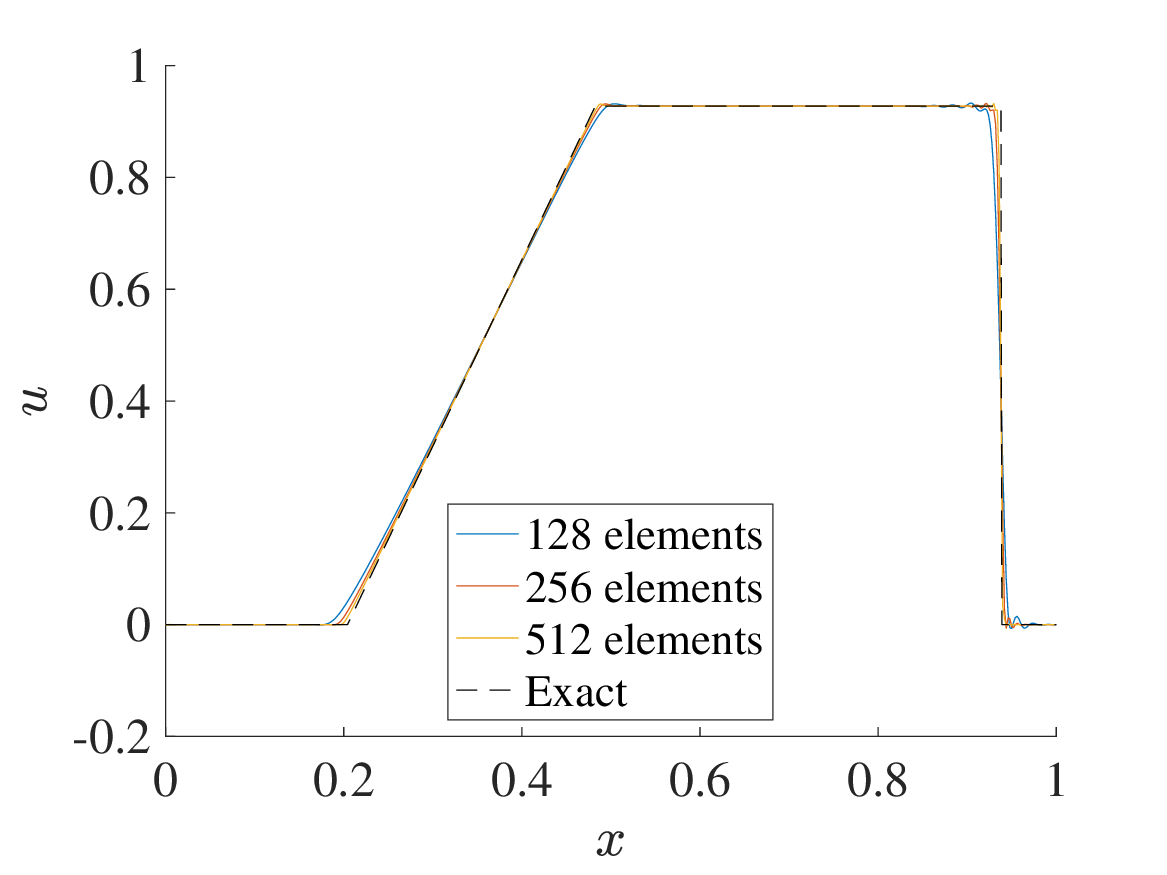}}\\
\subfloat[Pressure]{\label{sfig:sod_p_elemComp_lap}\includegraphics[width=.5\textwidth]{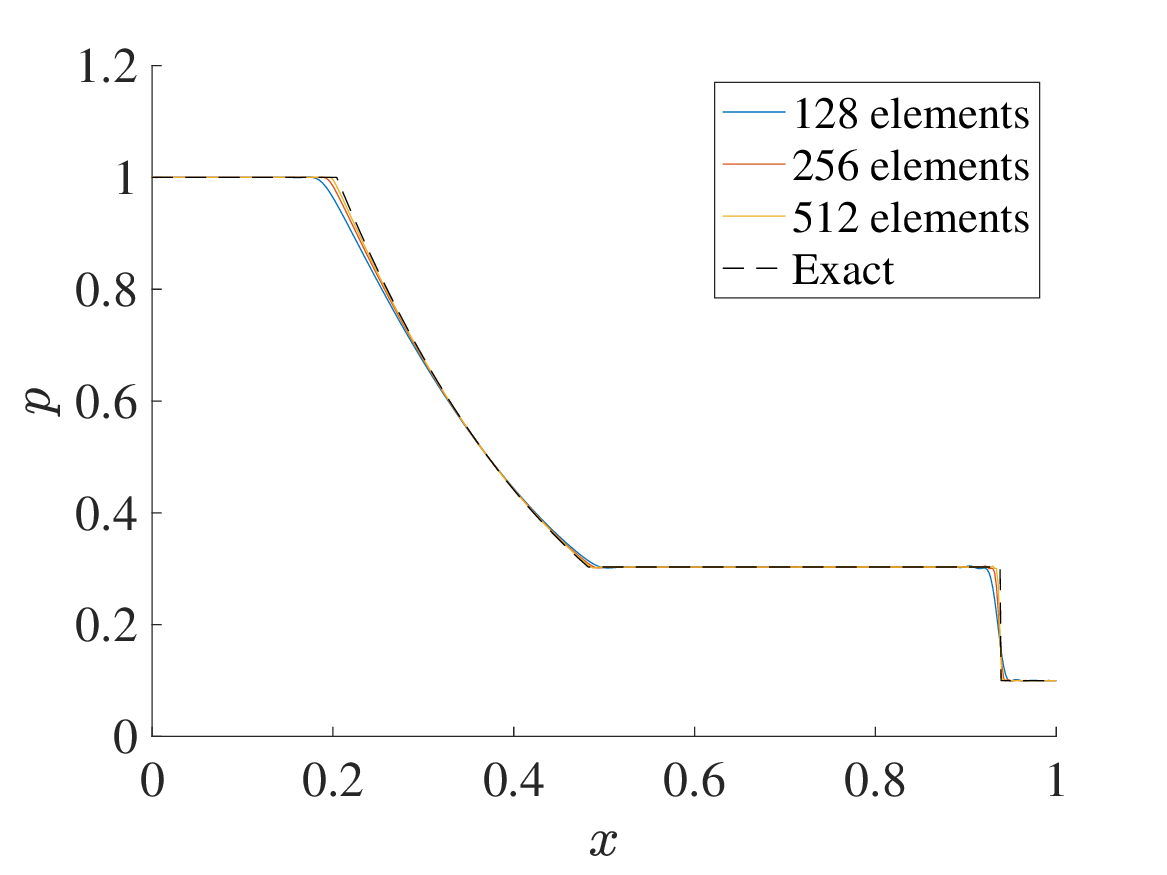}}\hfill
\subfloat[Residual-based viscosity]{\label{sfig:sod_nu_elemComp_lap}\includegraphics[width=.5\textwidth]{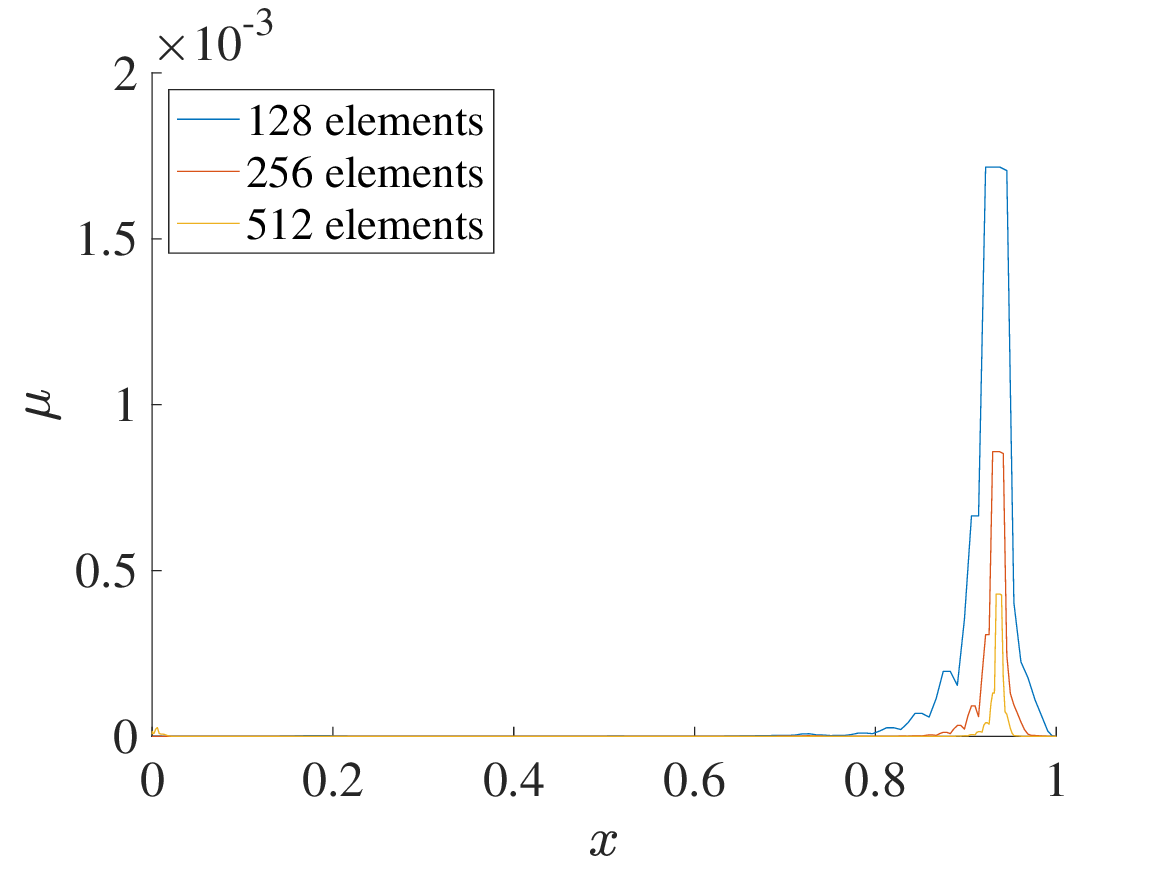}}\\
\caption{Sod problem solution with varying number of elements using Laplacian fluxes}
\label{fig:sod_elemComp_lap}
\end{figure}

\begin{figure}
\centering
\subfloat[Density]{\label{sfig:sod_r_elemComp_GP}\includegraphics[width=.5\textwidth]{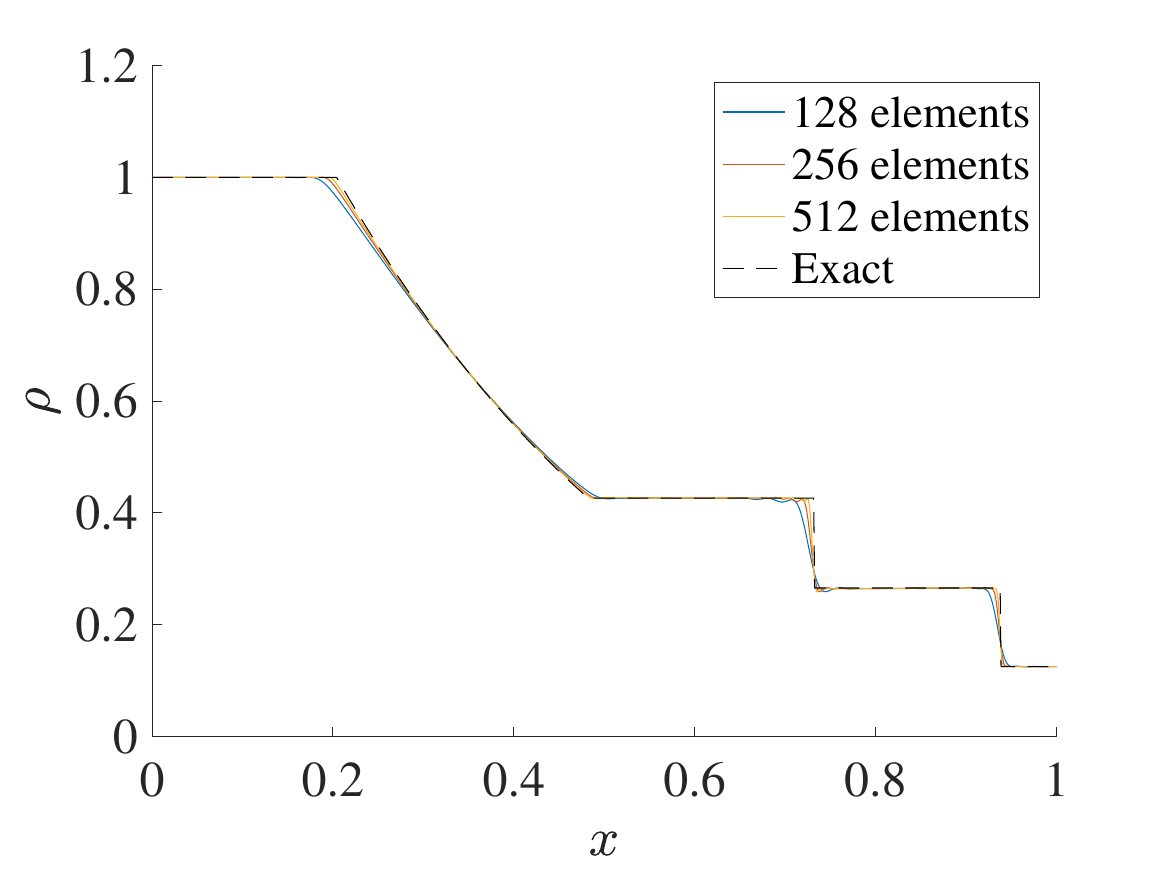}}\hfill
\subfloat[Velocity]{\label{sfig:sod_u_elemComp_GP}\includegraphics[width=.5\textwidth]{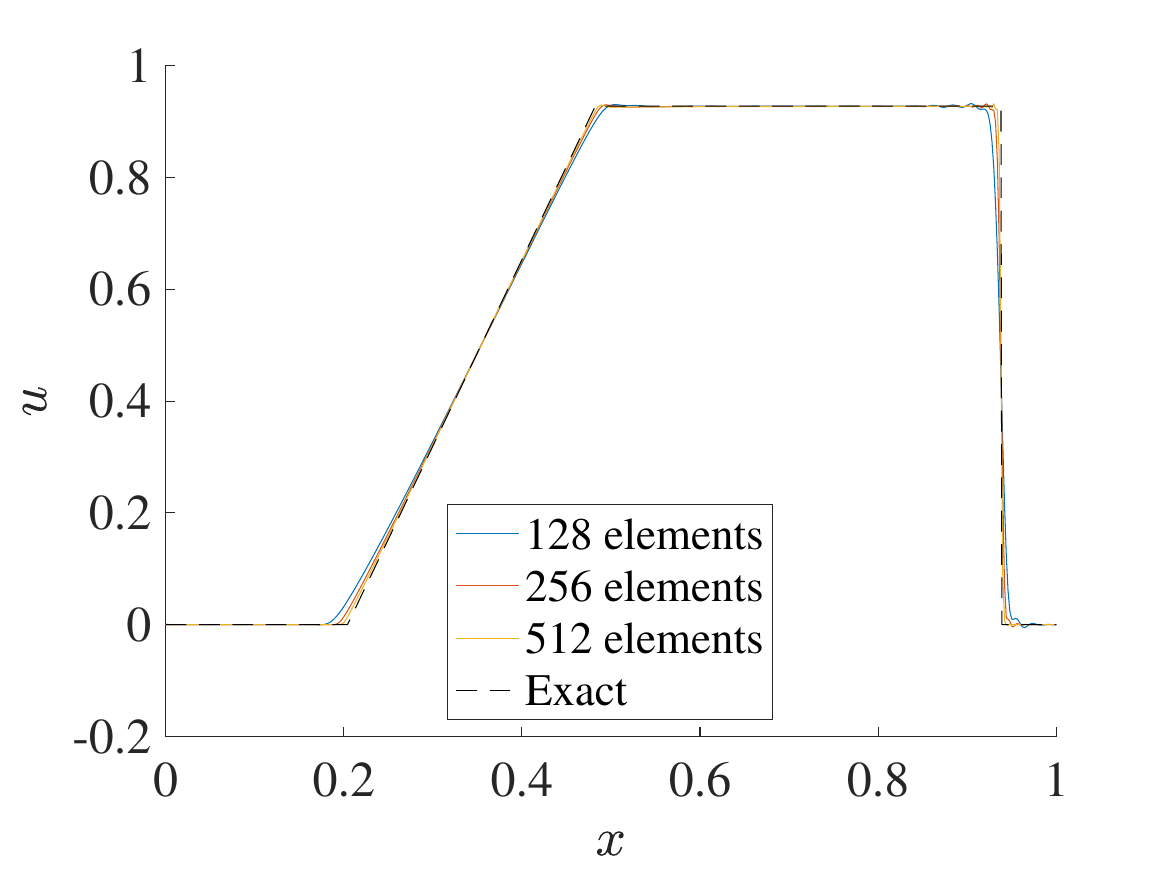}}\\
\subfloat[Pressure]{\label{sfig:sod_p_elemComp_GP}\includegraphics[width=.5\textwidth]{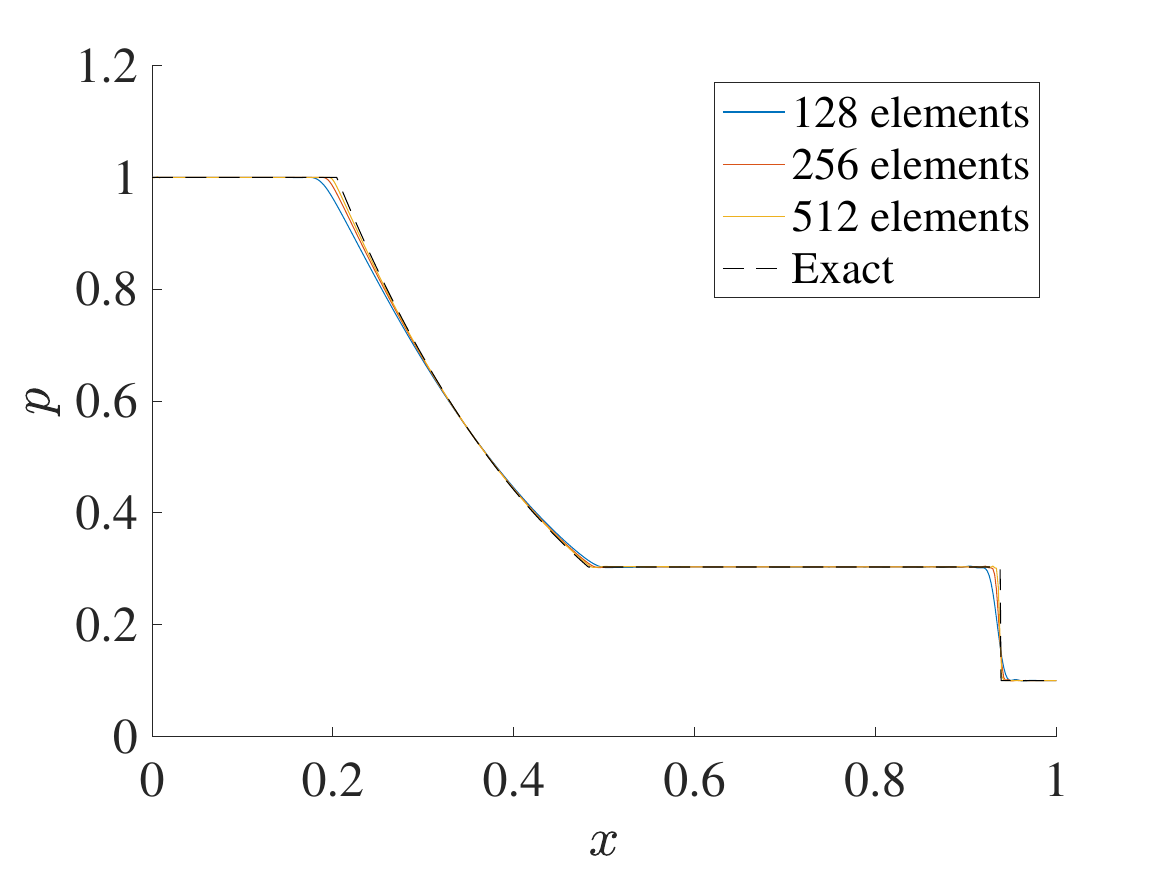}}\hfill
\subfloat[Residual-based viscosity]{\label{sfig:sod_nu_elemComp_GP}\includegraphics[width=.5\textwidth]{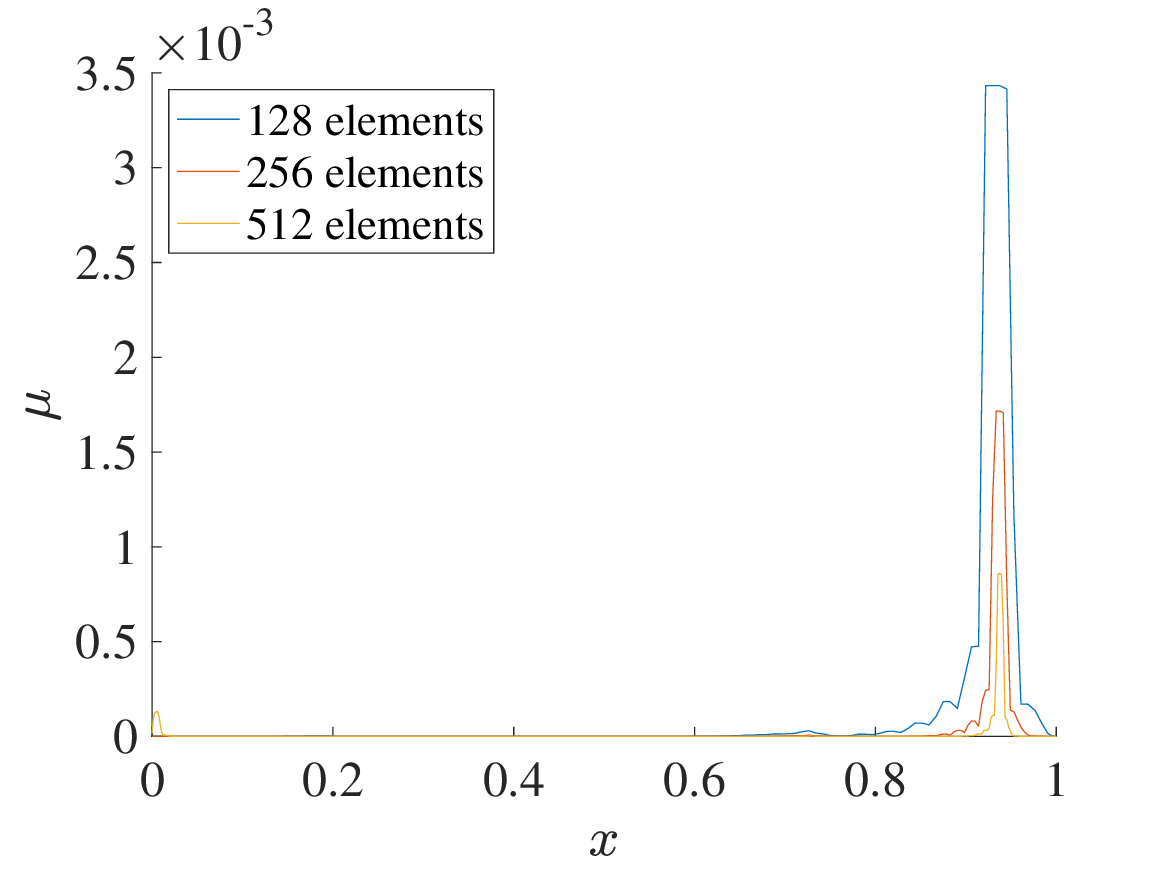}}\\
\caption{Sod problem solution with varying number of elements using Guermond-Popov fluxes}
\label{fig:sod_elemComp_GP}
\end{figure}

For a more qualitative comparison, Figure \ref{fig:sod_elemComp_lap} shows the density, velocity, pressure, and residual-based viscosity solutions for $k = 5$ and varying numbers of elements using the Laplacian fluxes, while Figure \ref{fig:sod_elemComp_GP} details the same results obtained using the Guermond-Popov regularization. The results match the exact solutions well in both cases, with the Laplacian flux results perhaps being slightly more diffuse near the contact. We also note that the residual-based viscosity is not active at the contact discontinuity, only at the shock location. 

\subsubsection{Shu-Osher Shock Tube}

Next we consider the Shu-Osher shock tube, which models the interaction of a Mach 3 shock with a sinusoidal density field \cite{shu1989efficient}. The computational domain is $x \in [0,10]$ and the initial condition is 

\begin{equation}
    \rho(x, 0) = 
    \begin{cases}
    3.857 & \textup{if } x < 1\\
    1+0.2\sin(5x) & \textup{Otherwise},
    \end{cases}
\end{equation}

\begin{equation}
     u(x, 0) = 
    \begin{cases}
    2.629 & \textup{if } x < 1\\
    0 & \textup{Otherwise},
    \end{cases}
\end{equation}

\begin{equation}
    p(x, 0) = 
    \begin{cases}
    10.333 & \textup{if } x < 1\\
    1& \textup{Otherwise}.
    \end{cases}
\end{equation}

\noindent The solution to this problem is characterized by the creation of high frequency density waves as the initial waves interact with the shock. This makes it a challenging test problem for numerical schemes, as the scheme must be dissipative enough to capture the shock without adding so much dissipation that the structure of the solution is destroyed. We simulate using a time step of $2 \times 10^{-5}$ to a final time of $t_f = 1.8$.

\begin{figure}
\centering
\subfloat[Density]{\label{sfig:Euler_SO_soln_lap}\includegraphics[width=.5\textwidth]{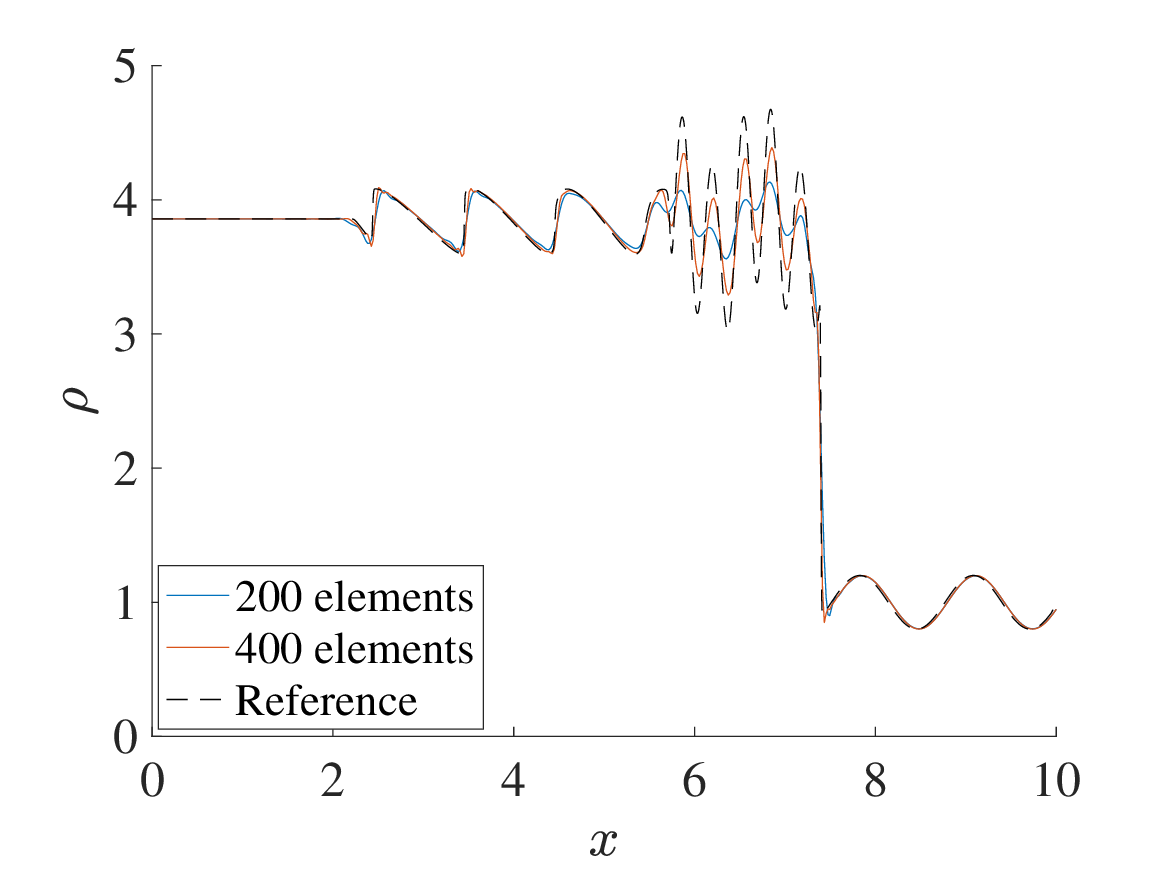}}\hfill
\subfloat[Zoom view]{\label{sfig:Euler_SO_soln_lap_zoom}\includegraphics[width=.5\textwidth]{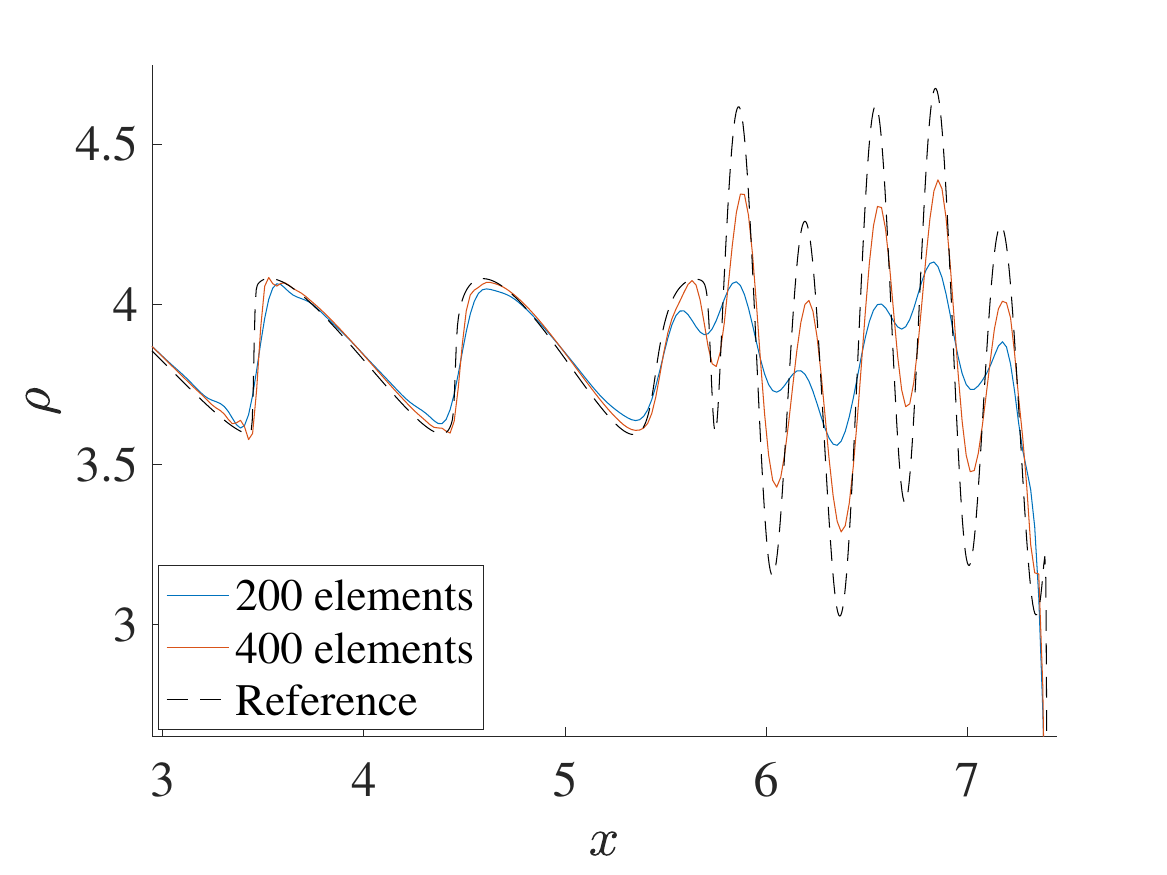}}\\
\caption{Shu-Osher problem solutions with varying number of elements using Laplacian fluxes}
\label{fig:Euler_SO_lap}
\end{figure}

\begin{figure}
\centering
\subfloat[Density]{\label{sfig:Euler_SO_soln_GP}\includegraphics[width=.5\textwidth]{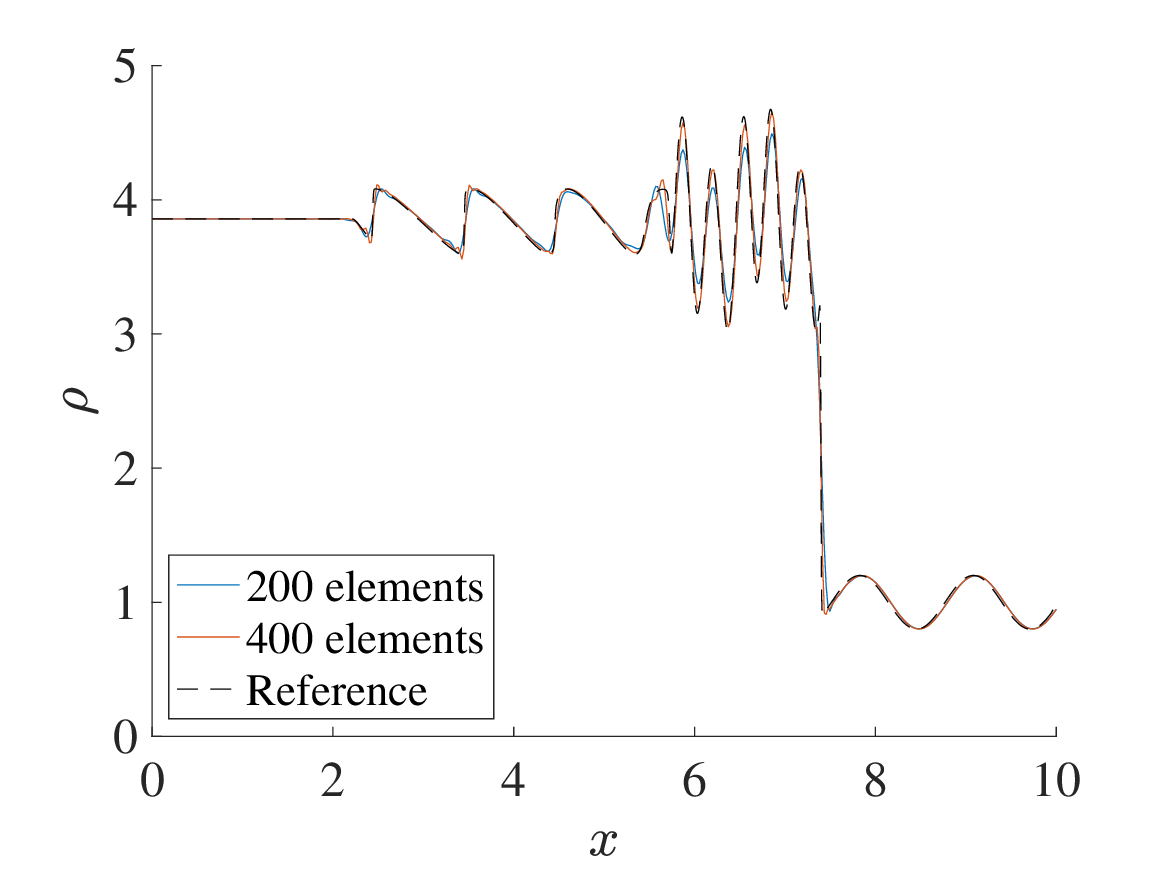}}\hfill
\subfloat[Zoom view]{\label{sfig:Euler_SO_soln_GP_zoom}\includegraphics[width=.5\textwidth]{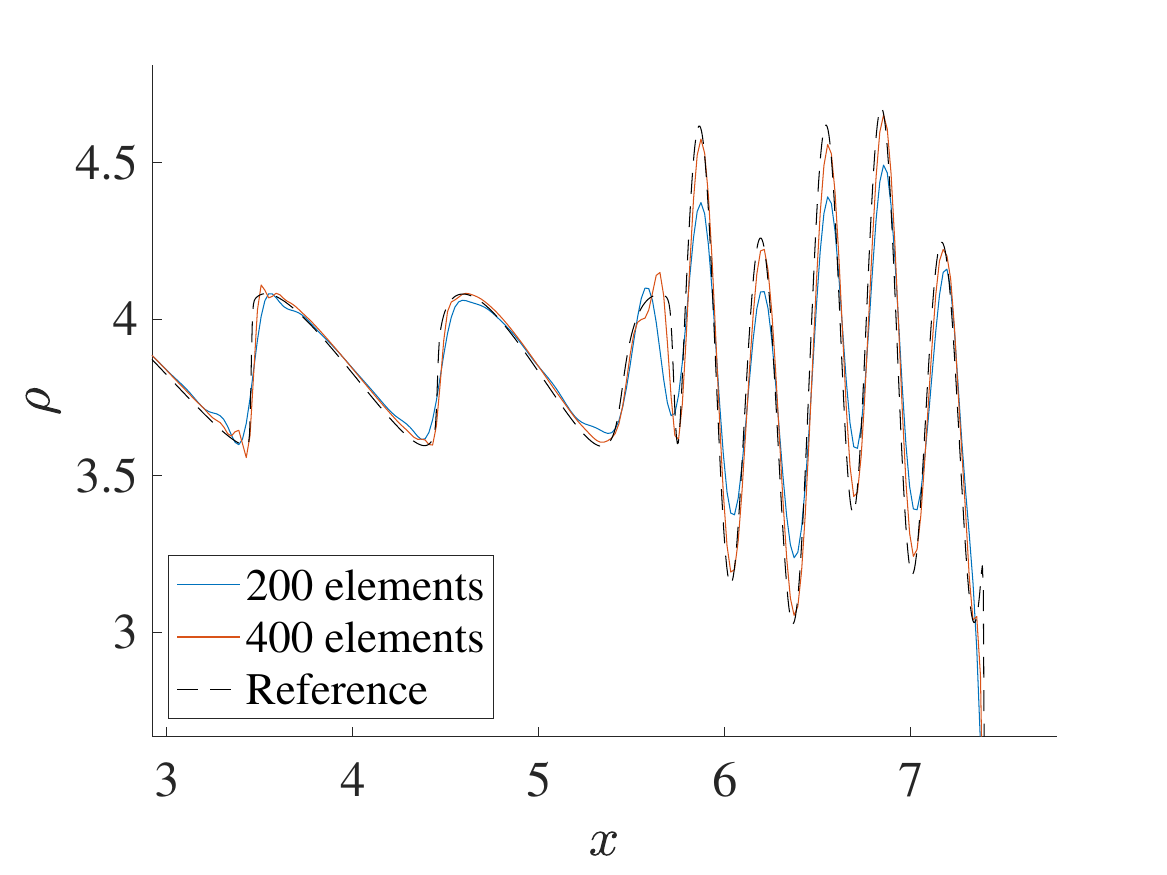}}\\
\caption{Shu-Osher problem solutions with varying number of elements using Guermond-Popov fluxes}
\label{fig:Euler_SO_GP}
\end{figure}

Figures \ref{fig:Euler_SO_lap} and \ref{fig:Euler_SO_GP} show the density fields obtained using each of the regularizations, with $k = 4$ and 2 different resolutions, as well as a reference solution obtained using a fifth order WENO-JS scheme on a 2000 cell grid. In this case there are much more pronounced differences between the Laplacian and Guermond-Popov results. The density field obtained using the Laplacian fluxes is much more damped than the Guermond-Popov results at the same resolutions. The Laplacian results look similar to the standard WENO results in \cite{fu2016teno}, which is known to be very diffusive. The Laplacian results can be sightly improved by weakening the nonlinear stabilization, however we have found that the improvement is minimal. Meanwhile, weakening the stabilization can reveal strong over and undershoots at the shock location.

The high frequency density oscillations in the 400 element, Guermond-Popov solution are almost identical to the reference solution. In this region the 200 element solution performs about as well as state of the art TENO schemes designed to minimize numerical disspation in cases such as this one with the same number of elements \cite{fu2016teno, fu2019tenoShort}. We do note, however, that there are small spurious oscillations which develop near the shocklets further upstream. In these regions the residual-based viscosity is still small, and suggests that perhaps further improvements to this stabilization mechanism may be warranted. 

\begin{figure}
\centering
\subfloat[Density]{\label{sfig:Euler_SO_soln_lin}\includegraphics[width=.5\textwidth]{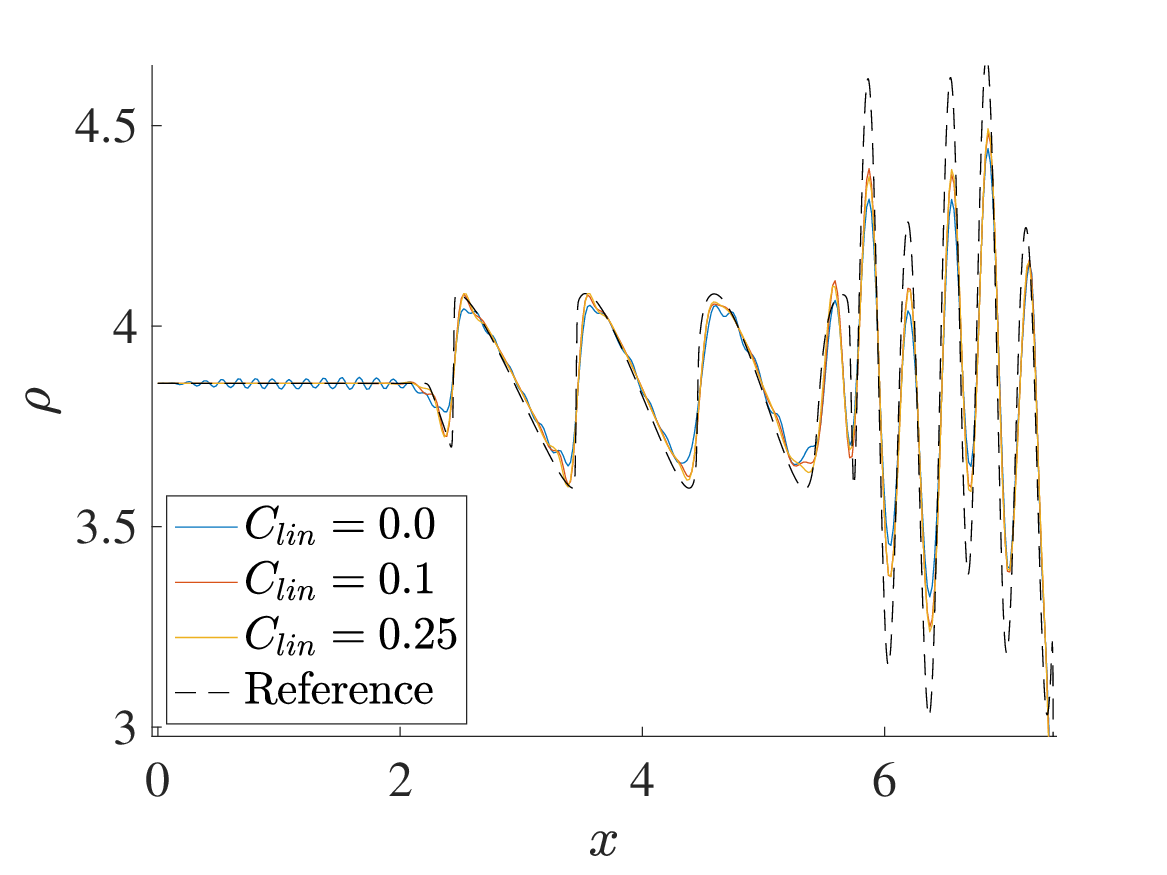}}\hfill
\subfloat[Viscosity]{\label{sfig:Euler_SO_nu_lin}\includegraphics[width=.5\textwidth]{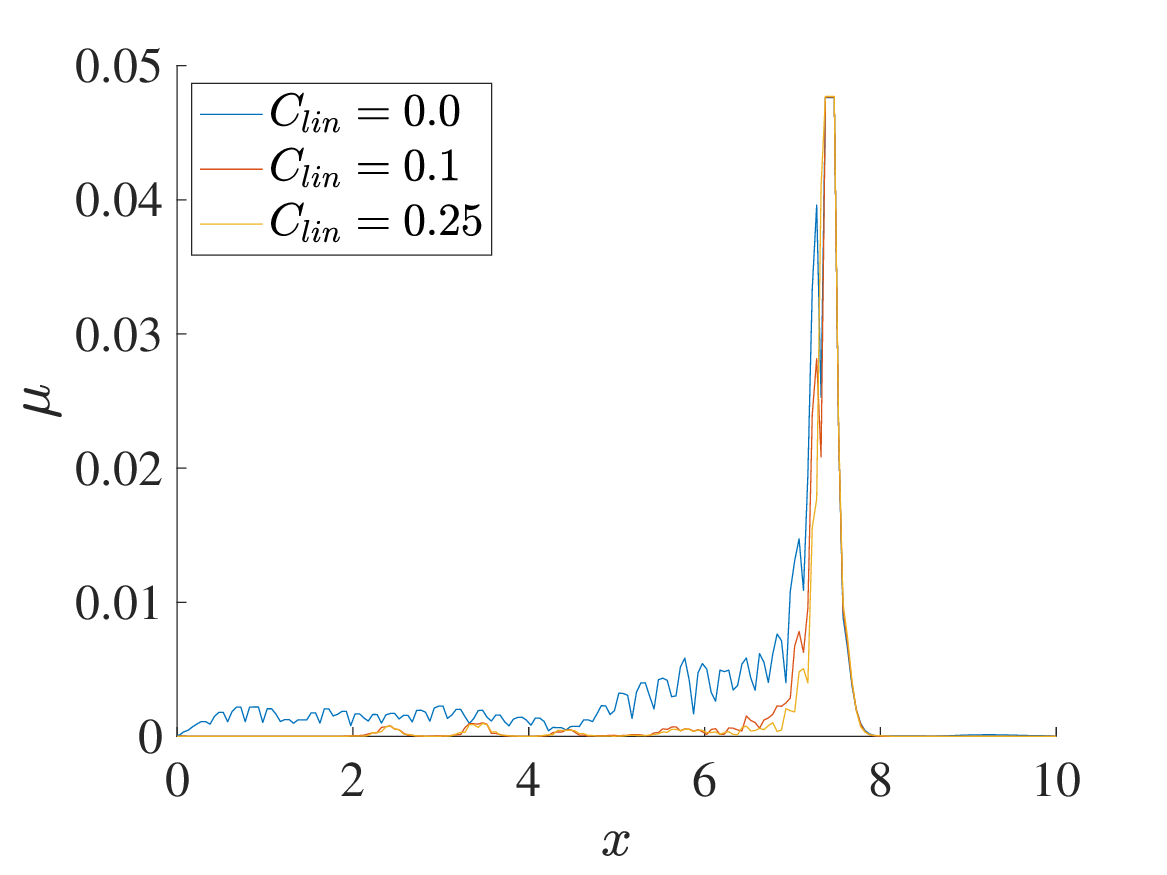}}\\
\caption{Shu-Osher problem solutions with varying linear stabilization strength and Guermond-Popov fluxes}
\label{fig:Euler_SO_lin}
\end{figure}

This test case also produces interesting results about the role of the linear stabilization in the collocation scheme. Figure \ref{fig:Euler_SO_lin} shows the zoomed-in density field as well as the full residual-based viscosity field for the Guermond-Popov scheme with a variety of linear stabilization strengths, including when it is removed completely. We can see the high frequency oscillations superimposed onto the density solution upwind of the shock, especially where the solution should be constant. We also see the residual-based viscosity activated over most of the domain, though it is ineffective at removing these oscillations, like in the linear advection case. Because of this, in the region of large density oscillations just after the shock we see that the results obtained without linear stabilization are actually slightly more dissipative (in terms of the larger structures) than results obtained with linear stabilization. Finally, we note that the results obtained with $C_{lin} = 0.1$ and $C_{lin} = 0.25$ are almost identical to one another.

\subsubsection{2D Riemann Problem: Case 12}

As a final test case, we consider the 2D Riemann problem titled Case 12 from \cite{liska2003comparison}. The domain of interest is $\mathbf{x}\in [0,1]^2$, with initial condition given by

\begin{equation}
   \rho(\mathbf{x}, 0) =  
   \begin{cases}
    4/5 & \textup{if } x < 1/2 \textup{ and } y < 1/2 \\
    1 & \textup{if } x < 1/2 \textup{ and } y \geq 1/2 \\
    1 & \textup{if } x \geq 1/2 \textup{ and } y < 1/2 \\
    17/32 & \textup{Otherwise}, 
    \end{cases}
\end{equation}

\begin{equation}
    \mathbf{u}(\mathbf{x}, 0) =  
    \begin{cases}
    (0,0) & \textup{if } x < 1/2 \textit{ and } y < 1/2 \\
    (3/\sqrt{17},0) & \textup{if } x < 1/2 \textup{ and } y \geq 1/2 \\
    (0,3/\sqrt{17}) & \textup{if } x \geq 1/2 \textup{ and } y < 1/2 \\
    (0,0) & \textup{Otherwise}, 
    \end{cases}
\end{equation}

\begin{equation}
    p(\mathbf{x}, 0) = 
    \begin{cases}
    1 & \textup{if } x < 1/2 \textup{ and } y < 1/2 \\
    1 & \textup{if } x < 1/2 \textup{ and } y \geq 1/2 \\
    1 & \textup{if } x \geq 1/2 \textup{ and } y < 1/2 \\
    2/5 & \textup{Otherwise}.
    \end{cases}
\end{equation}

\noindent Like the 2D Burgers example, we extend the domain to $\mathbf{x} \in [0,2]^2$, reflect the initial conditions such that the problem is symmetric, and use periodic boundary conditions. The simulation is advanced using a time step of $2 \times 10^{-4}$ to a final time of $t_f = 0.25$ so that we may compare the final solution against a number of published results. 

\begin{figure}
\centering
\subfloat[Solution]{\label{sfig:Euler_case12_soln_lap}\includegraphics[width=.5\textwidth]{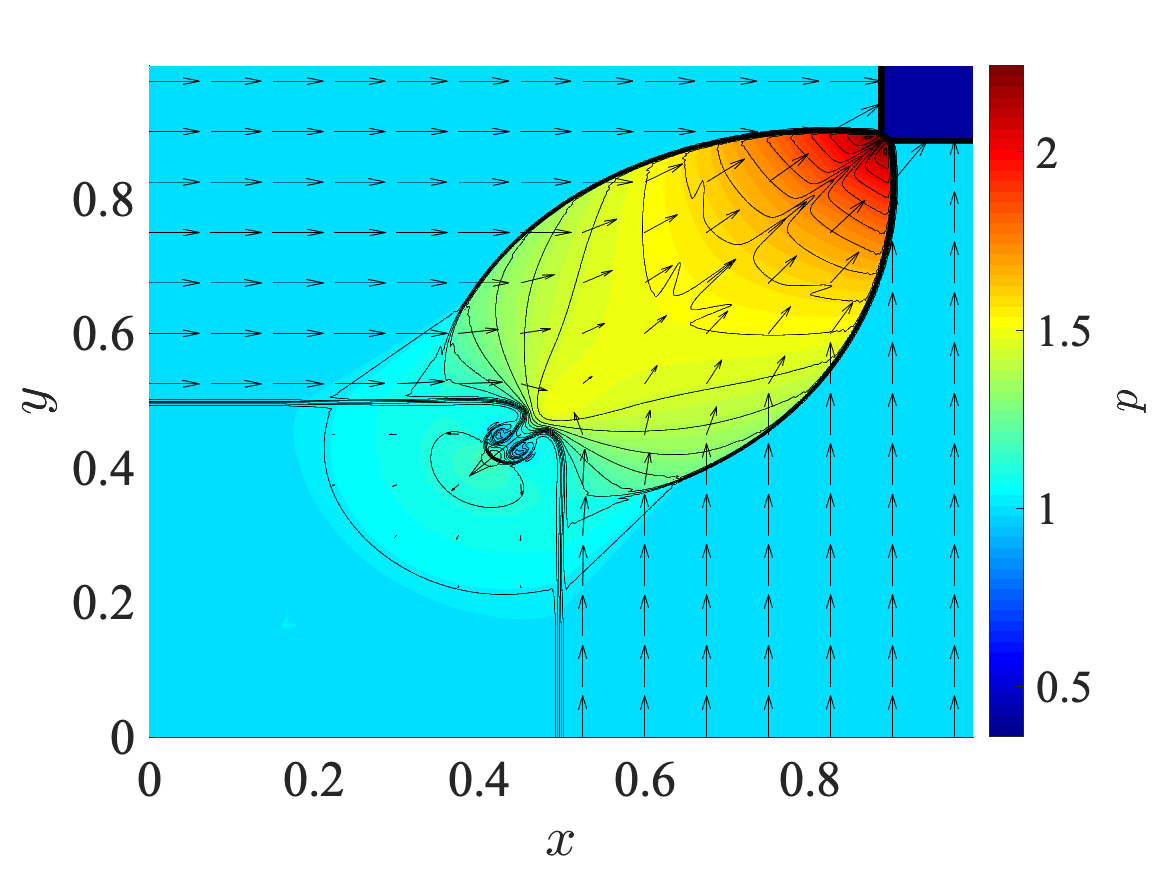}}\hfill
\subfloat[Residual-based viscosity]{\label{sfig:Euler_case12_nu_lap}\includegraphics[width=.5\textwidth]{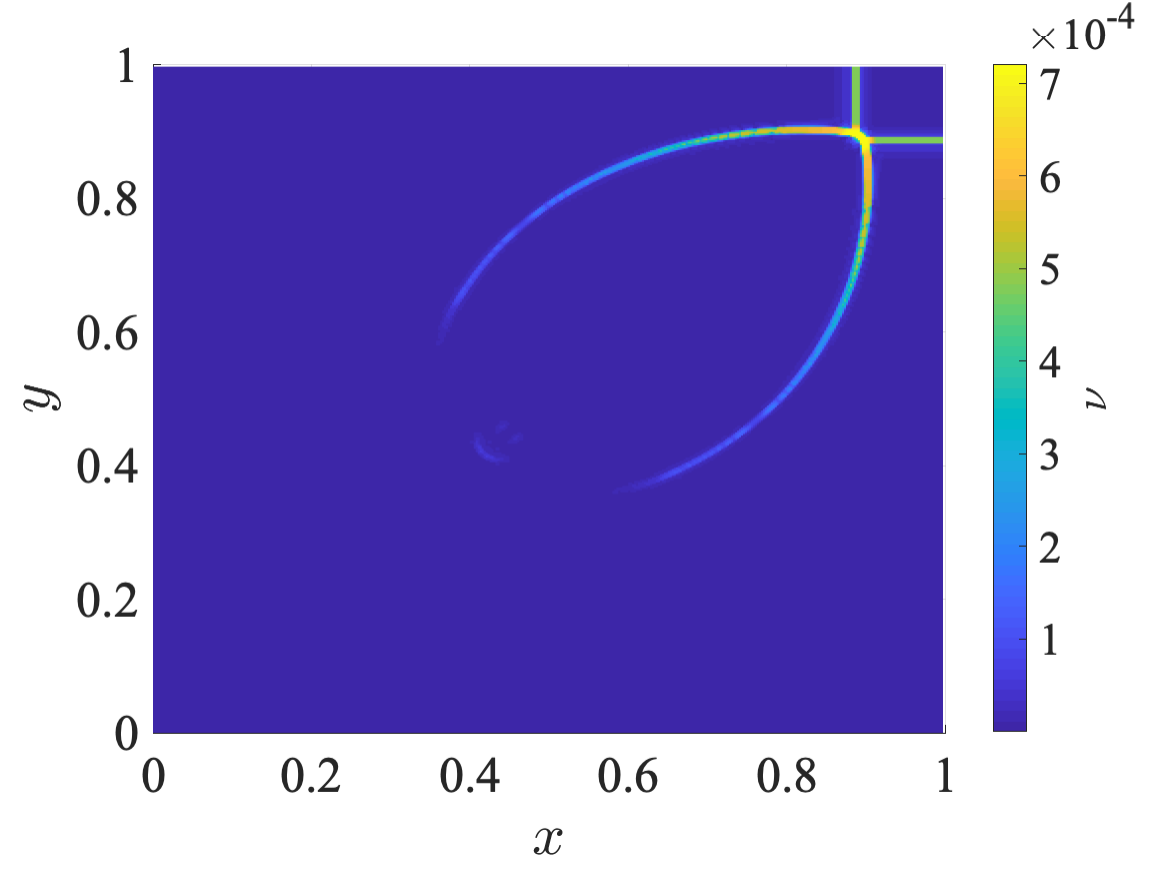}}\\
\caption{Case 12 solution and residual-based viscosity using Laplacian fluxes. Color of solution field is pressure, contours are density, and arrows are velocity.}
\label{fig:Euler_case12_lap}
\end{figure}

\begin{figure}
\centering
\subfloat[Solution]{\label{sfig:Euler_case12_soln_GP}\includegraphics[width=.5\textwidth]{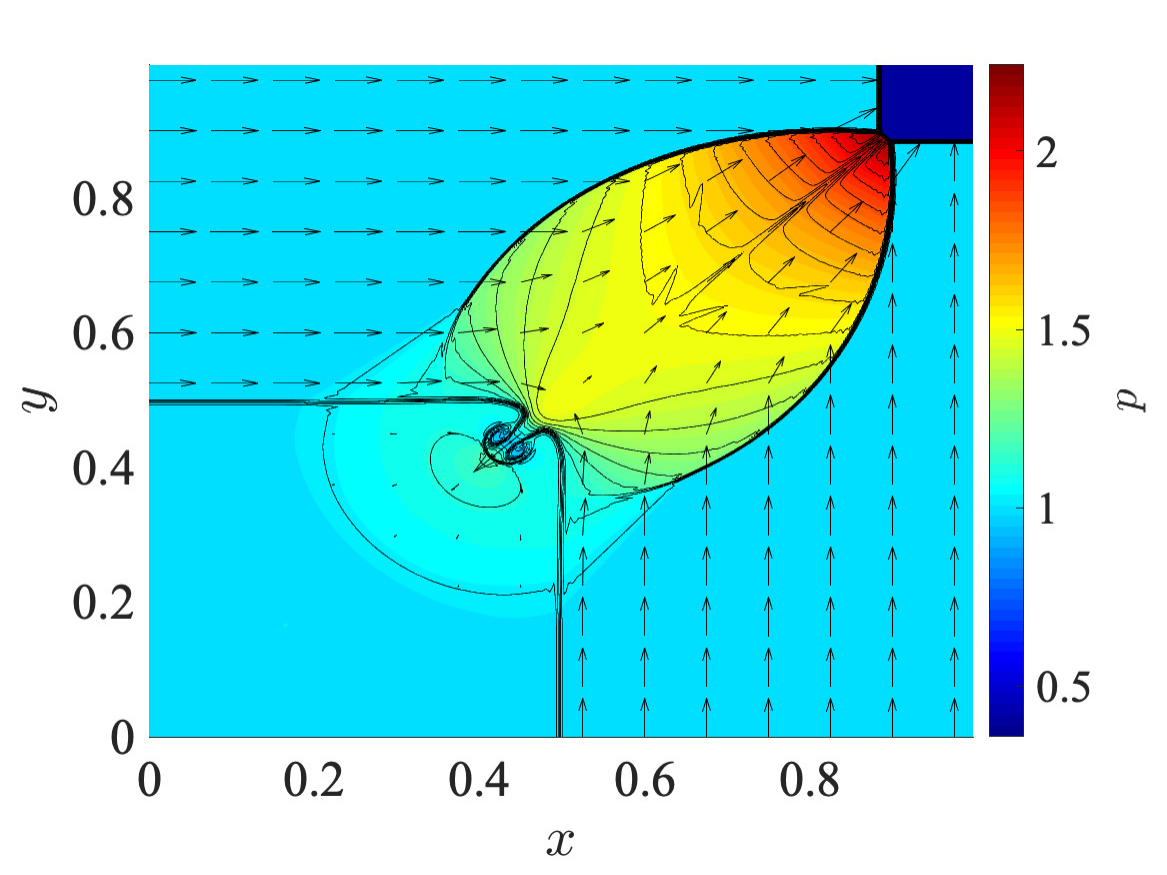}}\hfill
\subfloat[Residual-based viscosity]{\label{sfig:Euler_case12_nu_GP}\includegraphics[width=.5\textwidth]{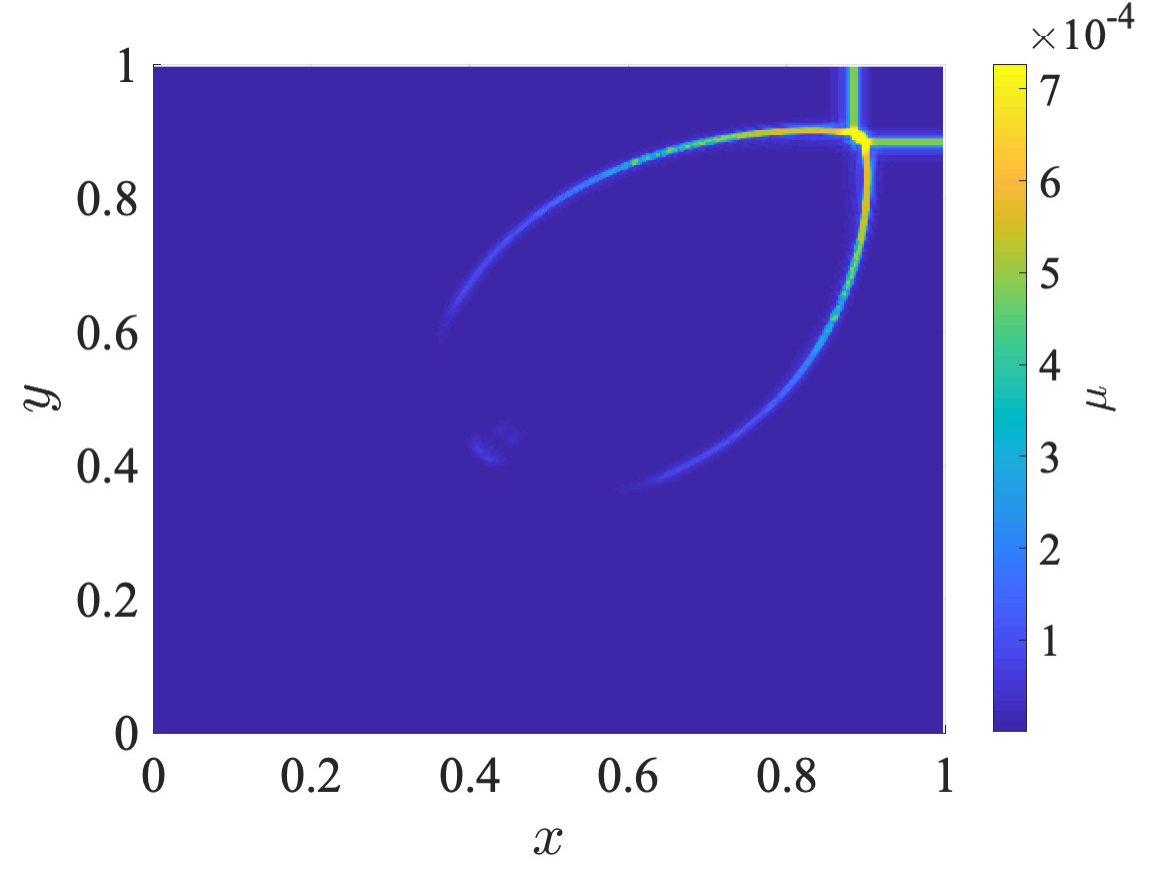}}\\
\caption{Case 12 solution and residual-based viscosity using Guermond-Popov fluxes. Color of solution field is pressure, contours are density, and arrows are velocity.}
\label{fig:Euler_case12_GP}
\end{figure}

The left of Figures \ref{fig:Euler_case12_lap} and \ref{fig:Euler_case12_GP} show the solutions obtained with $400^2$ elements over the domain $\mathbf{x} \in [0,1]^2$ and $k = 5$, where the contours represent density, the arrows represent the velocity field, and the colormap represents the pressure, as in \cite{liska2003comparison}. Note that this yields roughly the same number of degrees of freedom as the finite difference simulations in \cite{liska2003comparison}, the finite element solution in \cite{nazarov2017investigation}, and the RBF-FD solution in \cite{tominec2023RBF}. The right of Figures \ref{fig:Euler_case12_lap} and \ref{fig:Euler_case12_GP} show the resulting residual-based viscosity fields. Both solutions agree very well with the finite difference solutions of \cite{liska2003comparison}. Moreover, the density solutions are much less noisy and oscillatory than the solutions presented in \cite{tominec2023RBF, nazarov2017investigation}, where the former was computed using an RBF-FD discretization and the Laplacian fluxes, while the latter was computed using standard finite elements and the Guermond-Popov regularization.

The density peaks obtained with the Guermond-Popov fluxes are slightly sharper than those from the Laplacian solution, and these match well with those presented in \cite{fu2019tenoShort} generated using TENO schemes on a refined mesh of $1024^2$ cells. Consistent with the notion that the Guermond-Popov scheme is less dissipative, we do see some oscillations in the density field, especially near the stationary contact lines. These are much weaker in the Laplacian results. Within the Guermond-Popov scheme, these oscillations could be removed by increasing the viscosity, or perhaps using a definition of viscosity which activates more strongly in contacts, like the one used in \cite{nazarov2017investigation}.

\section{Conclusions}

In this work we have shown a method to construct spline collocation methods suitable for the simulation of hyperbolic conservation laws. By adding a residual-based shock capturing scheme along with a linear stabilization inspired by projection stabilization in FEA, our schemes are robust in the presence of shocks while maintaining the high-order approximation power of spline-based methods in the absence of discontinuities. Due to the nature of collocation schemes, these methods also have the potential for extremely efficient explicit time integration, due to the lack of any required numerical spatial integrations. Results obtained for a variety of conservation laws show the promise of the method as a simulation tool for compressible CFD. 

We believe this work opens up many interesting avenues for future research. One such interesting topic is the assessment of the effectiveness of different types of nonlinear stabilization techniques. While we have presented preliminary techniques in this work which seem to be effective, we do not claim that this is the only way one could achieve stable collocation schemes. Indeed we have seen in this work that different regularizations of the Euler equations can lead to different results. Different sensors which drive the stabilizing viscosity field could also be constructed, for example a WENO reconstruction sensor could be used \cite{kuzmin2023dissipation}. In addition, adapting other schemes such as positivity preserving limiters \cite{zhang2010positivity} to the collocation setting is also a very attractive direction. We are also interested in the application of this method to the compressible Navier-Stokes equations, where the scheme's effectiveness for efficient, scale-resolving simulations of turbulence could be studied in more detail. The stabilized scheme may also be able to function as an implicit Large Eddy Simulation technique, similar to \cite{nazarov2013AdaptiveLES}. Finally, while we have focused on simple domains in this work, results on more complicated domains could be obtained using, for example, NURBS basis functions and the isoparametric concept, or immersed collocation techniques which have recently been developed \cite{torre2023immersed}. 


\bibliography{ref}

\end{document}